\documentclass{article}

\usepackage{arxiv/arxiv}

\usepackage[utf8]{inputenc} % allow utf-8 input
\usepackage[T1]{fontenc}    % use 8-bit T1 fonts
\usepackage{hyperref}       % hyperlinks
\usepackage{url}            % simple URL typesetting
\usepackage{booktabs}       % professional-quality tables
\usepackage{amsfonts}       % blackboard math symbols
\usepackage{nicefrac}       % compact symbols for 1/2, etc.
\usepackage{microtype}      % microtypography
\usepackage{graphicx}
\usepackage[numbers]{natbib}
\usepackage{doi}
\usepackage{bbm}
\usepackage{amsmath}
\usepackage{subfig}
\usepackage{orcidlink}
\usepackage{enumitem}
\newlength{\jeroenlen}

\newtheorem{remark}{Remark}
\newcommand{\zh}{\mathbf{z}_h}
\newcommand{\yh}{\mathbf{y}_h}
\newcommand{\yinit}{\mathbf{y}_0}
\newcommand{\tildeyh}{\tilde{\mathbf{y}}_h}
\newcommand{\yhrec}{\mathbf{y}_{h,\mathrm{rec}}}
\newcommand{\yn}{\mathbf{y}_N}
\newcommand{\tildeyn}{\tilde{\mathbf{y}}_N}
\newcommand{\uh}{\mathbf{u}_h}
\newcommand{\tildeuh}{\tilde{\mathbf{u}}_h}
\newcommand{\uhrec}{\mathbf{u}_{h,\mathrm{rec}}}
\newcommand{\un}{\mathbf{u}_N} 
\newcommand{\tildeun}{\tilde{\mathbf{u}}_N}
\newcommand{\mus}{\boldsymbol{\mu}_s}

\newcommand{\R}{\mathbb{R}}

\newcommand{\acapo}{}

\title{Latent feedback control of distributed systems in multiple scenarios through deep learning-based reduced order models}

\author{\href{https://orcid.org/0000-0002-7111-5070}{\includegraphics[scale=0.06]{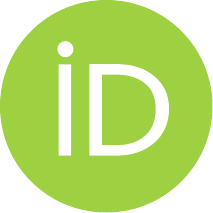}\hspace{1mm}Matteo Tomasetto} \\
	Department of Mechanical Engineering\\
	Politecnico di Milano, Milano, Italy \\
        \texttt{matteo.tomasetto@polimi.it}
	\And 
 \href{https://orcid.org/0000-0002-0476-4118}{\includegraphics[scale=0.06]{orcid.pdf}\hspace{1mm}Francesco Braghin} \\
	Department of Mechanical Engineering\\
	Politecnico di Milano, Milano, Italy \\ \texttt{francesco.braghin@polimi.it}
	\AND
	\href{https://orcid.org/0000-0001-8277-2802}{\includegraphics[scale=0.06]{orcid.pdf}\hspace{1mm}Andrea Manzoni} \\
        MOX - Department of Mathematics \\
        Politecnico di Milano, Milano, Italy \\
        \texttt{andrea1.manzoni@polimi.it}
        }

\date{}

\renewcommand{\headeright}{}
\renewcommand{\undertitle}{}
\renewcommand{\shorttitle}{Latent feedback control of distributed systems in multiple scenarios through DL-ROMs}

\hypersetup{
pdftitle={Latent feedback control of distributed systems in multiple scenarios through DL-ROMs},
pdfauthor={Matteo Tomasetto, Francesco Braghin, Andrea Manzoni}
}

\begin{document}
\maketitle

\begin{abstract}
Continuous monitoring and real-time control of high-dimensional distributed systems are often crucial in applications to ensure a desired physical behavior, without degrading stability and system performances. Traditional feedback control design that relies on full-order models, such as high-dimensional state-space representations or partial differential equations, fails to meet these requirements due to the delay in the control computation, which requires multiple expensive simulations of the physical system. The computational bottleneck is even more severe when considering parametrized systems, as new strategies have to be determined for every new scenario. To address these challenges, we propose a real-time closed-loop control strategy enhanced by nonlinear non-intrusive Deep Learning-based Reduced Order Models (DL-ROMs). Specifically, in the offline phase, {\em (i)} full-order state-control pairs are generated for different scenarios through the adjoint method, {\em (ii)} the essential features relevant for control design are extracted from the snapshots through a combination of Proper Orthogonal Decomposition (POD) and deep autoencoders, and {\em (iii)} the low-dimensional policy bridging latent control and state spaces is approximated with a feedforward neural network. After data generation and neural networks training, the optimal control actions are retrieved in real-time for any observed state and scenario. In addition, the dynamics may be approximated through a cheap surrogate model in order to close the loop at the latent level, thus continuously controlling the system in real-time even when full-order state measurements are missing. The effectiveness of the proposed method, in terms of computational speed, accuracy, and robustness against noisy data, is finally assessed on two different high-dimensional optimal transport problems, one of which also involving an underlying fluid flow.
\end{abstract}

% keywords can be removed

\section{Introduction}\label{sec:introduction}

Monitoring and controlling dynamical systems play a central role in Applied Sciences and Engineering. A control strategy enables to steer and influence the behaviour of a physical phenomenon in order to achieve desired outcomes. Fascinating examples in this direction are given by self-driving cars navigating traffic, spacecraft maneuvering through space, and emission reduction of industrial processes. Distributed dynamical systems are often described by high-dimensional state-space representations or Partial Differential Equations (PDEs). The space-time evolution of the state variable may be retrieved through full-order models, such as, e.g., finite differences, finite elements or finite volumes, yielding to a (possibly nonlinear) system of equations to be solved. This task is computationally demanding, especially when dealing with complex and stiff phenomena, since a large number of degrees of freedom are required to discretize and capture the state dynamics. The computational bottleneck becomes even more severe when {\em (i)} considering systems parametrized by a vector of scenario parameters $\mus$, as different independent resolutions are required for every scenario of interest, and when {\em (ii)} dealing with Optimal Control Problems (OCPs). Indeed, optimal control values are typically determined through multiple full-order resolutions of the physical model within an optimization procedure minimizing a problem-specific \textit{loss} or \textit{cost functional} \cite{Manzoni2021}, making real-time control unfeasible. However, delayed control strategies may degrade the overall optimality, stability and system
performances, which are often vital in applications. Our main goal is to overcome the computational barrier of traditional full-order solvers by devising an adaptive and instantaneous feedback control strategy in multiple scenarios for high-dimensional parametrized PDEs.
\\

Feedback control design is traditionally addressed with Dynamic Programming Principle \citep{Bellman1957} and Hamilton-Jacobi-Bellman equation (HJB) \cite{Bardi2009}. Despite attempts to circumvent and reduce the problem complexity in this framework -- see, e.g., \cite{Kunish2004, Alla2016, Schmidt2018, Alla2020, Albi2024} -- numerical schemes suffer from the curse of dimensionality and become computationally intractable as the state and control dimensions increase, thus limiting its applicability in real-time control scenarios.
\\

Instead of directly solving the HJB equation, the feedback controller may be determined exploiting Reinforcement Learning \cite{Sutton2018} or Deep Reinforcement Learning (DRL) \cite{Garnier2021, Vignon2023a, Rabault2020} algorithms. Specifically, the so-called policy or agent -- that is modeled through a neural network in DRL -- learns to predict the best possible control action from the current observed state. In practice, policy training is carried out through a continuous interaction with the dynamical system, maximizing the reward signal computed by looking at the system response to the applied control. This strategy has been explored in several applications, ranging from flow control \cite{Verma2018, Bucci2019, Ren2021a, Varela2022, Guastoni2023, Vignon2023, Vinuesa2022}, metamaterial design \cite{Shah2021, Mirzakhanloo2020, Rosafalco2023}, chemistry \cite{Zhou2019}, optimal navigation in turbulent flows \cite{Tonti2024, Biferale202} and swarm systems control \cite{Huttenrauch2019}. However, standard DRL approaches are limited to low-dimensional state and control spaces due to sample inefficiency and computationally demanding training phases, especially when dealing with parametrized problems. To overcome these drawbacks, efficient, robust and interpretable policies have been proposed by \cite{Zolman2024, Botteghi2024} leveraging sparsity and dictionary learning. Instead, \cite{Ma2018, Garnier2021} and \cite{Botteghi2021} exploit, respectively, autoencoders and a priori knowledge to cope with high-dimensional observations. Similar state compressions are employed by \cite{Luo2023} in the context of robust control, and by \cite{Ishize2023} to provide a fast linear quadratic regulator feedback controller. Moreover, distributed actions may be considered through Multi-Agent Reinforcement Learning (MARL) strategies \cite{Busoniu2008, Vasanth2024, Peitz2024}, where the dynamics is controlled by means of multiple local agents relying solely on local state values. However, the locality assumption compromises a smooth and non-local coordination among agents, that is typically crucial when controlling physical systems \cite{Jeon2024}.
\\

Another advanced technique for closed-loop control of dynamical systems is Model Predictive Control (MPC) \cite{Camacho2004}. In the MPC context, multiple open-loop optimization problems are solved as the dynamics evolves in time, iteratively updating the information on the current state looking at the system evolution. Rather than full-order models describing the dynamics at hand, MPC usually considers (possibly nonlinear) fast-evaluable surrogate models to quickly predict the state evolution over a prediction horizon starting from the current observed state. This approximation speeds up the optimization procedure aiming at computing the optimal control sequence in time by minimizing a suitable cost functional. Several techniques have been exploited to identify (possibly low-dimensional) surrogate models for the system dynamics -- such as neural networks \cite{Draeger1995, Aggelogiannaki2008, Chen2018, Bieker2019, Peitz2023, Antonelo2024}, POD \cite{Ghiglieri2014, Alla2015}, Sparse Identification of Nonlinear Dynamics (SINDy) \cite{Kaiser2018}, Eigensystem Realization Algorithm (ERA) \cite{Hickner2023} -- or to approximate the Koopman operator -- such as Dynamic Mode Decomposition (DMD) \cite{Korda2018, Peitz2019, Peitz2020, Klus2020} -- in the MPC framework. For further details on DMD, SINDy and Koopman theory see, e.g., \cite{Schmid2010, Brunton2016, Brunton2016a, Brunton2019, Brunton2022}. However, MPC may be not suitable for real-time applications with strict timing requirements and fast dynamics due to the multiple optimization problems involved online. Moreover, since the control is regarded as model input, it would be challenging to consider distributed control variables with a remarkably high number of degrees of freedom.
\\

In this work, we design a real-time feedback control strategy capable of steering the dynamics of distributed systems described in terms of parametrized PDEs in multiple scenarios. Differently from the DRL approach, we focus on an offline-online decomposition, in the direction of imitation learning or offline reinforcement learning \cite{Levine2020}, in order to reduce the dimensionality of distributed and boundary control actions, as well as high-dimensional state observations. As far as the dimensionality reduction is concerned, we consider nonlinear non-intrusive Deep Learning-based Reduced Order Models (DL-ROMs), where data are compressed through Proper Orthogonal Decomposition (POD) \cite{Hesthaven2018}, deep autoencoders \cite{Fresca2021, Franco2023}, or a combination thereof \cite{Fresca2022}, while a feedforward neural network approximates low-dimensional policy bridging state and control latent spaces. Nonlinearity and non-intrusiveness allow for greater flexibility and speed-up with respect to traditional reduced order modeling (ROM) techniques, such as the Reduced Basis (RB) method \cite{Hesthaven2015, Quarteroni2015, Brunton2019, Manzoni2019}. Indeed, despite several applications in the context of OCPs \cite{Kunisch1999, Leibfritz2006, Kunisch2008, negri2013reduced, Benner2014, Amsallem2015, Sinigaglia2022-ROM}, the RB method is not efficient when dealing with nonlinear or transport-dominated problems. 
\\

The proposed approach extends our previous work on real-time open-loop control for parametrized PDEs \cite{Tomasetto2024}, where we approximate the optimal state and control variables starting from the corresponding scenario parameters $\mus$ only, disregarding possible additional state measurements collected online. Whenever the latter data are available, it is possible to continuously inform the controller about the current system behaviour through a feedback loop, thus broadening the control strategy to challenging applications where the scenario parameters do not describe the whole problem variability. For example, in the optimal transport test cases detailed in Section~\ref{sec:test}, the scenario parameters identify the target location, while the current configuration is captured directly from the observed state, as it typically happens in realistic settings. Moreover, as shown in the test cases taken into account, the dimensionality reduction and the feedback signal are also helpful to gain robustness against uncertainties and deal with noisy state data, paving the way for sensor-based applications.
\\

Inspired by the MPC strategy, we propose a model closure at the latent level in order to control the dynamical system even when state data are not available online due to, e.g., sensor failures, delay in receiving the measurements or time-consuming simulations. Specifically, in our framework, the latent dynamics is approximated through a low-dimensional deep learning-based surrogate model, which quickly predicts the reduced state evolution starting from the current state and control information available at the latent level. This allows us to obtain a self-contained controller capable of dealing with applications where continuous monitoring is unfeasible, or where the computational burden of full-order simulations does not meet the real-time requirement.
\\

The paper is organized as follows. Section~\ref{sec:OCP} briefly reviews open and closed-loop optimal control problems constrained by parametrized PDEs. Section~\ref{sec:OCP-DL-ROM} presents the real-time deep learning-based reduced order feedback control, delving the steps required in the offline and online phases, as well as the latent feedback loop. Section~\ref{sec:test} shows the performance of the proposed approach when dealing with two challenging optimal transport problems. Section~\ref{sec:conclusions} discusses some possible future developments and extensions of the presented methodology.

\section{From open-loop to feedback control of parametrized systems}\label{sec:OCP}

This section briefly introduces the class of optimization problems investigated throughout this work, explaining the role of the main variables and setting up the notation. We consider Optimal Control Problems (OCPs) for distributed dynamical systems, whose state evolutions can be described in terms of Partial Differential Equations (PDEs). The main goal in this context is to steer the physical system at hand towards a target configuration, influencing the state behaviour by optimally tuning a suitable control variable. Mathematically speaking, this is achieved through the following constrained optimization problem:

\begin{equation}
\label{eq:OCP}
    \mbox{Given} \ \yinit \ \mbox{and} \ \mus, \quad \mbox{find} \quad  J_h\left(\yh(t), \uh(t); \mus \right) \to \min \quad \text{ s.t. } \quad 
    \begin{cases}
    {\bf G}_h(\yh(t), \uh(t);  \mus) = {\bf 0} & \text{in } (0,T]\\ 
    \yh(0) = \yinit
    \end{cases}
\end{equation}
where $\yh(t): [0,T] \to \R^{N_h^y}$ and $\uh(t): [0,T] \to \R^{N_h^u}$ denote, respectively, the semi-discrete state and control, obtained discretizing the PDE through numerical techniques, such as, e.g., the Finite Element Method (FEM) \citep{Quarteroni2017}. Specifically, after splitting the domain $\Omega$ -- that is the region where the physical phenomenon takes place -- in sub-elements, the (scalar, for the sake of simplicity) infinite-dimensional state $y(\mathbf{x},t): \Omega \times [0,T] \to \R$ and control $u(\mathbf{x},t): \Omega \times [0,T] \to \R$ are approximated in suitable finite-dimensional spaces $\mathcal{Y}_h$ and $\mathcal{U}_h$ through basis expansions, that is
\[
y(\mathbf{x},t) \approx \sum_{i=1}^{N_h^y} \mathbf{y}_{h,i}(t) \xi_i^y(\mathbf{x}); \qquad u(t) \approx \sum_{i=1}^{N_h^u} \mathbf{u}_{h,i}(t) \xi_i^u(\mathbf{x})
\]
where $\{ \mathbf{y}_{h,i}(t)\}_{i=0}^{N_h^y}$ and $\{\mathbf{u}_{h,i}(t)\}_{i=0}^{N_h^u}$ are the elements of the vectors $\yh(t)$ and $\uh(t)$, while $\{\xi_i^y(\mathbf{x}) \}_{i=0}^{N_h^y}$ and $\{\xi_i^u(\mathbf{x}) \}_{i=0}^{N_h^u}$ are the space-varying basis functions spanning $\mathcal{Y}_h$ and $\mathcal{U}_h$, respectively. The discretization parameter $h>0$ corresponds to the characteristic dimension of the sub-elements: the smaller $h$, the better the FEM approximation, at the price of more degrees of freedom and more expensive computations. In the following, we assume complete access to the state values within the domain $\Omega$. However, the proposed approach can be easily extended in order to deal with different observed quantities defined as
\[
\zh(t) = O_h(t)\yh(t) \in \R^{N^z_h} \quad \mbox{where} \quad O_h(t): \R^{N_h^y} \to \R^{N^z_h} \quad \forall t \in [0,T],
\] 
allowing for partially observable states or state-related markers of the actual physical configuration.

The real-valued \textit{loss} or \textit{cost functional} $J_h(\yh(t),\uh(t); \mus)$ in Equation~\eqref{eq:OCP} evaluates the effectiveness of a specific state-control trajectory with respect to the target configuration: the lower the cost, the higher the performance of that control action. The optimal control trajectory, which brings the dynamical system as close to the target as possible, and the corresponding optimal state, can be thus computed by minimizing the loss functional $J_h$. 

\begin{remark}
\label{rem:loss}
Whenever a target configuration $\mathbf{y}_d \in \R^{N_h^y}$ has to be reached at final time $t = T$, the following loss functional may be employed
\[
J_h(\yh(t),\uh(t); \mus) = \dfrac{1}{2}|| \yh(T) - \mathbf{y}_d ||^2 + \dfrac{\beta}{2} \int_0^T ||\uh(t)||^2 dt 
\]
where $\left\Vert \cdot \right\Vert$ stands for the Euclidean norm. If, instead, a target trajectory $\mathbf{y}_d(t)$ must be followed throughout the time interval $[0,T]$, we may consider
\[ J_h(\yh(t),\uh(t); \mus) = \dfrac{1}{2}\int_0^T || \yh(t) - \mathbf{y}_d(t) ||^2 dt + \dfrac{\beta}{2} \int_0^T || \uh(t) ||^2 dt. \]
Note that, as detailed in \cite{Manzoni2021}, regularization terms concerning the norm of the control and its gradient are generally considered to guarantee the well-posedness of the OCP, as well as to prevent unfeasible and overconsuming control actions. The parameter $\beta > 0$ has to be chosen in order to properly balance the terms in the cost functional: the lower $\beta$, the closer the optimal state to the target configuration, the bigger the optimal control norm, and the harder the optimization procedure.
\end{remark}

The constraint in Equation~\eqref{eq:OCP} is crucial to guarantee the physical admissibility of the optimal state-control pair. The state and the control variables are indeed coupled by the dynamical system \acapo ${\bf G}_h(\yh(t), \uh(t); \mus) = {\bf 0}$~with ${\bf G}_h \in \R^{N_h^y}$ -- for the sake of simplicity, the state system is here considered directly in the semi-discrete formulation, which is regarded as the high-fidelity Full-Order Model (FOM) -- governing the evolution of the state system within the domain $\Omega$ and in the time interval $(0,T]$. In this work, we take into account dynamics described in terms of (possibly nonlinear) PDEs parametrized by a vector of scenario parameters $\mus \in \mathcal{P} \subset \R^p$, where $\mathcal{P}$ denotes the parameter space. The presence of parameters identifying the scenario variability to address increases the OCP complexity, as we aim to optimally control the dynamics in real-time for multiple scenarios. For example, in the applications detailed in Section~\ref{sec:test}, $\mus$ represents the target endpoint of the state trajectory, while the velocity field driving the state is regarded as control. Therefore, changes in the scenario entail completely different state routes and optimal controls to apply. Note that, along with the initial condition $\yh(0) = \yinit$, suitable boundary conditions must be set on $\partial \Omega \times (0,T]$ in order to guarantee the well-posedness of the PDE and the corresponding OCP  \citep{Manzoni2021}. Besides the physical constraint in Equation~\eqref{eq:OCP}, additional constraints may be taken into account whenever the state and the control are subject to extra physical or practical requirements. For example, it is possible to properly define the set of admissible state and control values and require that 
\[\yh(t) \in \mathcal{Y}_{ad} \subset \R^{N_h^y}, \quad \uh(t) \in \mathcal{U}_{ad} \subset \R^{N_h^u} \quad \forall t \in [0,T].\]
In the following, without loss of generality, no additional constraints are considered, that is $\mathcal{Y}_{ad} = \R^{N_h^y}$ and $\mathcal{U}_{ad} = \R^{N_h^u}$. 

\begin{remark}
Consider the semi-discrete formulation of a nonlinear time-dependent PDE
\[ {\bf G}_h(\yh(t), \uh(t); \mus) =  M_h(\mus) \dfrac{d\mathbf{y}_h(t)}{dt} + A_h(\mus) \mathbf{y}_h(t) + C_h(\mathbf{y}_h(t); \mus) - \mathbf{f}_{h}(t;\mus) - B_h(\mathbf{u}_h(t); \mus) = {\bf 0}
\]
where $M_h(\mus)$ is the mass matrix, $A_h(\mus)$ is the stiffness matrix, $B_h(\mus)$ is the control matrix, $C_h(\mus)$ is the nonlinear term, and $\mathbf{f}_{h}(t;\mus)$ is the external source. Let $\{ t_j \}_{j=0}^{N_t}$ an evenly spaced grid with step $\Delta t$ discretizing the interval $[0,T]$. After approximating the time derivative through, e.g., semi-implicit Euler method, the $N_h^y$ basis expansion coefficients $\yh(t_{j+1})$ for $j = 0,...,N_t-1$ can be computed by solving the associated system of $N_h^y$ equations 
\[ [M_h(\mus) + \Delta t A_h(\mus)] \yh(t_{j+1}) = M_h(\mus) \yh(t_j) - \Delta t C_h(\yh(t_j); \mus) + \Delta t \mathbf{f}_{h}(t_ {j+1}, \mus) + \Delta t B_h(\uh(t_j); \mus)
\]
\end{remark}

\begin{figure}[!ht]
    \centering
    \includegraphics[scale = 0.55]{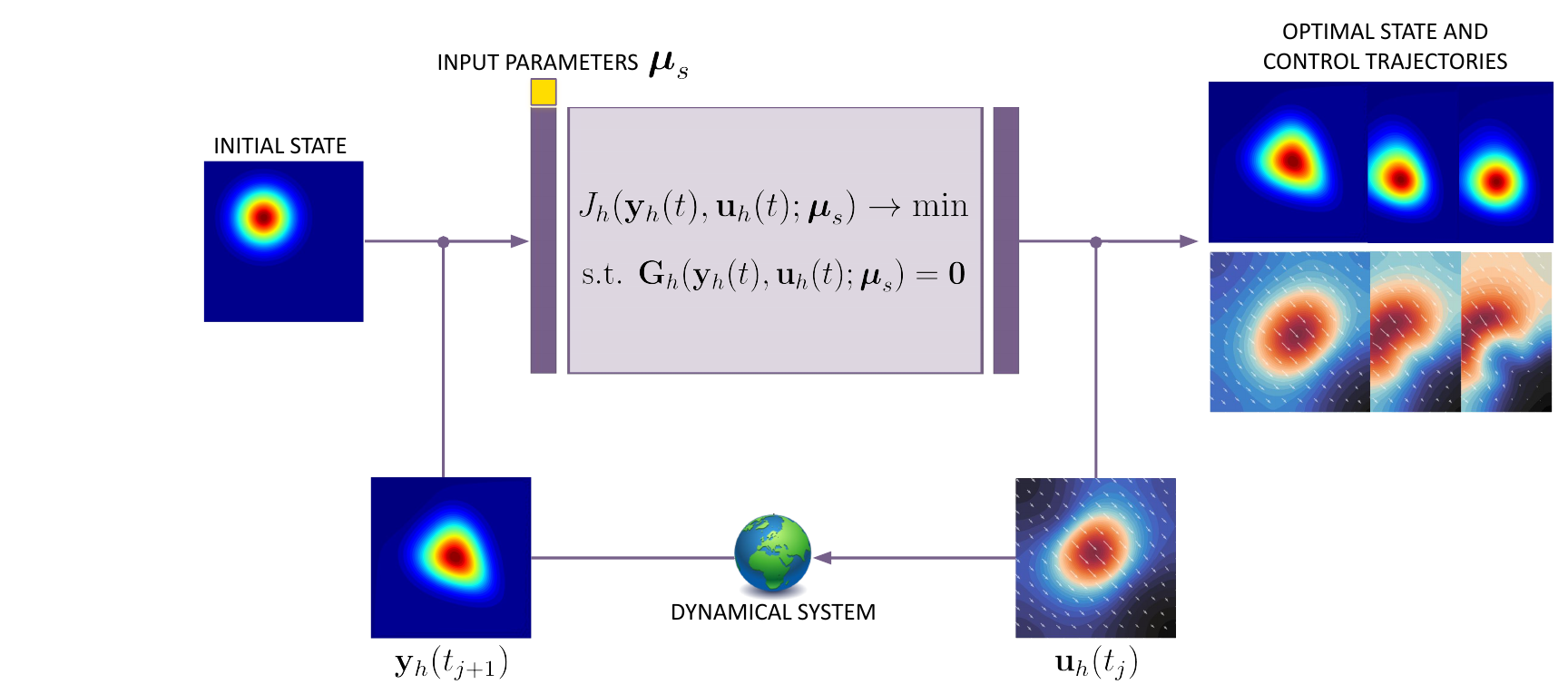}
    \caption{Feedback control scheme considering multiple optimal control problem resolutions and a high-fidelity full-order model of the dynamical system. Given the current state and the scenario parameters $\mus$, the open-loop optimal control and state trajectories are retrieved by solving the optimal control problem. The feedback signal is then recovered by measuring or simulating the state at the next time step.}
    \label{fig:feedback-loop}
\end{figure}

The OCP in Equation~\eqref{eq:OCP} provides, once solved, the open-loop optimal control trajectory related to the initial state $\yinit$ and the scenario $\mus$ under investigation. See, e.g., \cite{Manzoni2021} for a complete presentation and analysis of open-loop OCPs. Nevertheless, the open-loop solution is vulnerable to disturbances and uncertainties in the system, and may compromise system stability and performance. To address this drawback, feedback control strategies may be exploited, in which the considered OCP is solved repeatedly over the time horizon, continuously updating the initial state with the current observed state, as depicted in Figure~\ref{fig:feedback-loop} with reference to the first test case detailed in Section~\ref{sec:test}. In practice, the delay between the availability of the current state and the computation of the subsequent optimal action is unavoidable, making this strategy not suitable for applications with strict timing requirements. The time lag becomes even larger when online sensor data are lacking and synthetic state data must be simulated at each time step, especially in case of a large full-order dimension $N_h^y$. The computational bottleneck may be overcome by fast-evaluable surrogate models approximating the full-order dynamics or by merely considering low-dimensional variables, as often happens in the Model Predictive Control (MPC) context, as detailed in Section~\ref{sec:introduction}. However, in order to design a real-time feedback controller capable of dealing with high-dimensional variables, as well as complex nonlinear dynamics, we aim to avoid any optimization procedure online, concentrating the computationally expensive steps in an offline phase to be performed only once.

Similarly to the Reinforcement Learning (RL) framework, it is convenient to rephrase the problem focusing on the approximation of the so-called policy function
\begin{equation}
\label{eq:policy}
    \pi_h: \R^{N_h^y} \times \mathcal{P} \to \R^{N_h^u}, \quad \uh = \pi_h \left(\yh, \mus \right) \quad \mbox{where} \quad 
    \uh = \arg\min J_h\left(\yh(\uh), \uh; \mus \right)
\end{equation}
where, for the sake of compactness, the reduced cost functional $J_h\left(\yh(\uh), \uh; \mus \right)$ directly considers the state-control dependence $\yh(\uh)$ given implicitly by the PDE ${\bf G}_h(\yh, \uh; \mus) = {\bf 0}$. The policy function corresponds to the one time-step solution map of the OCP in Equation~\eqref{eq:OCP}: indeed, starting from the current observed state $\yh$ and the scenario parameters $\mus$, it returns the corresponding subsequent optimal control minimizing the loss $J_h$. Note that, since the policy returns only the current optimal control in place of the entire trajectory over time, the time $t$ is not made explicit anymore. To recover the optimal control at every time instant, multiple policy evaluations are required with a continuous state update. As detailed in the next section, the policy input-output dimensions may be high-dimensional -- especially when dealing with distributed state and control variables showing large $N_h^y$ and $N_h^u$ full-order dimensions -- and model order reduction techniques are crucial to retrieve faster and lighter map.

\section{Deep learning-based reduced order feedback control}\label{sec:OCP-DL-ROM}

This section presents the non-intrusive deep learning-based ROM strategy aiming at providing the optimal control action starting from an observed state and a scenario encoded in a vector of input parameters $\mus$. Our approach focuses on an offline-online decomposition, where {\em (i)} the offline phase is concerned with the computationally expensive steps needed to build and calibrate a low-dimensional fast-evaluable surrogate model for the policy function, such as data generation, dimensionality reduction and neural networks training, while {\em (ii)} real-time optimal control actions may be inferred by evaluating the constructed policy in the online phase. This pipeline is inspired by the non-intrusive ROMs proposed by \cite{Hesthaven2018, Fresca2021, Fresca2022} -- namely POD-NN, DL-ROM and POD-DL-ROM -- and adapted to a feedback control framework. The proposed offline-online decomposition is crucial in a twofold manner: first of all, optimal control snapshots are generated and properly reduced offline, thus enabling us to easily handle distributed or boundary controls. This would not be easily implemented within, e.g., RL strategies, where the learning requires a continuous interaction between the agent and the dynamical system. Moreover, in contrast with the MPC framework, the optimization procedures are performed during data generation (offline), thus making it possible to retrieve the optimal control actions (online) through fast policy evaluations. The next sections explore in details the different steps required to build and employ the proposed method.

\subsection{Offline phase}

The rationale of the proposed approach is to approximate the policy function in Equation~\eqref{eq:policy} in a supervised manner through a feedforward neural network. To do so, {\em (i)} state-control pairs have to be generated exploring different initial conditions $\mathbf{y}_0$ and different scenarios $\mus$ in the parameter space $\mathcal{P}$, and {\em (ii)} neural network training must be performed to synthesize a surrogate policy that accurately predicts the optimal control starting from the corresponding observed state. Since the full-order state and control dimensions $N_h^y$ and $N_h^u$ may be remarkably high, an intermediate step is taken into account to compress the snapshots dimensionality and sharply reduce the input-output policy layers.

\paragraph{Data generation}

The first preparatory procedure to build a policy approximation through supervised learning strategies is the generation of policy input-output pairs, that are optimal state and control trajectories. In order to acquire a robust policy capable of controlling the dynamical system across a wide variety of configurations, it is crucial to extensively explore different initial states and scenarios. While sampling strategies are usually enough to explore the (typically low-dimensional) parameter space $\mathcal{P}$, smart choices are required to sample from the (possibly high-dimensional) state space $\R^{N_h^y}$. For instance, in the applications detailed in Section~\ref{sec:test}, we parametrize the initial state field, so that state variability is reduced, while still covering all the possible relevant cases of interest that may occur. Every optimal state and control trajectory can be computed by solving the OCP in Equation~\eqref{eq:OCP} while considering $\yinit = \yinit^{(i)}$ and $\mus = \mus^{(i)}$ for $i=1,...,N_s$. To this aim, it is possible to solve the system of Karush-Kuhn-Tucker (KKT) optimality conditions derived from Equation~\eqref{eq:OCP} via Lagrange multipliers' method, that is
\begin{equation}
\label{eq:KKT}
\left\{
\begin{array}{ll}
\nabla_y J_h (\yh(t), \uh(t);\mus) 
 + (\partial_y {\bf G}_{h} (\yh(t), \uh(t);\mus))^{\top }\mathbf{p}_h(t)= {\bf 0}& \text{ \ \ (adjoint equation)}  \smallskip \\
 \mathbf{p}_h(T) = \nabla_y J_h (\yh(T), \uh(T); \mus) & \text{ \ \ (final condition)}  \smallskip \\
 \nabla_u J_h (\yh(t), \uh(t);\mus) + (\partial_u {\bf G}_{h} (\yh(t), \uh(t);\mus))^{\top }{\bf p}_h(t)  ={\bf 0}   & \text{ \ \ (optimality condition)} \smallskip \\
 {\bf G}_h\left( \yh(t), \uh(t); \mus \right) ={\bf 0} &  \text{ \ \ (state equation)}. \smallskip \\
 \yh(0) = \yinit & \text{ \ \ (initial condition)} \smallskip%
\end{array}%
\right.
\end{equation}
where ${\bf p}_h(t) \in \mathbb{R}^{N_h^p}$ for $t \in [0,T]$ is the semi-discrete adjoint vector.  For a complete presentation of the OCP solving methods see, e.g., \cite{Manzoni2021}.

\begin{remark}
In case of final target observation, which corresponds to the first example presented in Remark~\ref{rem:loss}, the final condition for the adjoint equation ends up being
\[
\mathbf{p}_h(T) = \yh(T) - \mathbf{y}_d.
\]
Instead, whenever the target state observation is distributed in $[0,T]$ as in the second example of Remark~\ref{rem:loss}, the final condition becomes
\[ 
\mathbf{p}_h(T) = {\bf 0}.
\]
\end{remark}

The system of optimality conditions in Equation~\eqref{eq:KKT} returns -- after being discretized through FEM and solved via, e.g., the Newton method -- the optimal state and control trajectories 
\[\{\yh^{(i)}(t_0), \yh^{(i)}(t_1),...,\yh^{(i)}(t_{N_t-1}), \yh^{(i)}(t_{N_t}) \}_{i = 1,...,N_s}, \qquad \{ \uh^{(i)}(t_0),\uh^{(i)}(t_1),...,\uh^{(i)}(t_{N_t-2}),\uh^{(i)}(t_{N_t-1}) \}_{i = 1,...,N_s}
\]
where $\{ t_j \}_{j=0}^{N_t}$ introduces a uniform partition of the time horizon $[0,T]$. Starting from the generated trajectories, it is then possible to assemble the $N_t N_s$ input-output pairs useful for training and testing the policy surrogate model, that is
\begin{equation*}
\{ (\yh^{(i)}(t_j), \mus^{(i)}), \uh^{(i)}(t_j) \}_{
\begin{subarray}{l}
i = 1,...,N_s
\\
j = 0,...,N_t-1
\end{subarray} }
\quad \text{ where } \quad  (\yh^{(i)}(t_j), \mus^{(i)}) \in \R^{N_h^y} \times \mathcal{P}, \quad \uh^{(i)}(t_j) \in R^{N_h^u}.
\end{equation*}
In the following, without loss of generality, whenever the dependence on time is not made explicit, the snapshots are reordered and denoted by
\[
\{\yh^{(k)}, \mus^{(k)}, \uh^{(k)}\}_{k = 1,...,N_t N_s}.
\]
 
\paragraph{Dimensionality reduction}

Reduced Order Models (ROMs) are very helpful when dealing with parametrized or many-query problems, as the number of degrees of freedom into play may be dramatically reduced by extracting the (few) essential features relevant for the problem, ending up with faster (but still very accurate) methods. The dimensionality reduction is typically performed by projecting the available snapshots onto lower-dimensional subspaces, which may be linear -- such as in the Proper Orthogonal Decomposition (POD) framework -- or nonlinear -- as when autoencoders (AEs) are employed.

\begin{itemize}
    
    \item \textbf{Proper Orthogonal Decomposition}: full-order state and control snapshots may be projected onto linear subspaces of dimensions, respectively, $N_y \ll N_h^y$ and $N_u \ll N_h^u$ by POD \cite{Sirovich1987}. The $N_y$ and $N_u$ basis elements (also known as POD modes) spanning the two latent subspaces are computed by the Singular Values Decomposition (SVD) of the state and control snapshots matrices, and are then collected column-wise in the matrices $\mathbb{V}_y \in \R^{N_h^y \times N_y}$ and $\mathbb{V}_u \in \R^{N_h^u \times N_u}$. The resulting projections thus read as follows:
    \begin{alignat*}{2}
    &\yn =  \mathbb{V}_y^{\top} \yh, \qquad \yhrec =  \mathbb{V}_y \yn
    \\
    &\un = \mathbb{V}_u^{\top} \uh, \qquad \uhrec = \mathbb{V}_u \un
    \end{alignat*}
    where $\yhrec \in \R^{N_h^y}$ and $\uhrec \in \R^{N_h^u}$ are the full-order state and control projections approximating $\yh$ and $\uh$ up to a reconstruction error. Therefore, instead of dealing with the full-order state and control variables, it is possible to focus on a set of low-dimensional features extracted from their basis expansions coefficients $\yn \in \R^{N_y}$ and $\un \in \R^{N_u}$, respectively. Note that the latent dimensions $N_y$ and $N_u$ are usually selected by a trade-off between having small latent subspaces and low reconstruction errors \cite{Hesthaven2015, Quarteroni2015}. 
    
    \item \textbf{Autoencoders}: linear ROMs may require a remarkably high number of POD modes in complex settings characterized by, e.g., involved scenario-state or scenario-control dependencies, as well as nonlinear or transport-dominated PDEs, compromising the speed up of the method. To recover more effective embeddings, projections onto nonlinear subspaces may be exploited, under the form    
    \begin{alignat*}{2}
    &\yn = \varphi_E^y(\yh), \qquad \yhrec = \varphi_D^y(\yn)
    \\
    &\un = \varphi_E^u(\uh), \qquad \uhrec = \varphi_D^u(\un)
    \end{alignat*}
    where $\varphi_E^y:\R^{N_h^y} \to \R^{N_y}$, $\varphi_E^u:\R^{N_h^u} \to \R^{N_u}$, $\varphi_D^y:\R^{N_y} \to \R^{N_h^y}$ and $\varphi_D^u:\R^{N_u} \to \R^{N_h^u}$ are nonlinear functions compressing high-dimensional snapshots into latent representations and viceversa. Throughout this work, we model encoding and decoding mappings through the so-called autoencoders, which are neural networks showing a bottleneck architecture. The nonlinearity feature is due to the activation functions exploited within the autoencoders, such as the leaky ReLU function, as taken into account in the following. Note that, as proposed by \cite{Fresca2021, Franco2023} in the framework of Deep Learning-based Reduced Order Models (DL-ROMs), convolutional autoencoders may be employed for a remarkable dimensionality reduction.
 
    \item \textbf{POD+AE}: the speed of POD and the power of autoencoders can also be combined, as initially proposed by \cite{Fresca2022} in the POD-DL-ROM framework. In this context -- here referred to as POD+AE -- autoencoders are trained to further compress the dimensionality of snapshots that has been preliminarily reduced by POD, that is,
    \begin{alignat*}{2}
    &\yn = \varphi_E^y(\mathbb{V}_y^{\top}\yh), \quad &&\yhrec = \mathbb{V}_y \varphi_D^y(\yn)
    \\
    &\un = \varphi_E^u(\mathbb{V}_u^{\top} \uh), \quad &&\uhrec = \mathbb{V}_u \varphi_D^u(\un).
    \end{alignat*}
    In this way we can obtain lighter NNs architecture and faster trainings, without compromising the overall accuracy.

\end{itemize}

\paragraph{Low-dimensional policy} After selecting a suitable compression technique to reduce state and control data, it is possible to define a feedforward neural network approximating the policy function at the latent level, that is,
\[
\pi_N: \R^{N_y} \times \mathcal{P} \to \R^{N_u}, \quad \tildeun = \pi_N \left(\yn, \mus \right),
\]   
where $\tildeun \in \R^{N_u}$ is the approximation of the control embedding $\un \in \R^{N_u}$ given by the surrogate policy $\pi_N$. Differently from the full-order policy function $\pi_h$ in Equation~\eqref{eq:policy}, the input-output spaces are now reduced. Lighter neural network architectures are crucial to speed up both training and evaluations, especially when dealing with real-time applications involving high-dimensional variables.

\paragraph{Neural network training}
The weights and biases of the aforementioned neural networks -- namely, the autoencoders and the policy function -- are trained all at once within a single optimization procedure. Doing so, an implicit link is established between the state and control reductions and the low-dimensional policy, obtaining latent coordinates that are meaningful in view of the policy-based control. In particular, the following cumulative loss function is taken into account
\[
J_{\mathrm{NN}} = 
\lambda_1 J_{\mathrm{rec}}^y + \lambda_2 J_{\mathrm{rec}}^u + J_{\pi_N}.
\]
The reconstruction errors $J_{\mathrm{rec}}^y$ and $J_{\mathrm{rec}}^u$ are given by

\begin{equation}
J_{\mathrm{rec}}^y = \begin{cases}
0 & \text{in case a POD-NN is used,}
\\
J_{\mathrm{AE}}^y & \text{in case a DL-ROM is used,}
\\
J_{\mathrm{POD+AE}}^y & \text{in case a POD-DL-ROM is used;}
\end{cases} \qquad 
J_{\mathrm{rec}}^u = \begin{cases}
0 & \text{in case a POD-NN is used,}
\\
J_{\mathrm{AE}}^u & \text{in case a DL-ROM is used,}
\\
J_{\mathrm{POD+AE}}^u & \text{in case a POD-DL-ROM is used,}
\end{cases}
\label{eq:Jrec}
\end{equation}
where the mean square reconstruction errors
\begin{align*}
J_{\mathrm{AE}}^y &= \dfrac{1}{|I_{\text{train}}|} \sum_{k \in I_{\text{train}}}  \left\lVert \yh^{(k)} - \yhrec^{(k)} \right\rVert^2 = \dfrac{1}{|I_{\text{train}}|} \sum_{k \in I_{\text{train}}} \left\lVert \yh^{(k)} - \varphi_D^y(\varphi_E^y(\yh^{(k)})) \right\rVert^2
\\
J_{\mathrm{AE}}^u &= \dfrac{1}{|I_{\text{train}}|} \sum_{k \in I_{\text{train}}}  \left\lVert \uh^{(k)} - \uhrec^{(k)} \right\rVert^2 = \dfrac{1}{|I_{\text{train}}|} \sum_{k \in I_{\text{train}}} 
\left\lVert \uh^{(k)} -\varphi_D^u(\varphi_E^u(\uh^{(k)})) \right\rVert^2
\\
J_{\mathrm{POD+AE}}^y &= \dfrac{1}{|I_{\text{train}}|} \sum_{k \in I_{\text{train}}}  \left\lVert \mathbb{V}_y^{\top} \yh^{(k)} -\mathbb{V}_y^{\top} \yhrec^{(k)} \right\rVert^2 = \dfrac{1}{|I_{\text{train}}|} \sum_{k \in I_{\text{train}}}  \left\lVert \mathbb{V}_y^{\top} \yh^{(k)} - \varphi_D^y(\varphi_E^y(\mathbb{V}_y^{\top} \yh^{(k)})) \right\rVert^2
\\
J_{\mathrm{POD+AE}}^u &= \dfrac{1}{|I_{\text{train}}|} \sum_{k \in I_{\text{train}}} \left\lVert \mathbb{V}_u^{\top} \uh^{(k)} - \mathbb{V}_u^{\top} \uhrec^{(k)} \right\rVert^2 = \dfrac{1}{|I_{\text{train}}|} \sum_{k \in I_{\text{train}}} \left\lVert \mathbb{V}_u^{\top} \uh^{(k)} - \varphi_D^u(\varphi_E^u(\mathbb{V}_u^{\top} \uh^{(k)})) \right\rVert^2
\end{align*}
are computed on a set of training snapshots with indices $I_{\text{train}} \subset \{ 1,...,N_t N_s \}$ with, typically, $|I_{\text{train}}| \approx 0.8 N_t N_s$. Moreover, the mean square prediction error entailed by the policy on training data is computed as
\begin{align}
\begin{split}
J_{\pi_N} &= \dfrac{1}{|I_{\text{train}}|} \sum_{k \in I_{\text{train}}} \underset{\text{Latent space error}}{\underbrace{\left\lVert \un^{(k)} - \tildeun^{(k)} \right \rVert^2}} + \lambda_3 \underset{\text{After-decoding error}}{\underbrace{\left\lVert \varphi_D^u(\un^{(k)}) - \varphi_D^u(\tildeun^{(k)}) \right \rVert^2}} 
\\
&= \dfrac{1}{|I_{\text{train}}|} \sum_{k \in I_{\text{train}}} \left\lVert \un^{(k)} - \pi_N(\yn^{(k)}, \mus^{(k)}) \right \rVert^2 + \lambda_3 \left\lVert \varphi_D^u(\un^{(k)}) - \varphi_D^u(\pi_N(\yn^{(k)}, \mus^{(k)})) \right \rVert^2. 
\end{split}
\label{eq:loss_policy}
\end{align}

\begin{figure}[!ht]
    \centering
    \includegraphics[scale = 0.55]{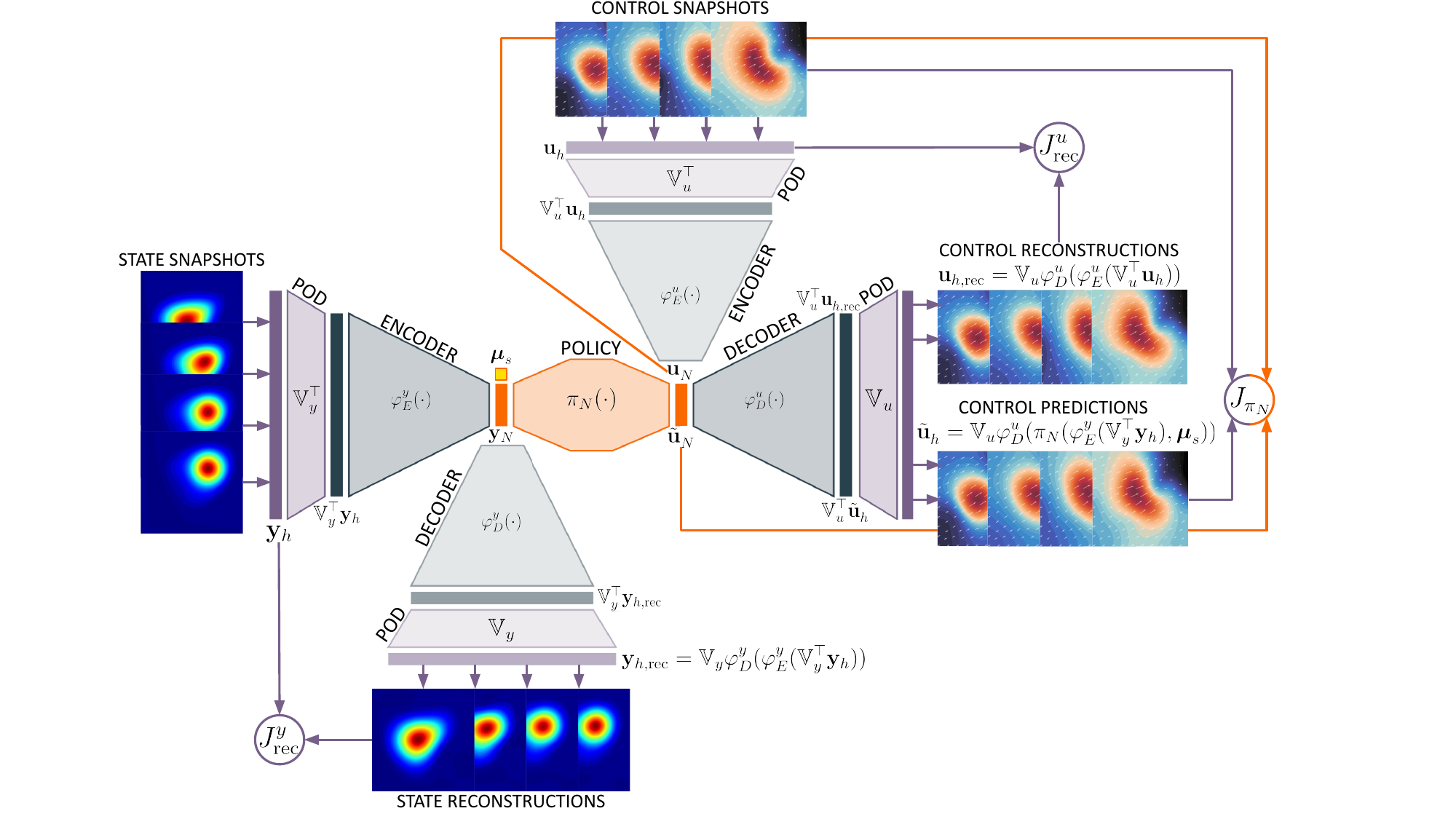}
    \caption{Offline phase of the deep learning-based reduced order feedback controller. After generating optimal state and control snapshots through the adjoint method, the state and control autoencoders, namely $\varphi_D^y(\varphi_E^y(\cdot))$ and $\varphi_D^u(\varphi_E^u(\cdot))$, and the surrogate policy $\pi_N$ are trained minimizing the cumulative loss function $J_{\mathrm{NN}} = \lambda_1 J_{\mathrm{rec}}^y + \lambda_2 J_{\mathrm{rec}}^u + J_{\pi_N}$.}
    \label{fig:OCP-DL-ROM-offline}
\end{figure}

The after-decoding loss term in Equation~\eqref{eq:loss_policy} --
which is considered whenever the control is reduced through AE or POD+AE -- is helpful in training the control decoder in accordance with the surrogate policy, thus achieving more accurate decodings of the policy predictions and better results at the full-order levels. The hyperparameters $\lambda_1, \lambda_2, \lambda_3 \in \R$ must be chosen properly to balance the magnitudes of the different terms within $J_{\mathrm{NN}}$. A summary of the offline phase is available in Figure~\ref{fig:OCP-DL-ROM-offline} when considering the snapshots of the first test case detailed in Section~\ref{sec:test} and when taking into account POD+AE as reduction strategy both for the state and the control.

\subsection{Online phase}
\begin{figure}[!ht]
    \centering
    \includegraphics[scale = 0.55]{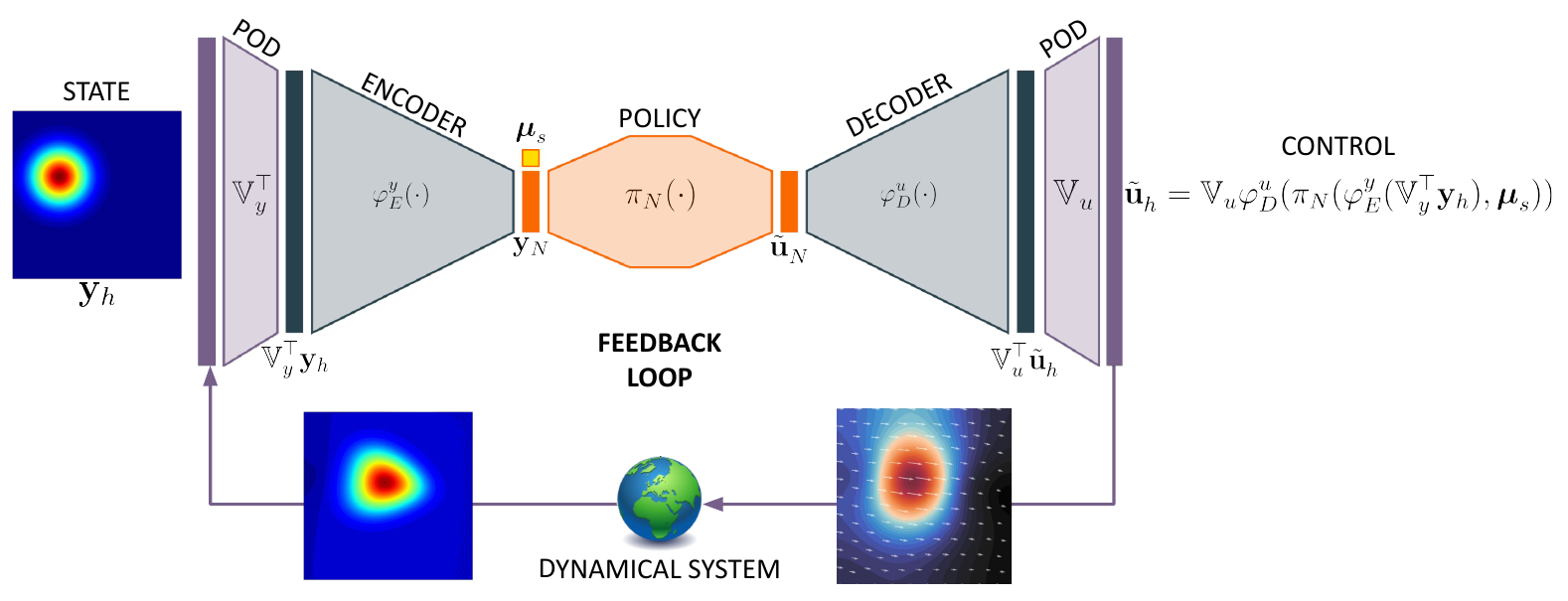}
    \caption{Online phase of the deep learning-based reduced order feedback controller. The optimal full-order control action corresponding to the observed state $\yh$ in a scenario described by input parameters $\mus$ is inferred online through forward passes of $\pi_N$ and the encoding-decoding mappings.}
    \label{fig:OCP-DL-ROM-online}
\end{figure}

Whenever a new state $\yh^{\mathrm{new}}$ is observed in a scenario $\mus^{\mathrm{new}}$ unseen during training, the corresponding optimal control action is inferred through a forward pass of the low-dimensional policy and the encoding-decoding maps, that is
\[
\tilde{\mathbf{u}}_{h}^{\mathrm{new}} = \phi_D^u(\tilde{\mathbf{u}}_{N}^{\mathrm{new}}) = \phi_D^u(\pi_N(\phi_E^y(\yh^{\mathrm{new}}), \mus^{\mathrm{new}}))
\]
where
\begin{equation}
\begin{split}
\phi_E^y(\yh^{\mathrm{new}}) = \begin{cases}
\mathbb{V}^{\top}_y \yh^{\mathrm{new}} & \text{in case a POD-NN is used,}
\\
\varphi_E^y(\yh^{\mathrm{new}}) & \text{in case a DL-ROM is used,}
\\
\varphi_E^y(\mathbb{V}^{\top}_y \yh^{\mathrm{new}}) & \text{in case a POD-DL-ROM is used;}
\end{cases}
\\
\phi_D^u(\tilde{\mathbf{u}}_{N}^{\mathrm{new}}) = \begin{cases}
\mathbb{V}_u \tilde{\mathbf{u}}_{N}^{\mathrm{new}} & \text{in case a POD-NN is used,}
\\
\varphi_D^u(\tilde{\mathbf{u}}_{N}^{\mathrm{new}}) & \text{in case a DL-ROM is used,}
\\
\mathbb{V}_u \varphi_D^u(\tilde{\mathbf{u}}_{N}^{\mathrm{new}}) & \text{in case a POD-DL-ROM is used.}
\end{cases}
\end{split}
\label{eq:Phi}
\end{equation}
In practice, starting from the observed high-dimensional state directly measured from the dynamical system at hand, we predict the distributed control responsible for steering the system towards the optimal configuration. To this aim, state reduction and control decoding are crucial to move back and forth from the full-order dimensions to the latent level where the surrogate policy is employed. New control predictions can be then retrieved in loop whenever new state measurements are available. A sketch of the proposed feedback loop is depicted in Figure~\ref{fig:OCP-DL-ROM-online} when considering the snapshots of the first test case detailed in Section~\ref{sec:test} and when taking into account POD+AE as reduction strategy both for the state and the control. Differently from intrusive techniques -- such as, e.g., the Reduced Basis method \cite{Hesthaven2015, Quarteroni2015} where the parameter-to-solution map is retrieved by properly projecting and solving the system of optimality conditions in Equation~\eqref{eq:KKT} -- the proposed deep-learning based reduced order modeling strategy is non-intrusive, resulting in a very general and flexible tool capable of retrieving the optimal control action in real-time for a wide range of control problems.

\subsection{Latent feedback loop}
\label{subsec:latentfeedbackloop}

The deep-learning based feedback control strategy introduced so far relies on the real-time availability of full-order state information in the online phase. Indeed, when dealing with feedback controllers, the optimal control action is inferred by looking at quantities observed from the dynamical system at hand. Specifically, the input of the surrogate policy $\pi_N$ is the current full-order state variable in its reduced form, which provides a comprehensive overview of the configuration and evolution of the system. However, continuous monitoring of the dynamical system may be unfeasible, especially when dealing with high-dimensional data. Indeed, the computational burden to generate synthetic state data online through high-fidelity FOM simulations -- as taken into account throughout the test cases in Section~\ref{sec:test} -- may not meet strict timing requirements, thereby compromising the feedback controller performances. Instead, whenever the state snapshot is measured through sensors widespread in the domain, continuous monitoring may be compromised by, e.g., sensor malfunctions, delays in data processing or temporary absence of signal for communications. To overcome these limitations, here we propose a feedback loop closure at the latent level that allows the model to rapidly predict the control action even in absence of state measurements online. Specifically, we aim to train a neural network to approximate the forward model (i.e. the time advancing scheme of the dynamical system) at the latent level, that is
\[
\varphi_N: \R^{N_y} \times \R^{N_u} \times \mathcal{P} \to \R^{N_y}, \quad \tildeyn (t_{j}) = \varphi_N \left(\yn (t_{j-1}), \un (t_{j-1}), \mus \right) \qquad \forall j = 1,...,N_t.
\]

As previously discussed, to preserve a better consistency between the dimensionality reduction and the maps defined in the latent spaces, we employ a single optimization step to train all the neural networks into play, that are the autoencoders, the policy $\pi_N$ and the forward map $\varphi_N$. Specifically, the cumulative loss function now becomes
\[
J_{\mathrm{NN}} = 
\lambda_1 J_{\mathrm{rec}}^y + \lambda_2 J_{\mathrm{rec}}^u + J_{\pi_N} + J_{\varphi_N}
\]
where
\begin{equation}
    \label{eq:latentlooperror}
    J_{\varphi_N} = \dfrac{1}{|I'_{\text{train}}|} \sum_{(i,j) \in I'_{\text{train}}} \lambda_4 \underset{\text{Prediction-from-data error}}{\underbrace{\left\lVert \yn^{(i)}(t_{j}) - \tildeyn^{(i)}(t_{j}) \right \rVert^2}} 
+ \lambda_5\underset{\text{Prediction-from-policy error}}{\underbrace{\left\lVert \yn^{(i)}(t_{j}) - \hat{\mathbf{y}}_n^{(i)}(t_{j}) \right \rVert^2}} +
\lambda_6 \underset{\text{After-decoding error}}{\underbrace{\left\lVert \varphi_D^y(\yn^{(i)}(t_{j})) - \varphi_D^y(\tildeyn^{(i)}(t_{j})) \right \rVert^2}} 
\end{equation}
while $J_{\mathrm{rec}}^y, J_{\mathrm{rec}}^u$ and $J_{\pi_N}$ are defined in Equation~\eqref{eq:Jrec} and Equation~\eqref{eq:loss_policy}, respectively. The state predictions appearing in Equation~\eqref{eq:latentlooperror} are the approximations of future state values starting from, respectively, the available training control data and policy-based control predictions, that are 
\[
\tildeyn^{(i)}(t_{j}) = \varphi_N(\yn^{(i)}(t_{j-1}), \un^{(i)}(t_{j-1}), \mus^{(i)}),
\qquad 
\hat{\mathbf{y}}_N^{(i)}(t_{j}) = \varphi_N(\yn^{(i)}(t_{j-1}), \pi_N(\yn^{(i)}(t_{j-1}), \mus^{(i)}), \mus^{(i)}) \qquad \forall j = 1,...,N_t, 
\]
while $I'_{\text{train}} \subset \{1,...,N_s\} \times \{ 1,...,N_t\}$ denotes the set of training indices. The prediction-from-data and prediction-from-policy error terms in Equation~\eqref{eq:latentlooperror} are useful to achieve accurate forward-in-time states at the latent level, both starting from control data and policy outputs. Instead, the after-decoding error term is helpful in obtaining acceptable results also at higher dimensions whenever AE or POD+AE reduction strategies are exploited. Here, three additional hyperparameters $\lambda_4, \lambda_5, \lambda_6 \in \R$ are considered to balance the terms in the cumulative loss function.

\begin{figure}[!ht]
    \centering
    \includegraphics[scale = 0.53]{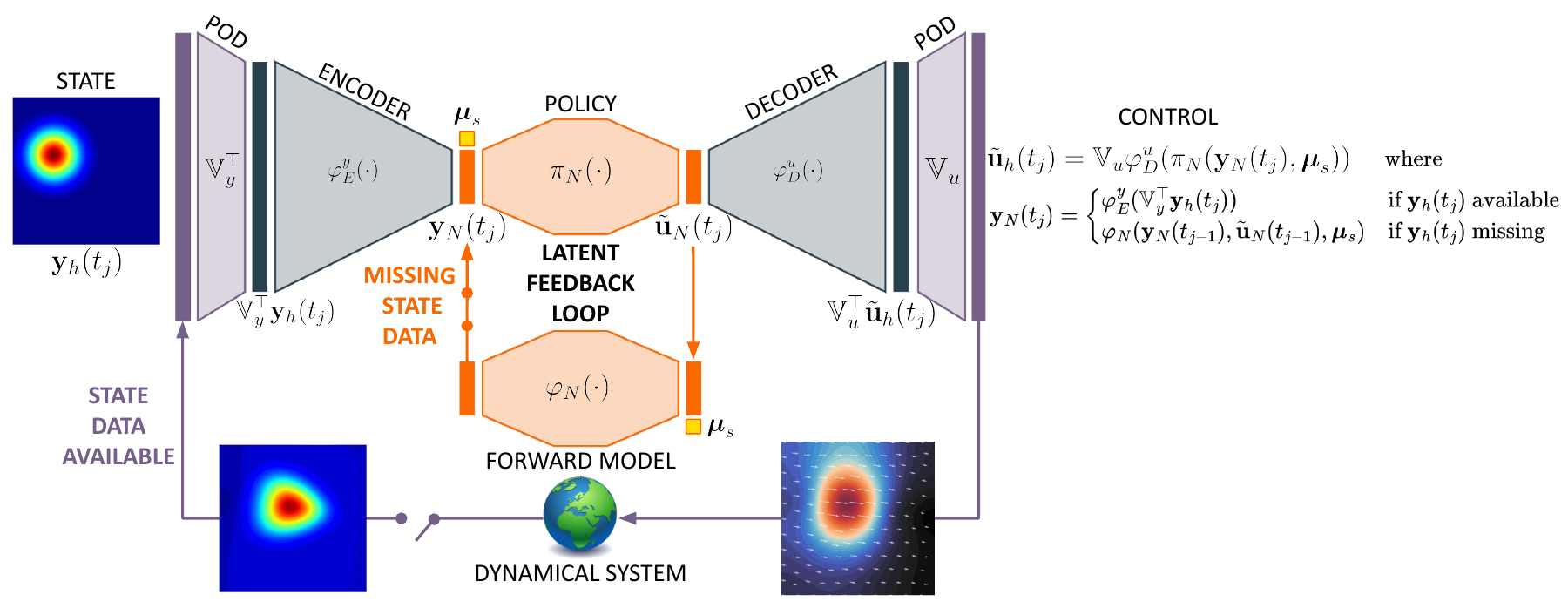}
    \caption{Online phase of the deep learning-based reduced order feedback controller with latent feedback loop. The optimal full-order control action corresponding to the observed state $\yh$ in a scenario described by input parameters $\mus$ is inferred online through forward-passes of $\pi_N$ and the encoding-decoding mappings. Whenever full-order state data $\yh$ are not available online, the trained deep learning-based forward model $\varphi_N$ is exploited to predict the state evolution, allowing for a continuous prediction of the control action.}
    \label{fig:latent-loop}
\end{figure}

The proposed feedback latent loop is depicted in Figure~\ref{fig:latent-loop} when considering the snapshots of the first test case detailed in Section~\ref{sec:test} and when taking into account POD+AE as reduction strategy both for the state and the control. Whenever full-order state data are available online, either in the form of sensor measurements or synthetic data simulated via a high-fidelity FOM, these may be exploited to infer the related optimal control action to be applied on the dynamical system, that is
\[
\tilde{\mathbf{u}}_h (t_j) = \phi_D^u(\pi_N(\yn(t_j), \mus)) \quad \text{where} \quad \yn(t_j) = \phi_E^y(\yh(t_j)) \qquad \forall j = 0,...,N_t-1
\]
with $\phi_D^u$ and $\phi_E^y$ defined in Equation~\eqref{eq:Phi}. Whether the full-order state snapshot is missing due to, e.g., latency concerns or sensors failures, state predictions at the latent level are provided by a forward pass of $\varphi_N$, and optimal control can still be predicted up to an approximation error through
\[
\tilde{\mathbf{u}}_h (t_j) = \phi_D^u(\pi_N(\tilde{\mathbf{y}}_N(t_{j}), \mus)) \quad \text{where} \quad \tilde{\mathbf{y}}_N(t_{j}) = \varphi_N(\yn(t_{j-1}), \tilde{\mathbf{u}}_N(t_{j-1}), \mus) \qquad \forall j = 1,...,N_t-1,
\]
resulting in a continuous control of the dynamical system.

%%%%%%%%%%%%%%%%%%%%%%%%%%%%%%%%%%%%%%%%%%%%%%%%%%%%%%%%

\section{Numerical results}\label{sec:test}

This section details the numerical tests performed to assess the effectiveness of the proposed approach. Specifically, we focus on two different time-dependent optimal transport problems in two spatial dimensions where, by optimally tuning the velocity field responsible for the state movement over space and time, the state has to be steered from its starting configuration to the desired final destination. In this context, our aim is to synthesize a feedback controller capable of dealing with both different initial settings and different target locations to be chosen online. To do so, we explore the initial state variability, as well as different target locations, in the data generation procedure. To better visualize the problem, note that the state variable may represent the density of particles systems described through a mean-field model, such as swarms of autonomous robots \cite{Sinigaglia2022a, Sinigaglia2022b, Sinigaglia2022c} delivering goods to the target location, while avoiding collisions with boundaries and obstacles along the route. 

Throughout this section, the following mean relative errors are employed to fairly evaluate the prediction accuracy of the proposed approach:
\[
\varepsilon^y_{\mathrm{rel}} = \dfrac{1}{|I_{\text{test}}|} \sum_{k \in I_{\text{test}}} 
 \dfrac{\lVert \yh^{(k)} - \tildeyh^{(k)} \rVert}{\lVert \yh^{(k)} \rVert}, \qquad \varepsilon^u_{\mathrm{rel}} = \dfrac{1}{|I_{\text{test}}|} \sum_{k \in I_{\text{test}}} 
 \dfrac{\lVert \uh^{(k)} - \tildeuh^{(k)}\rVert}{\lVert \uh^{(k)}\rVert}
\]
where $I_{\text{test}} = \{1,...,N_t N_s \} \setminus I_{\text{train}}$ is the set containing the indices of test data, that are the snapshots exploited only for evaluation purposes. Note that, when evaluating the reconstruction capabilities of the chosen reduction method, $\yhrec^{(k)}$ and $\uhrec^{(k)}$ are taken into account in place of $\tildeyh^{(k)}$ and $\tildeuh^{(k)}$, respectively.

\subsection{Optimal transport in a vacuum}
\label{subsec:vacuum}

This section is devoted to the application of the proposed real-time feedback controller to an optimal transport problem in a vacuum -- that is, for the sake of simplicity, we neglect the effects of surrounding fluid, such as air or water, in the space-time domain $\Omega \times (0,T]$, with final time $T>0$. The state dynamics is described by the Fokker-Planck equation (also known as Kolmogorov forward equation)
\begin{equation}
\label{eq:FP}
\begin{cases}
      \dfrac{\partial y}{\partial t} + \nabla \cdot (- \nu \nabla y + \mathbf{u} y) = 0 \qquad
      & \text{in} \ \Omega \times (0,T]
      \\
      (-\nu \nabla y + \mathbf{u} y) \cdot \mathbf{n} = 0
      \qquad
      & \text{on} \ \partial \Omega  \times (0,T]
      \\
       y(0) = y_0(\mu_1^0, \mu_2^0)
       \qquad
       & \text{in} \ \Omega \times \{t = 0\}
\end{cases}
\end{equation}
where $\Omega = (-1,1)^2$ is the 2D domain with space coordinates denoted by $x_1$ and $x_2$, $\nu$ is the diffusion coefficient -- here set equal to $0.001$ in order to focus mainly on the transport effect -- and $\mathbf{n}$ is the outward normal versor to the boundary $\partial \Omega$. Starting from a Gaussian density with variance equal to $0.05$ and centered at $(\mu_1^0, \mu_2^0)$, that is
\[
y_0(\mu_1^0, \mu_2^0) = \dfrac{10}{\pi} \exp{(- 10(x_1 - \mu_1^0)^2 - 10(x_2 - \mu_2^0)^2)},
\]
we aim to steer the state $y:\Omega \times [0,T] \to \R$ toward a final target destination exploiting the velocity field $\mathbf{u}:\Omega \times [0,T] \to \mathbb R^2$ as control action. Note that the parametrization of the initial state $y_0 = y_0(\mu_1^0, \mu_2^0)$ is crucial to reduce the state variability while exploring the state space in the data generation. By doing so, we focus on a specific but meaningful set of starting configurations, which correspond to steady particles systems at rest awaiting instructions.

\begin{remark}
Let $M(t) = \int_{\Omega} y(t) d\Omega$ be the total mass obtained by integrating the state variable in the domain of interest. As shown by, e.g., \cite{Sinigaglia2022a}, Equation~\eqref{eq:FP} entails mass conservation, as required when considering robotic swarms moving in $\Omega \times [0,T]$. Indeed, thanks to the divergence theorem and to the boundary condition selected, we obtain 
\[
\dfrac{d}{dt}M(t) = \int_{\Omega} \dfrac{\partial y(t)}{\partial t} d\Omega = \int_{\Omega} - \nabla \cdot (- \nu \nabla y + \mathbf{u} y) d\Omega = \int_{\partial \Omega} - (-\nu \nabla y + \mathbf{u} y) \cdot \mathbf{n} d \Gamma = 0
\]
\end{remark}

The domain $\Omega$ -- which does not show inner obstacles, for the sake of simplicity -- is discretized through \texttt{gmsh} utilities \citep{gmsh} yielding a conformal mesh with triangular elements and $7569$ nodes. The semi-discrete formulation of Equation~\eqref{eq:FP} is derived by FEM, taking into account continuous and piecewise linear finite element basis functions both for the state and the control, ending up with remarkably high number of degrees of freedom $N_h^y = 7569$ and $N_h^u = 15138$. To solve the high-fidelity FOM, time discretization is made over a uniform grid spanning $[0,T]$ with time step $\Delta t = 0.25$, where the final time $T$ is set equal to $1$.

\begin{figure}[!ht]
\centering
\subfloat{\includegraphics[scale = 0.23]{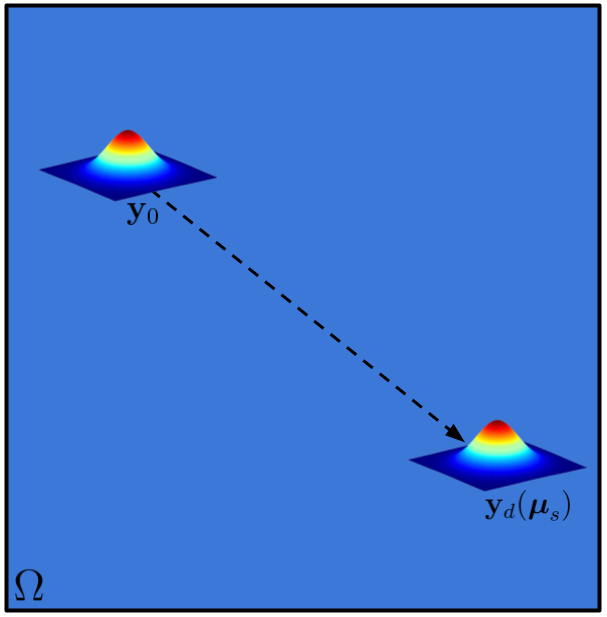}} \qquad \qquad \subfloat{\includegraphics[scale = 0.5]{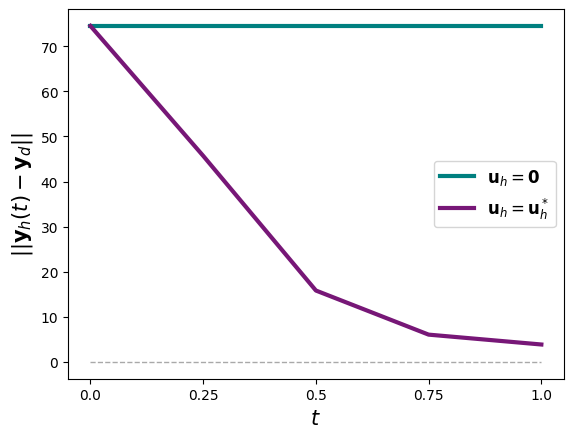}}

\subfloat{\includegraphics[height = 0.195\linewidth]{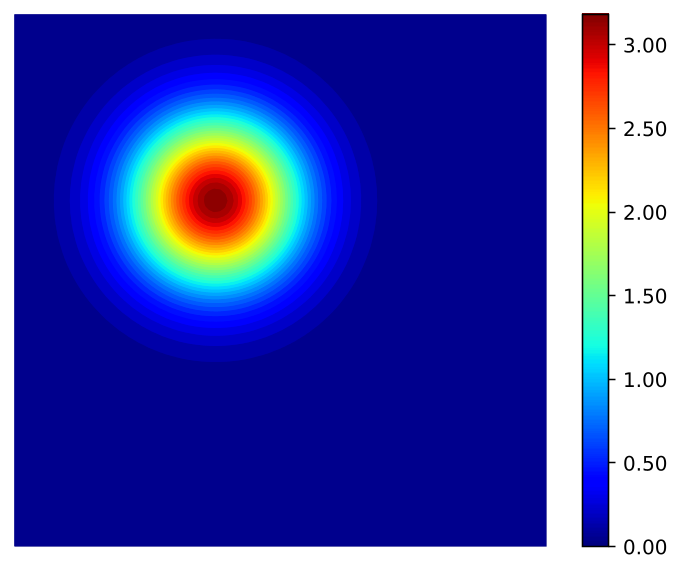}} \quad    \subfloat{\includegraphics[height = 0.195\linewidth]{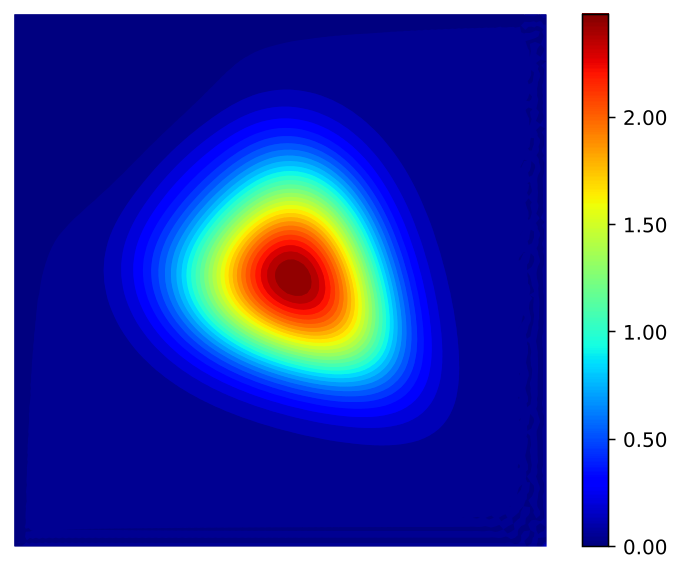}} \quad
\subfloat{\includegraphics[height = 0.195\linewidth]{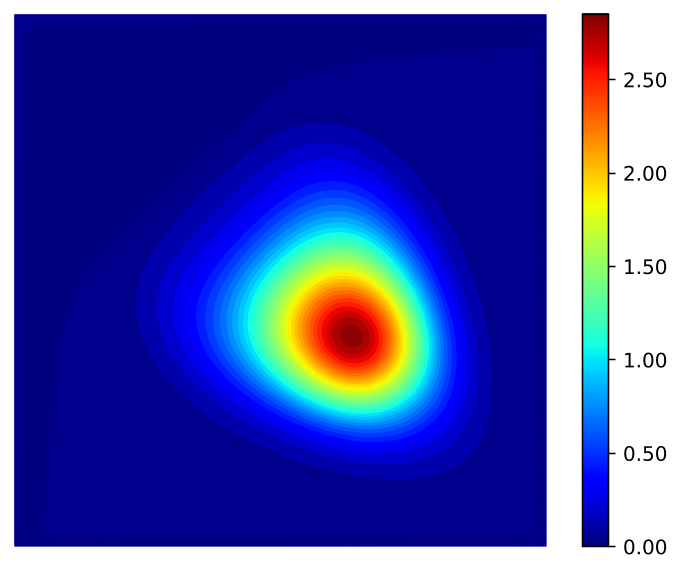}} \quad
\subfloat{\includegraphics[height = 0.195\linewidth]{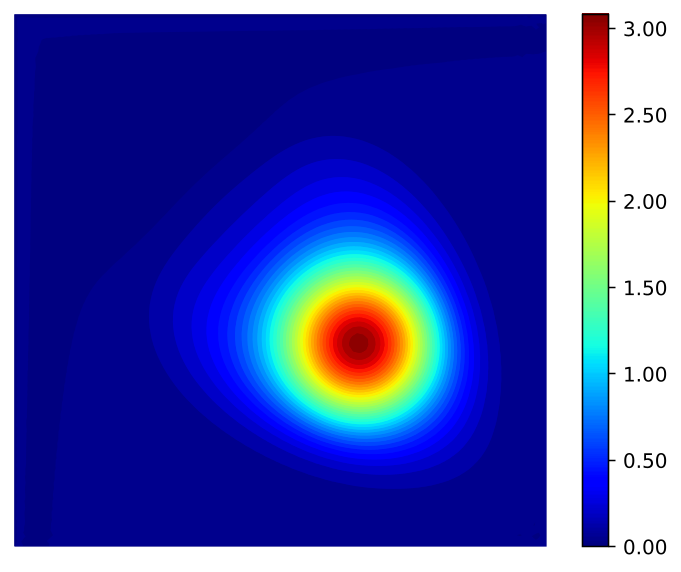}}
    
\subfloat{\includegraphics[height = 0.195\linewidth]{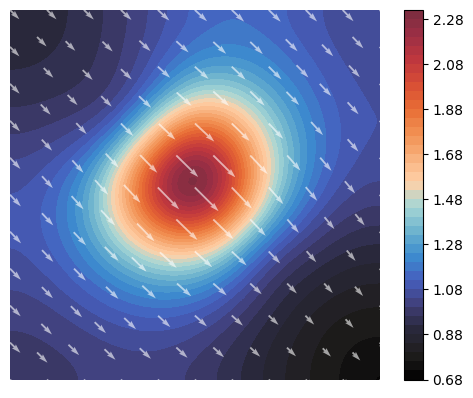}} \quad
\subfloat{\includegraphics[height = 0.195\linewidth]{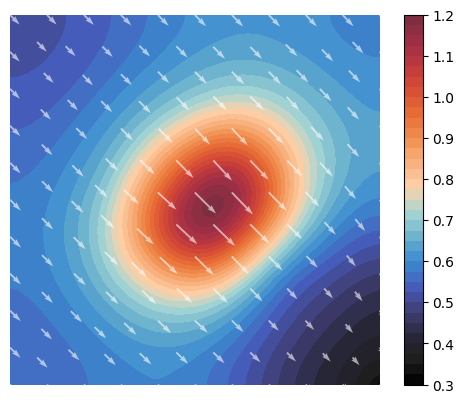}} \quad
\subfloat{\includegraphics[height = 0.195\linewidth]{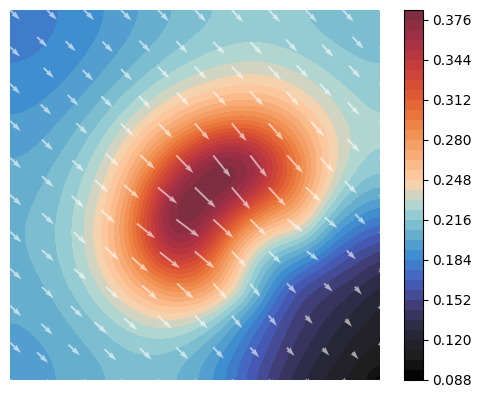}} \quad \subfloat{\includegraphics[height = 0.195\linewidth]{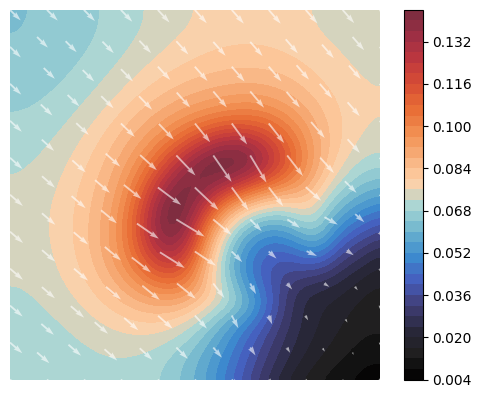}} 

\caption{\textit{Test 1.1}. Optimal transport in a vacuum. Top left: representation of an optimal state trajectory in a vacuum within the domain $\Omega$, where $\yinit$ stands for the initial density centered at $(\mu_1^0, \mu_2^0) =(-0.45, 0.21)$, while $\mathbf{y}_d(\mus)$ represents the target configuration centered at $(\mu_1^d, \mu_2^d) = (0.29, -0.24)$. Top right: discrepancy between current state $\yh(t)$ and target configuration $\mathbf{y}_d$ centered at $(\mu_1^d, \mu_2^d) = (0.29, -0.24)$ at different time instants in the uncontrolled ($\uh = \mathbf{0}$) and optimal ($\uh = \uh^*$) settings. Other panels: space-varying optimal state and control at $t=0,0.25,0.5,0.75$ related to the scenario parameters $\mus = (0.29, -0.24)$. The control velocity fields on $\Omega$ are depicted through vector fields, with the underlying colours corresponding to their magnitude.}
\label{fig:vacuum_setting}
\end{figure}

The optimal transport field capable of steering the state towards a target location in $\Omega$ can be found through an optimization procedure. In particular, among all the physically-admissible state-control pairs satisfying Equation~\eqref{eq:FP}, we aim to find the minimizer of the cost functional 
\begin{equation}
\label{eq:J}
J(y,\mathbf{u}; \mus) = \frac{1}{2} \int_0^T \int_{\Omega} (y - y_d)^2 d\Omega dt + \int_0^T \int_{\partial \Omega} y^2 d\Gamma dt + \frac{\beta}{2} \int_0^T \int_{\Omega} ||\mathbf{u}||^2 d\Omega dt + \frac{\beta_g}{2}  \int_0^T \int_{\Omega} ||\nabla \mathbf{u}||^2 d\Omega dt 
\end{equation}
measuring the discrepancy between the current observed state and the target state $y_d(\mu_1^d, \mu_2^d)$, which is the Gaussian density
\[
y_d(\mu_1^d, \mu_2^d) = \dfrac{10}{\pi} \exp{(- 10(x_1 - \mu_1^d)^2 - 10(x_2 - \mu_2^d)^2)}
\]
centered at $(\mu_1^d, \mu_2^d)$ and showing the same covariance as the initial configuration $y_0$. Since we aim to rapidly control the state dynamics considering different destinations online, the mean coordinates of the target location are regarded as scenario parameters, that is $\mus = (\mu_1^d, \mu_2^d)$. Instead, since the state information is directly captured within the feedback loop, the mean coordinates of the initial state $(\mu_1^0, \mu_2^0)$ are not regarded as scenario parameters, so that their knowledge is not required in the online phase. The boundary integral in Equation~\eqref{eq:J} is useful to avoid collisions with the domain boundary, as necessary when dealing with particles systems, while the regularizing terms concerning the norms of the control and its gradient penalize overconsuming and irregular control actions, which may be unfeasible in practice. To properly balance the magnitudes of the different terms appearing in the loss functional, we set $\beta = \beta_g = 0.2$. Note that the discrete cost functional $J_h = J_h(\yh, \uh; \mus)$ introduced in Equation~\eqref{eq:OCP} can be easily recovered starting from Equation~\eqref{eq:J} thanks to discretization techniques such as FEM. To further visualize the problem setting, Figure~\ref{fig:vacuum_setting} shows an example of optimal state trajectory related to $(\mu_1^0, \mu_2^0) =(-0.45, 0.21)$ and $\mus = (0.29, -0.24)$, along with the corresponding optimal control actions applied. Specifically, it is possible to assess that the optimal transport steers the state density towards the target destination over time, that is
\[||\yh(t) - \mathbf{y}_d|| \underset{t \to T}{\longrightarrow} 0\]

As described in Section~\ref{sec:OCP-DL-ROM}, the starting point in building the deep learning-based reduced order feedback controller is data generation. In particular, we consider $100$ random scenarios sampled in the parameter space $\mathcal{P} = (0.0, 0.5) \times (-0.5, 0.5)$, i.e. we take into account different endpoints placed on the right-hand side of the domain, avoiding positions too close to the boundaries. Moreover, to account for the initial state variability, we sample the starting position coordinates on the left-hand side of $\Omega$, that is $(\mu_1^0, \mu_2^0) \in (-0.5, 0.0) \times (-0.5, 0.5)$.  For every combination of starting-final state positions, we thus compute the optimal state and control trajectory at $N_t = \frac{T}{\Delta t} = 4$ time instants through \texttt{dolfin-adjoint} \citep{Mitusch2019}, which provides an OCP solver in \texttt{fenics} \citep{fenics} exploiting FEM and the adjoint method, with an average computational time equal to $15$ minutes per scenario. In particular, we select L-BFGS-B as optimization algorithm, while considering a tolerance and maximum number of iterations equal to, respectively, $10^{-6}$ and $500$. Note that, leveraging the horizontal symmetry of the problem, it is possible to double the dataset cardinality for free, ending up with a total of $N_s = 200$ scenarios and $N_s N_t = 800$ state-scenario-control triplets investigated. Out of all the available trajectories, training and test sets are obtained through a $80:20$ split. In particular, while training data are exploited for reduction and neural networks training purposes, test data are used only to evaluate the generalization capabilities of our controller.   

\begin{figure}[!ht]
\centering
\begin{minipage}{0.44\linewidth}
\centering
\includegraphics[scale = 0.5]{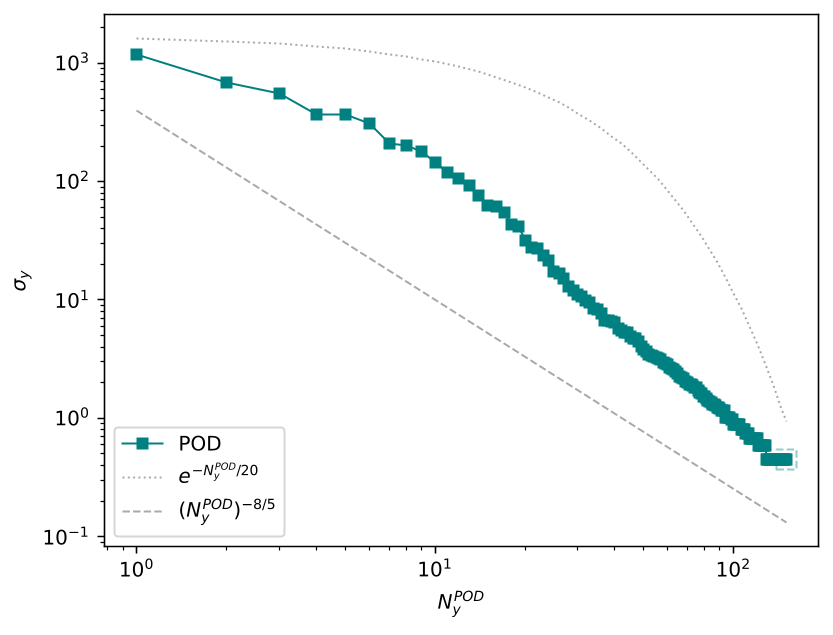}
\end{minipage}
\begin{minipage}{0.55\linewidth}
\centering
\hfill
\vspace{1.05cm}
    \subfloat{\includegraphics[height = 0.38\linewidth]{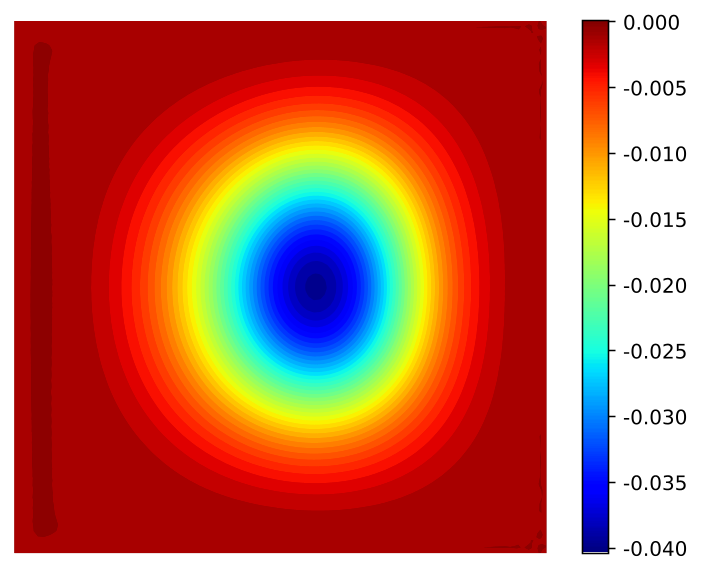}} \,
    \subfloat{\includegraphics[height = 0.38\linewidth]{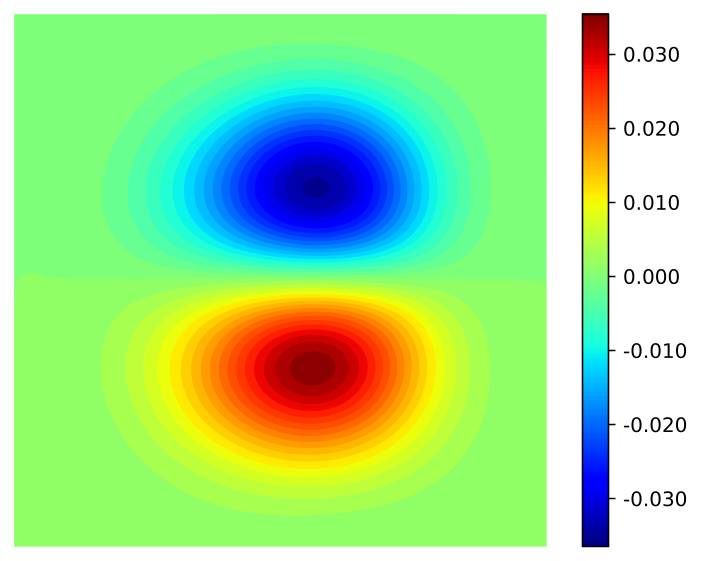}}
\end{minipage}

\begin{minipage}[c]{0.44\linewidth}
\centering
    \includegraphics[scale = 0.5]{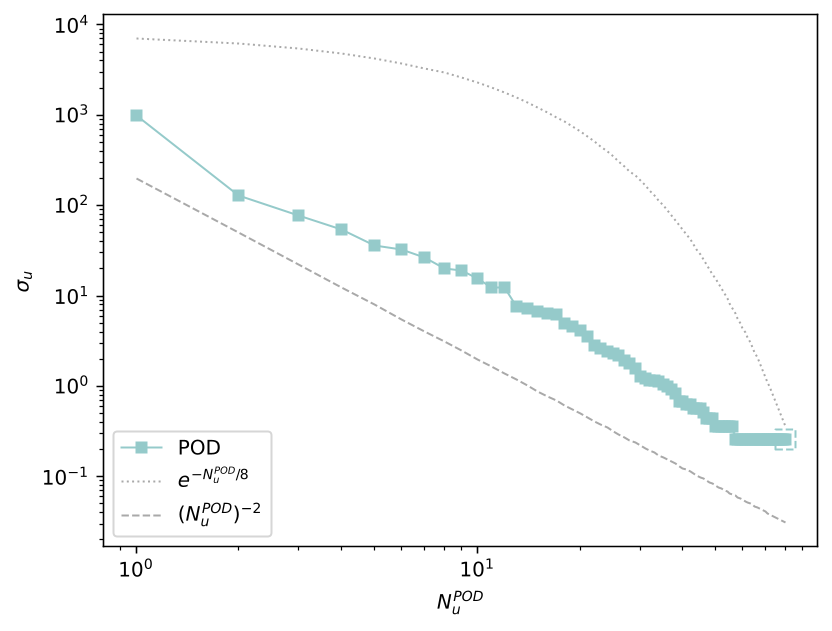}
\end{minipage}
\begin{minipage}[c]{0.55\linewidth}
\centering
\hfill
\vspace{1.05cm}
    \subfloat{\includegraphics[height = 0.38\linewidth]{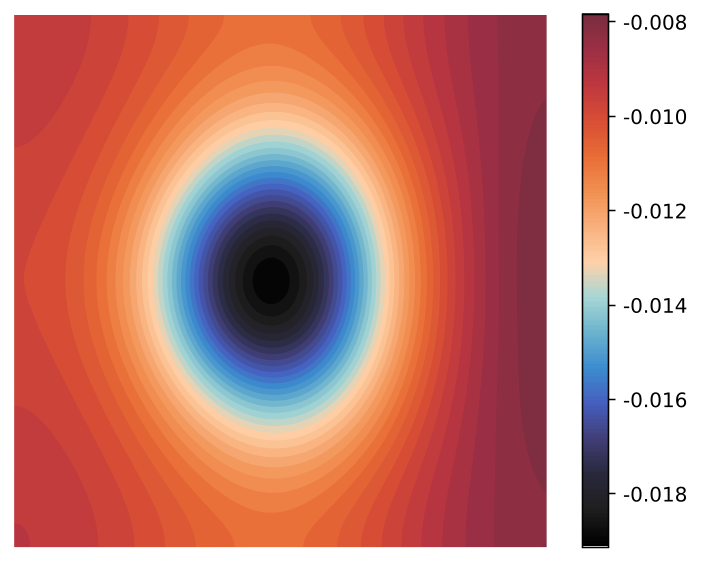}} \,
    \subfloat{\includegraphics[height = 0.38\linewidth]{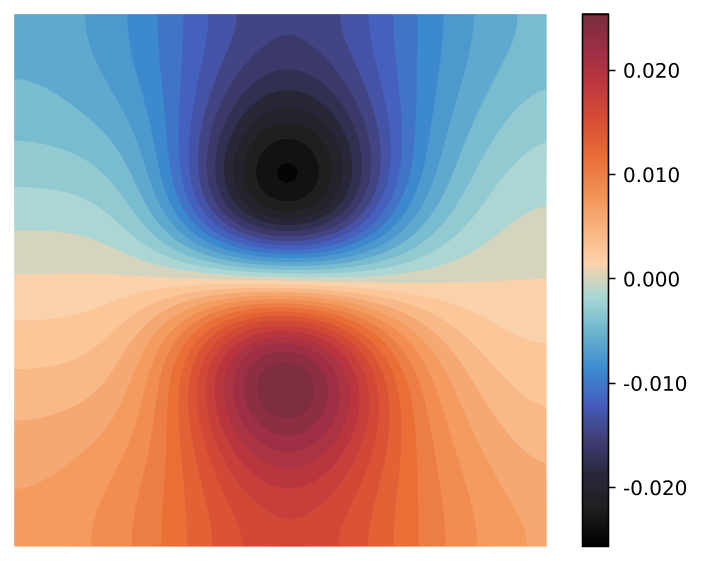}}
\end{minipage}

\begin{minipage}[c]{0.44\linewidth}
\centering
    \includegraphics[scale = 0.5]{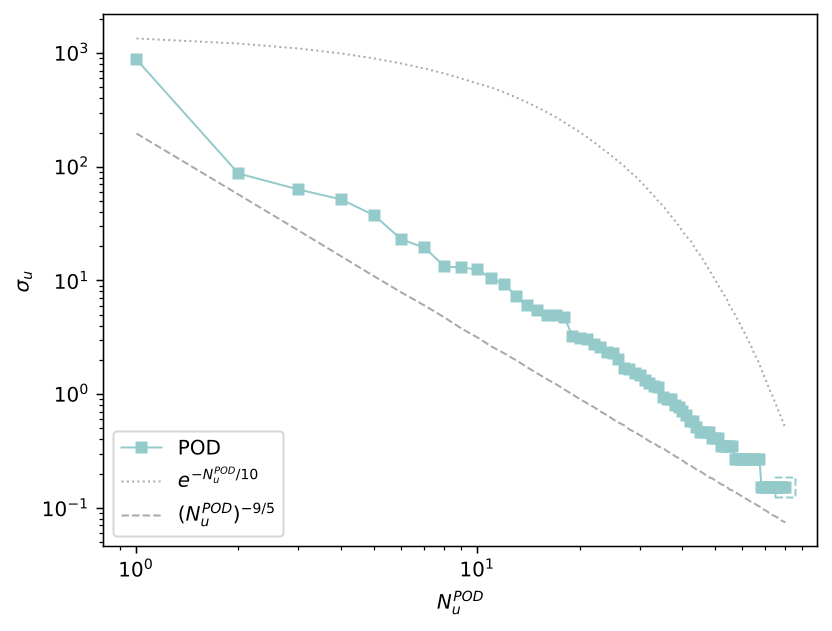}
\end{minipage}
\begin{minipage}[c]{0.55\linewidth}
\centering
\hfill
\vspace{1.05cm}
    \subfloat{\includegraphics[height = 0.38\linewidth]{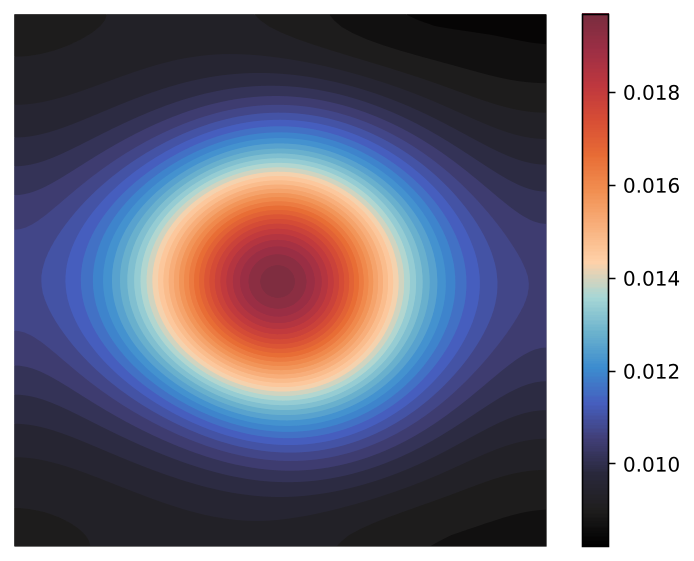}} \,
    \subfloat{\includegraphics[height = 0.38\linewidth]{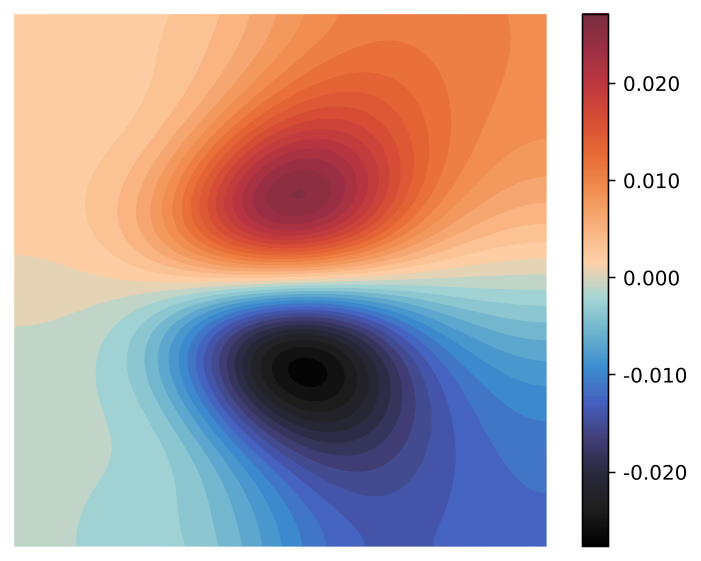}}
\end{minipage}

\caption{\textit{Test 1.1}. Optimal transport in a vacuum. Singular values decay in log-log scale along with the two most energetic POD modes related to the state (top), $x_1$ component (center) and $x_2$ component (bottom) of the control.}
\label{fig:PODvacuum}
\end{figure}

The second necessary step to construct our feedback controller is dimensionality reduction. First of all, we apply a linear projection through POD, looking at the singular values decays and the reconstruction errors in order to select the number of modes to retain. In particular, by selecting $150$ and $160$ modes for state and control, respectively, we end up with $\varepsilon^y_{\mathrm{rel}} = 0.21\%$ and $\varepsilon^u_{\mathrm{rel}} = 0.36\%$. Note that, being the control variable a vector field, POD is applied component-wise, thus taking into account $80$ POD modes for each spatial component. Figure~\ref{fig:PODvacuum} displays the singular values decays and the two most energetic POD modes related to state and control reductions. In particular, due to the significant transport effect, the singular values show slow polynomial decays, thus requiring a remarkably high number of POD modes in order to correctly reconstruct the state and control fields.

\begin{figure}
\centering
\subfloat{\includegraphics[height = 0.195\linewidth]{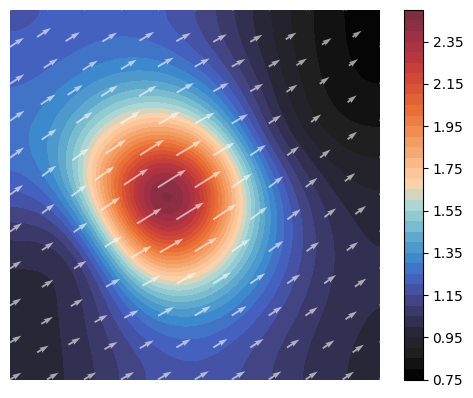}} \quad
\subfloat{\includegraphics[height = 0.195\linewidth]{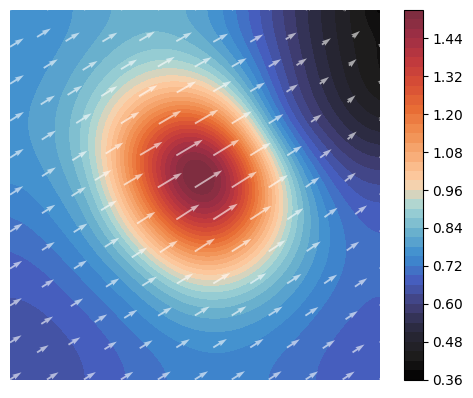}} \quad
\subfloat{\includegraphics[height = 0.195\linewidth]{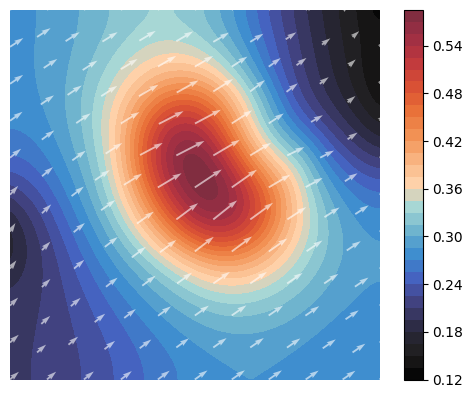}} \quad
\subfloat{\includegraphics[height = 0.195\linewidth]{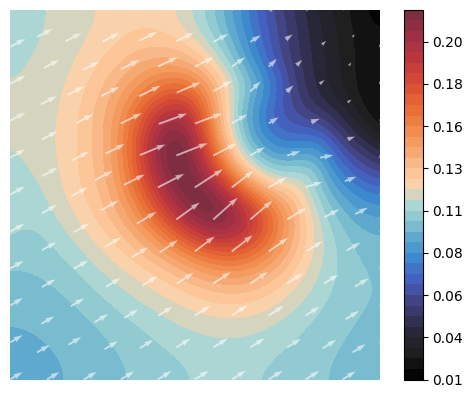}}

\subfloat{\includegraphics[height = 0.195\linewidth]{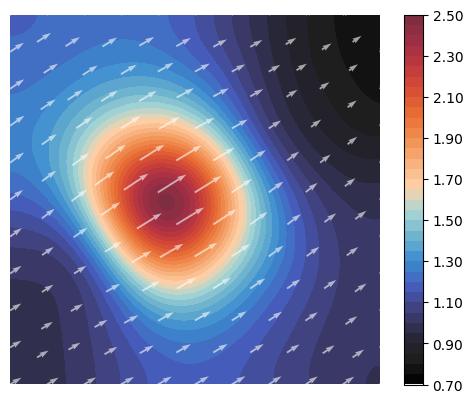}} \quad
\subfloat{\includegraphics[height = 0.195\linewidth]{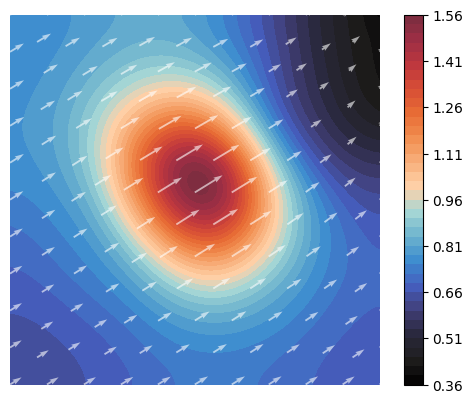}} \quad
\subfloat{\includegraphics[height = 0.195\linewidth]{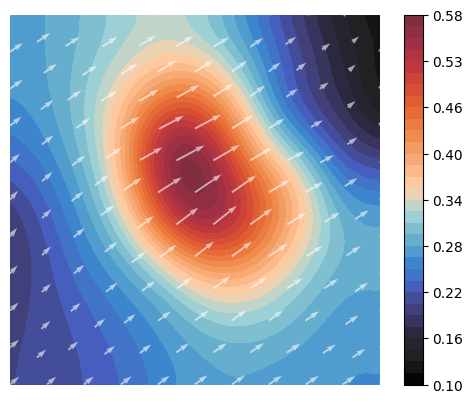}} \quad
\subfloat{\includegraphics[height = 0.195\linewidth]{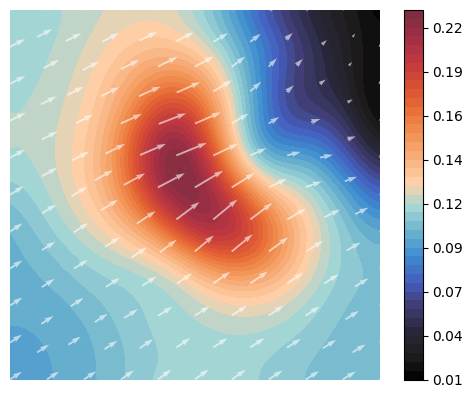}}

\subfloat{\includegraphics[height = 0.195\linewidth]{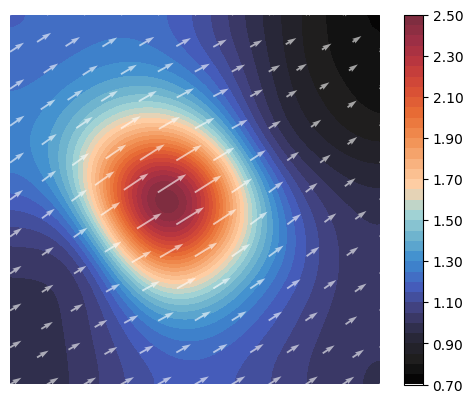}} \quad
\subfloat{\includegraphics[height = 0.195\linewidth]{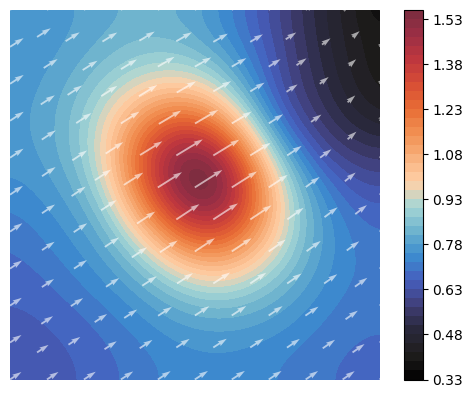}} \quad
\subfloat{\includegraphics[height = 0.195\linewidth]{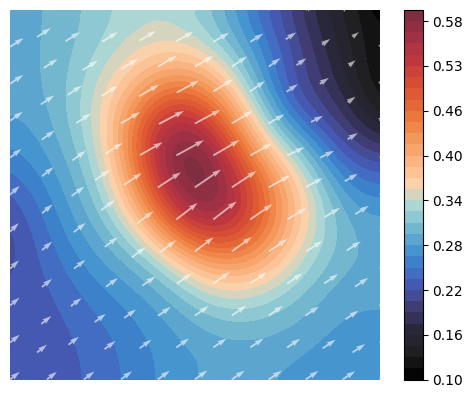}} \quad
\subfloat{\includegraphics[height = 0.195\linewidth]{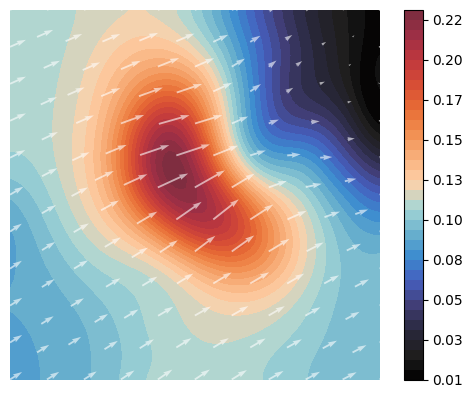}}

\subfloat{\includegraphics[height = 0.195\linewidth]{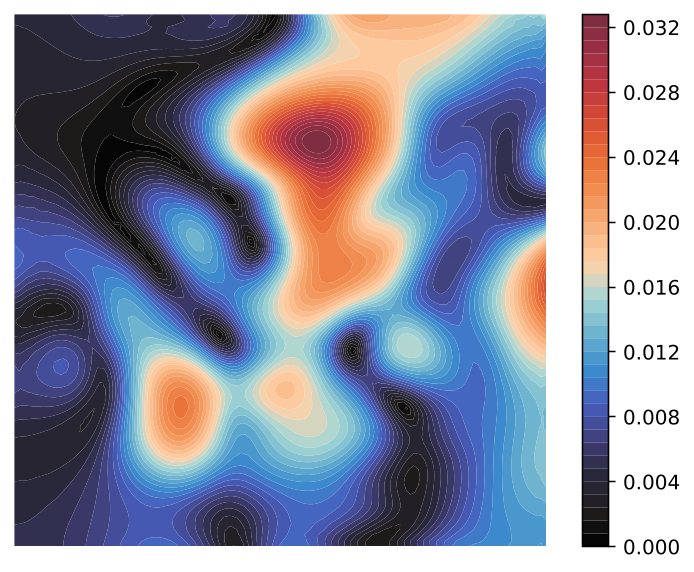}} \hspace{0.55em}
\subfloat{\includegraphics[height = 0.195\linewidth]{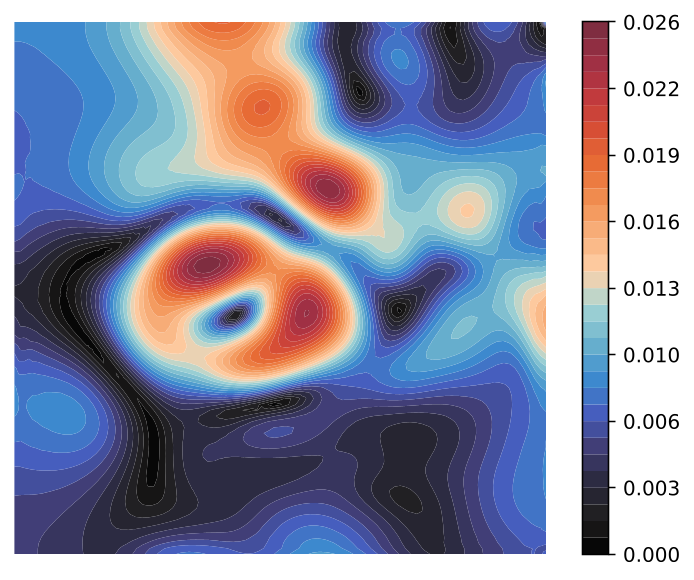}} \hspace{0.55em}
 \subfloat{\includegraphics[height = 0.195\linewidth]{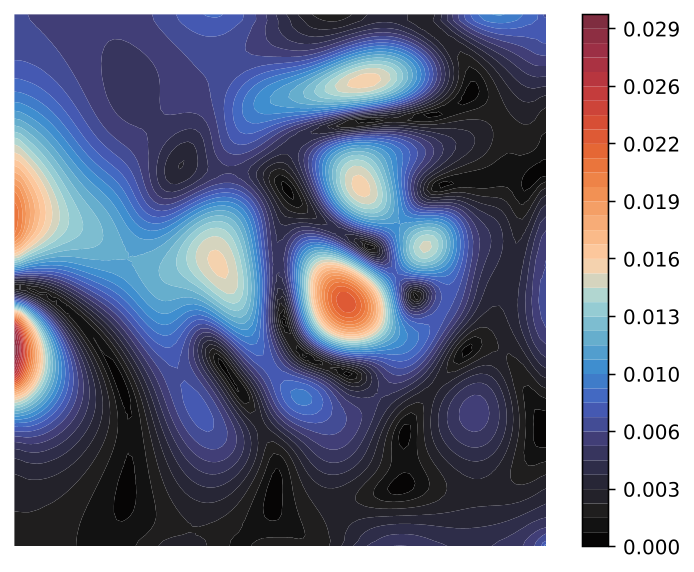}} \hspace{0.55em}
\subfloat{\includegraphics[height = 0.195\linewidth]{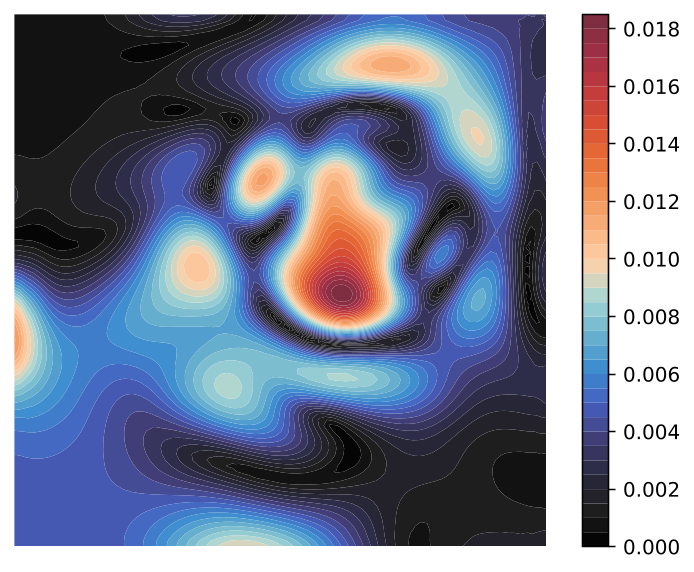}}

\subfloat{\includegraphics[height = 0.195\linewidth]{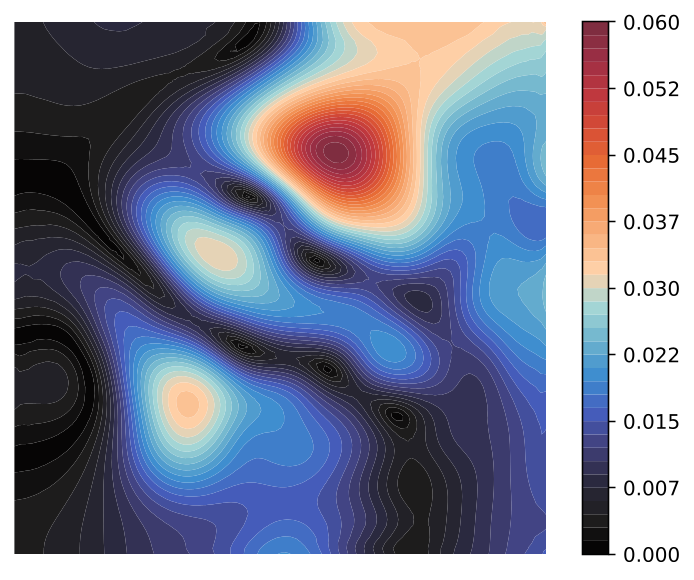}} \hspace{0.55em}
\subfloat{\includegraphics[height = 0.195\linewidth]{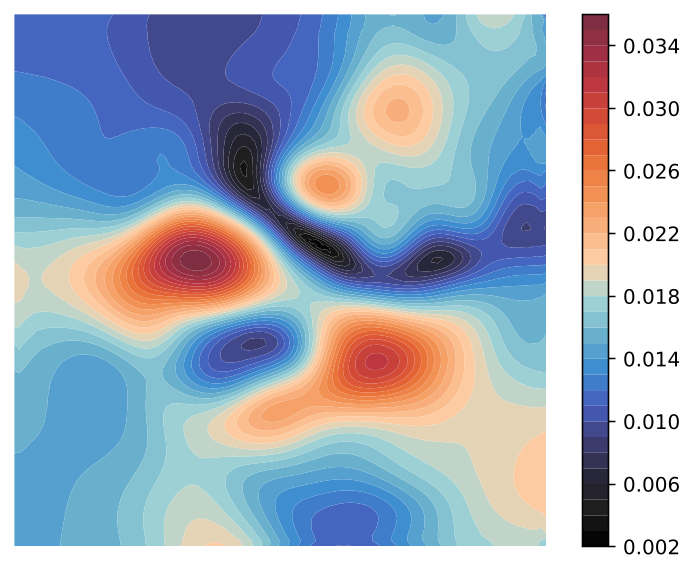}} \hspace{0.55em}
\subfloat{\includegraphics[height = 0.195\linewidth]{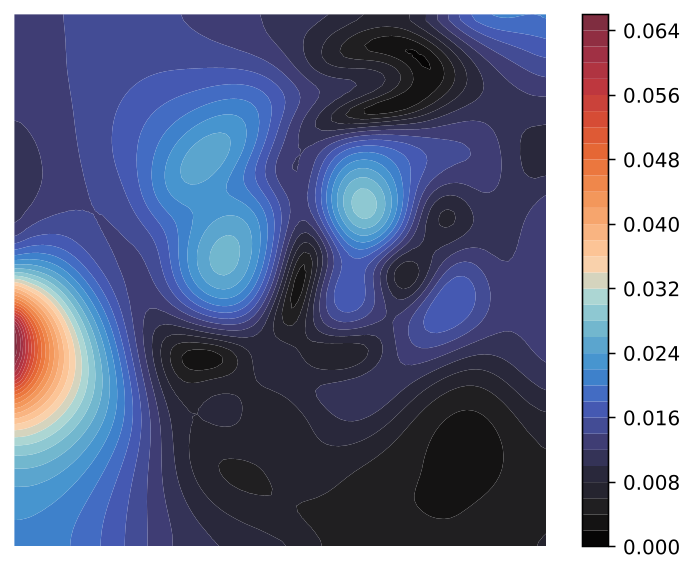}} \hspace{0.55em}
\subfloat{\includegraphics[height = 0.195\linewidth]{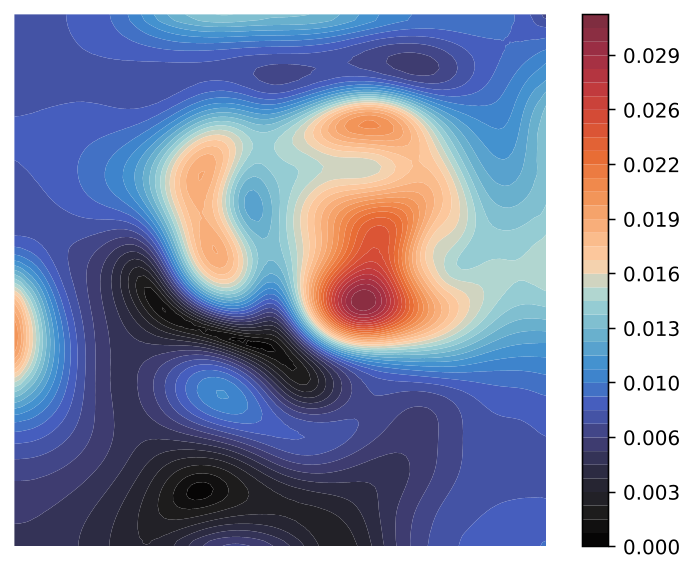}}

\caption{\textit{Test 1.1}. Optimal transport in a vacuum. High-fidelity optimal control trajectory (first row), POD+AE reconstructions (second row), policy predictions (third row), POD+AE reconstruction errors (fourth row) and policy prediction errors (fifth row) at $t = 0, 0.25, 0.5, 0.75$ related to the test scenario parameters $\mus = (0.31, 0.27)$. The control velocity fields on $\Omega$ are depicted through vector fields, with the underlying colours corresponding to their magnitude.}
\label{fig:vacuum_control}
\end{figure}

The state and control POD coefficients can be further compressed through an autoencoder-based projection onto nonlinear subspaces, resulting in the POD+AE reduction strategy introduced in Section~\ref{sec:OCP-DL-ROM}. As far as the state autoencoder architecture is concerned, the latent dimension is set equal to $N_y = 10$, while the encoder and decoder consist, respectively, of $1$ and $2$ hidden layers with $100$ neurons each and with leaky Relu as activation function. A similar structure is taken into account for the control reduction, where the latent dimension is increased to $N_u = 18$ and the number of neurons per layer in the decoder is doubled. In order to predict the optimal control action starting from an observed state configuration, a policy surrogate model $\pi_N$ is required. In particular, to this aim, we take into account a deep feedforward neural network showing $3$ hidden layers with $50$ neurons each and leaky Relu as activation function. Note that no activation functions are considered in the output layer of the networks to avoid restrictions on output values. After initializing the weights through the strategy proposed by \cite{He2015}, we train the networks in $1$ hour and $11$ minutes minimizing the cumulative loss function $J_{\mathrm{NN}}$ introduced in Section~\ref{sec:OCP-DL-ROM} through the L-BFGS optimization algorithm, while considering $\lambda_1 = \lambda_2 = \lambda_3 = 0.01$. The POD+AE reconstruction errors entailed by the use of the state and control autoencoders end up being, respectively, $\varepsilon^y_{\mathrm{rel}} = 3.20\%$ and $\varepsilon^u_{\mathrm{rel}} = 5.04\%$. Instead, the relative prediction error of the policy surrogate model on the test data at the latent level -- that is, the discrepancy between $\un$ and $\tildeun = \pi_N(\yn, \mus)$ -- is equal to $4.28\%$, while it increases to $7.09\%$ after decoding through POD+AE -- that is, the error between $\uh$ and $\tildeuh = \mathbb{V}_u \varphi_D^u(\tildeun)$. Figure~\ref{fig:vacuum_control} displays a control test trajectory reconstructed by POD+AE (second row), that is
\[
\uhrec = \mathbb{V}_u \varphi_D^u(\varphi_E^u(\mathbb{V}_u^{\top} \uh))
\]
and the corresponding reconstruction errors (fourth row). In particular, we highlight that the latent coordinates, whose dimensionality is $841$ times smaller than $N_h^u$, perfectly capture the full-order features of the control velocity. By a visual inspection of the third row of Figure~\ref{fig:vacuum_control}, along with the reported errors in the last row of the same figure, it is also possible to assess the high accuracy of the optimal control predictions
\[
\tildeuh = \mathbb{V}_u \varphi_D^u(\pi_N(\yn, \mus))
\]
provided by the low-dimensional policy.

\begin{figure}
\centering
\subfloat{\includegraphics[height = 0.195\linewidth]{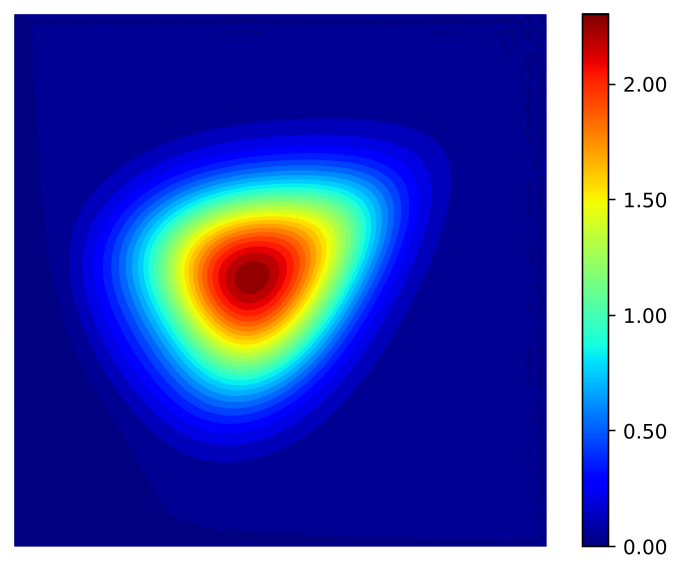}} \quad
\subfloat{\includegraphics[height = 0.195\linewidth]{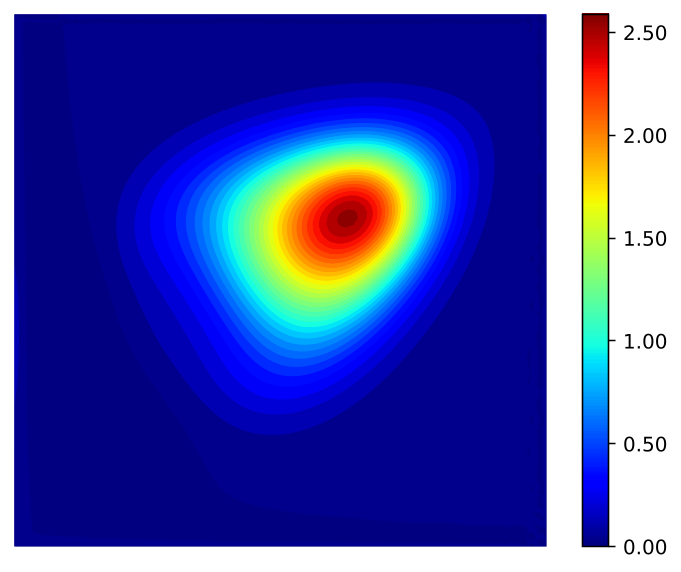}} \quad
\subfloat{\includegraphics[height = 0.195\linewidth]{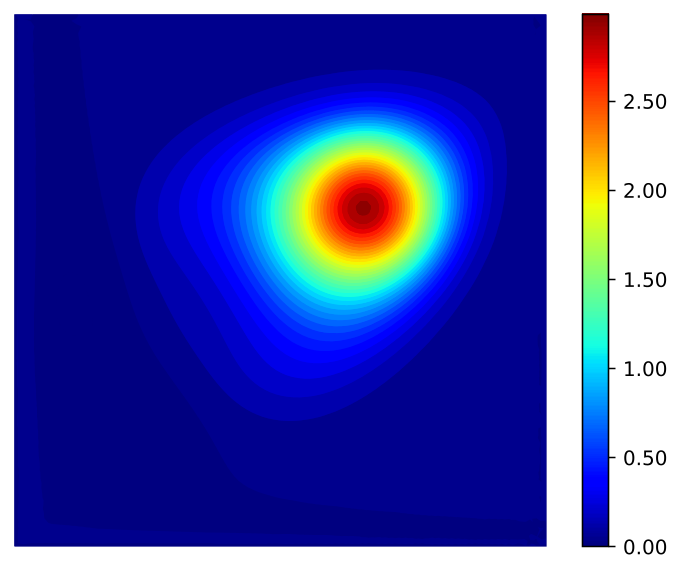}} \quad
\subfloat{\includegraphics[height = 0.195\linewidth]{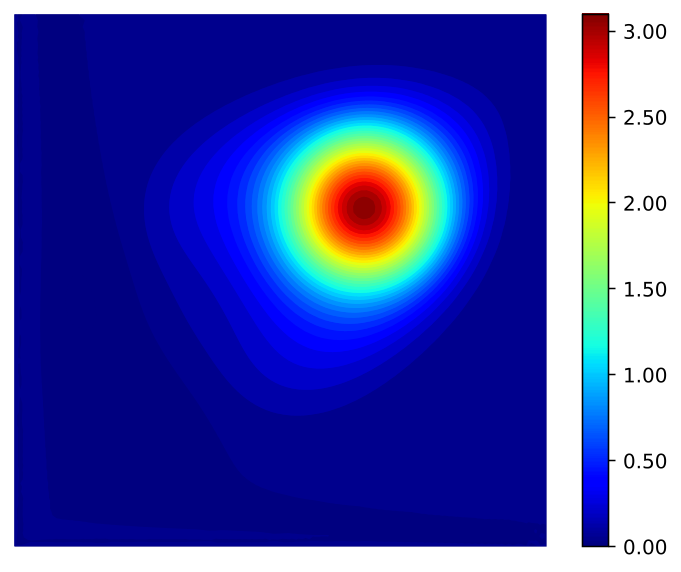}}

\subfloat{\includegraphics[height = 0.195\linewidth]{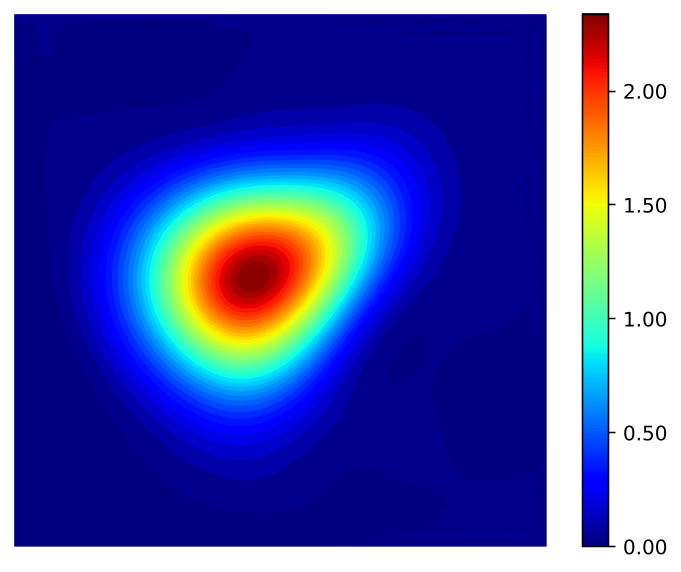}} \quad
\subfloat{\includegraphics[height = 0.195\linewidth]{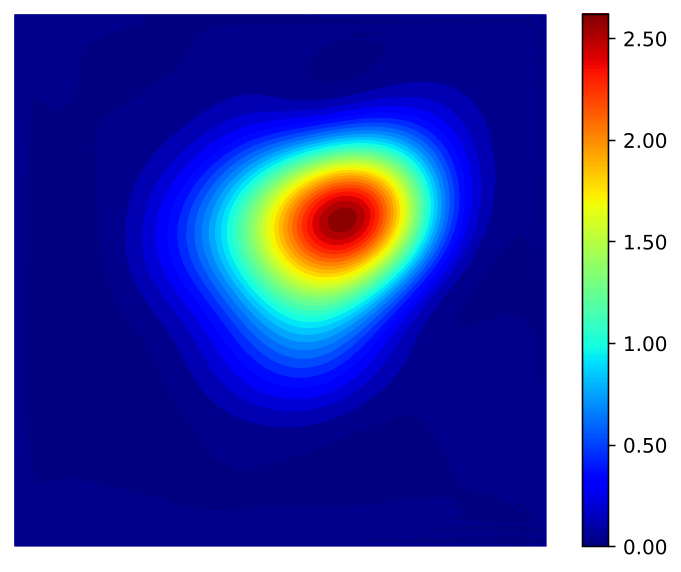}} \quad
\subfloat{\includegraphics[height = 0.195\linewidth]{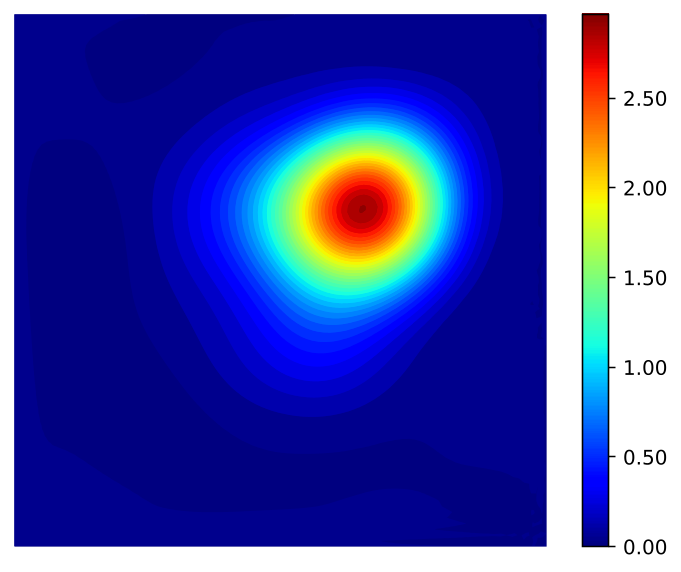}} \quad
\subfloat{\includegraphics[height = 0.195\linewidth]{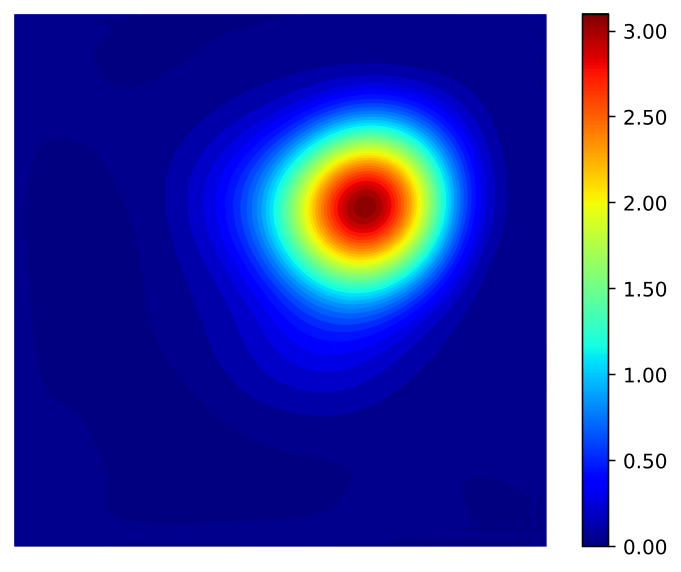}}

\subfloat{\includegraphics[height = 0.195\linewidth]{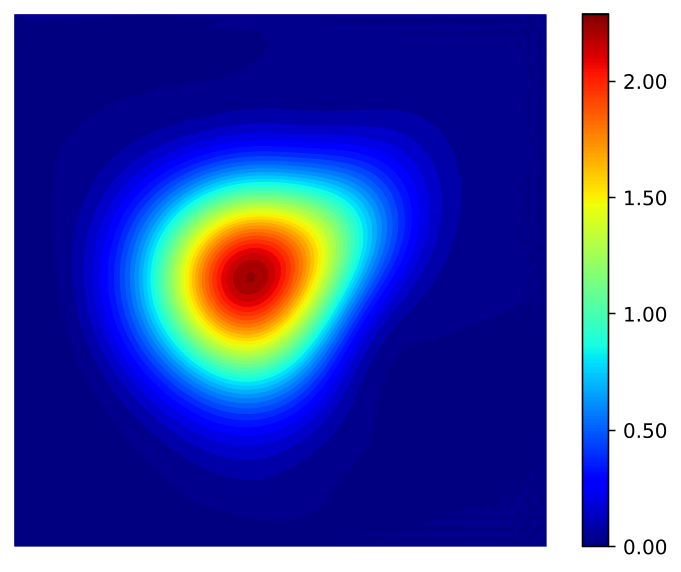}} \quad
\subfloat{\includegraphics[height = 0.195\linewidth]{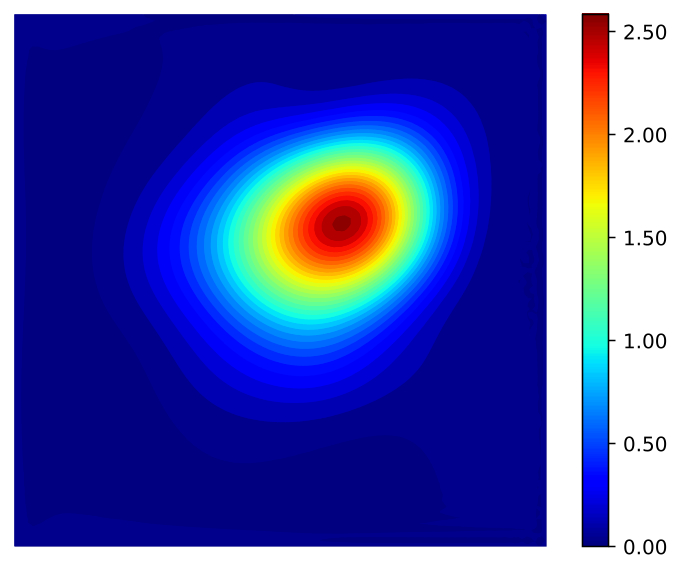}} \quad
\subfloat{\includegraphics[height = 0.195\linewidth]{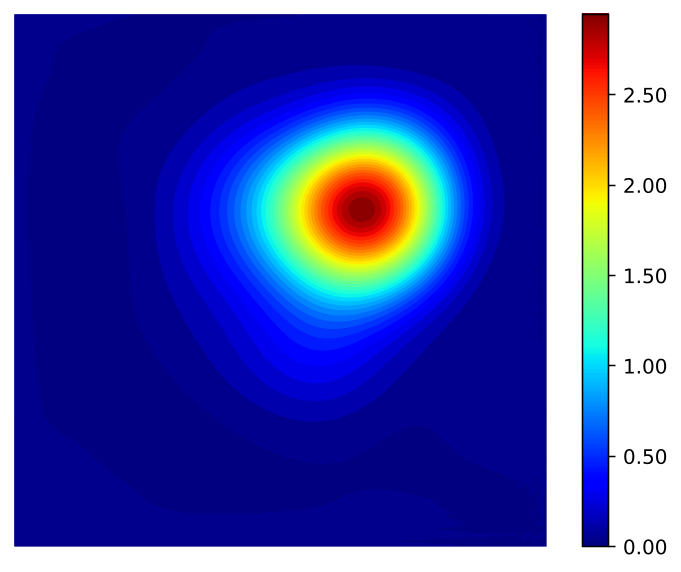}} \quad
\subfloat{\includegraphics[height = 0.195\linewidth]{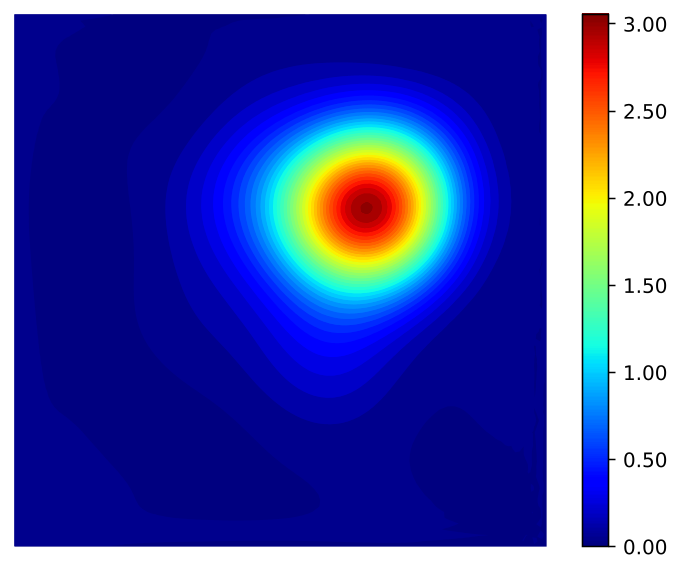}}

\subfloat{\includegraphics[height = 0.195\linewidth]{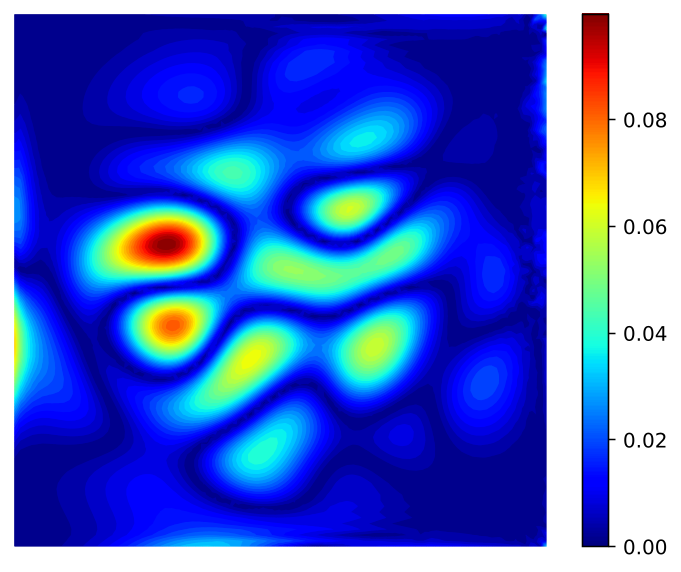}} \quad
\subfloat{\includegraphics[height = 0.195\linewidth]{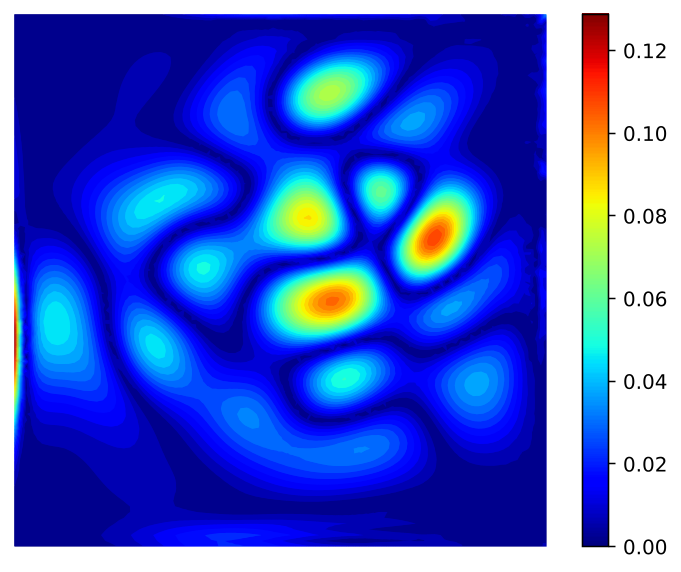}} \quad
\subfloat{\includegraphics[height = 0.195\linewidth]{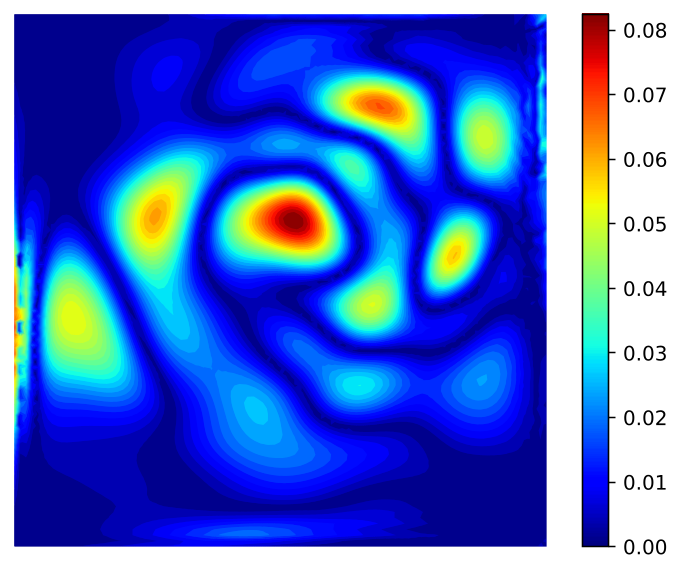}} \quad
\subfloat{\includegraphics[height = 0.195\linewidth]{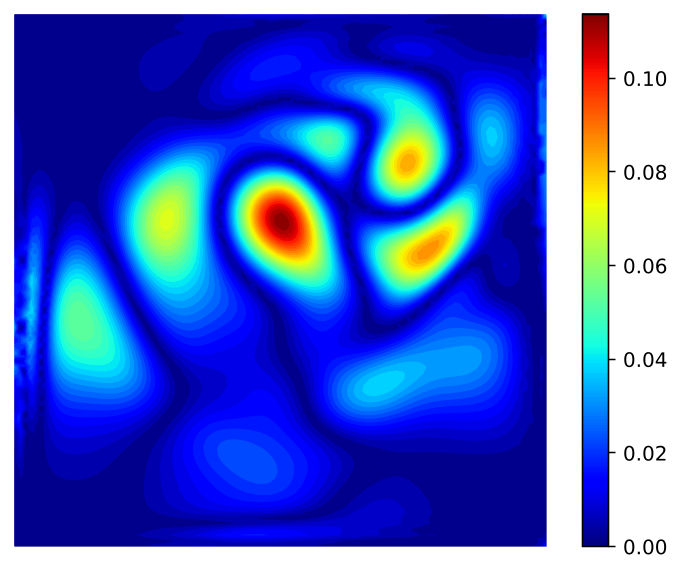}}

\subfloat{\includegraphics[height = 0.195\linewidth]{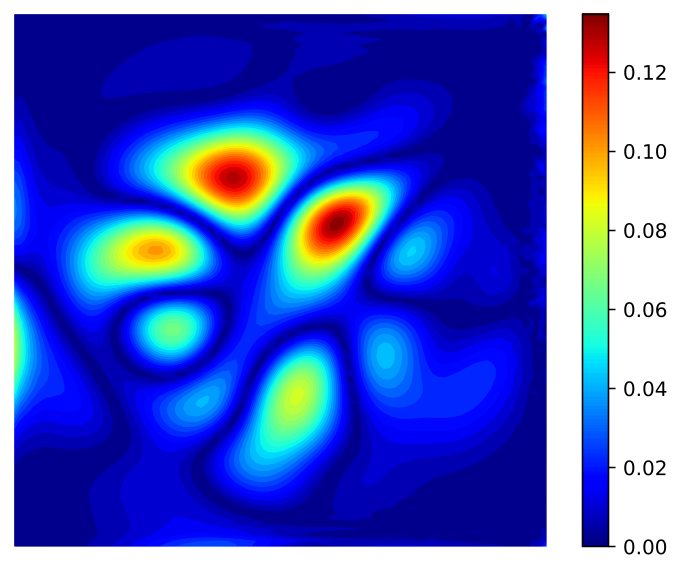}} \quad
\subfloat{\includegraphics[height = 0.195\linewidth]{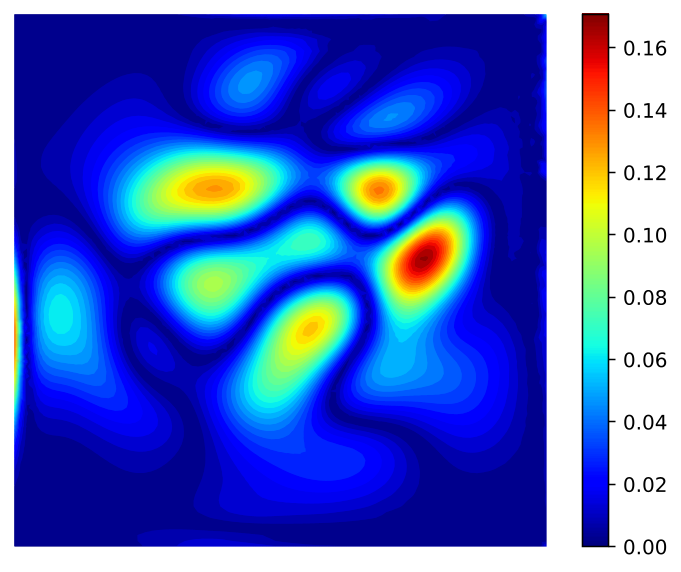}} \quad
\subfloat{\includegraphics[height = 0.195\linewidth]{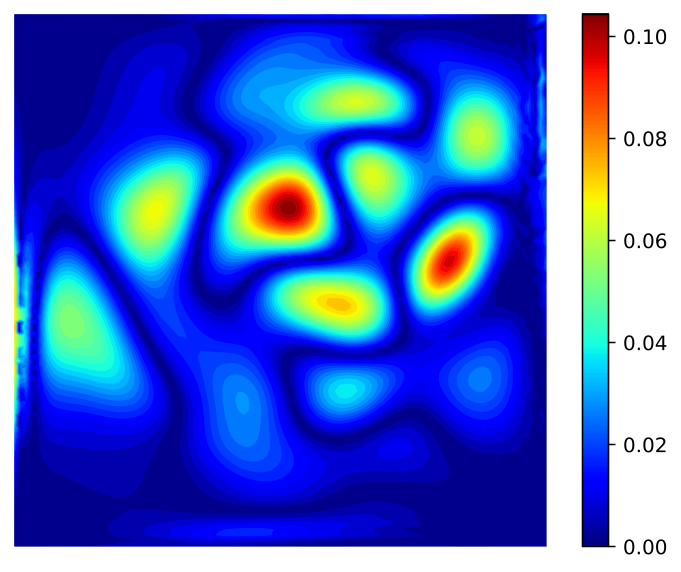}} \quad
\subfloat{\includegraphics[height = 0.195\linewidth]{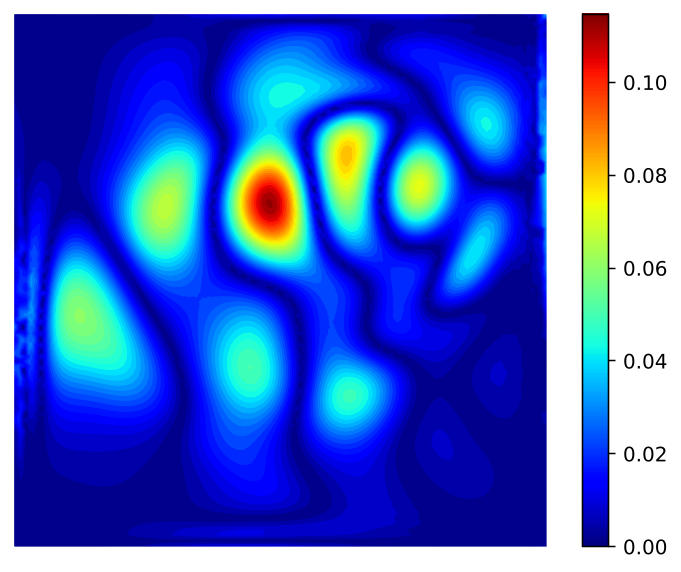}}

\caption{\textit{Test 1.1}. Optimal transport in a vacuum. High-fidelity optimal state trajectory (first row), POD+AE reconstructions (second row), forward model predictions (third row), POD+AE reconstruction errors (fourth row) and forward model prediction errors (fifth row) at $t = 0.25, 0.5, 0.75, 1.0$ related to the test scenario parameters $\mus = (0.31, 0.27)$.}
\label{fig:vacuum_state}
\end{figure}

In case continuous monitoring of the state density is unfeasible or the computational burden to simulate synthetic data does not meet strict timing requirements, we can exploit the latent feedback loop introduced in Section~\ref{subsec:latentfeedbackloop}. In particular, the surrogate model $\varphi_N$ for the dynamics at the latent level is retrieved by a deep feedforward neural network that has $3$ hidden layers with $50$ neurons each and takes into account leaky Relu as activation function. We train the neural networks --  that are the state and control autoencoders, the policy $\pi_N$ and the forward model $\varphi_N$ -- in $1$ hour and $20$ minutes by minimizing the cumulative loss function introduced in Section~\ref{subsec:latentfeedbackloop} with $\lambda_1 = \lambda_2 = \lambda_3 = \lambda_6 = 0.001$ and $\lambda_4 = \lambda_5 = 1$, exploiting L-BFGS as optimization algorithm. Note that, while $\varphi_N$ is initialized through the strategy proposed in \cite{He2015}, we exploit the previously trained networks as initializations for the autoencoders and the policy. The POD+AE reconstruction errors committed by the state and control autoencoders are now equal to $\varepsilon^y_{\mathrm{rel}} = 3.96\%$ and $\varepsilon^u_{\mathrm{rel}} = 4.61\%$, while the policy prediction error after-decoding is equal to $6.98\%$. Instead, the prediction-from-data and the prediction-from-policy relative errors entailed by $\varphi_N$ on test data are equal to, respectively, $2.09 \%$ and $1.37 \%$ at the latent level, while they increase to $7.49\%$ and $5.45 \%$ after POD+AE decoding. Figure~\ref{fig:vacuum_state} shows a state test trajectory reconstructed by POD+AE (second row), that is
\[
\yhrec = \mathbb{V}_y \varphi_D^y(\varphi_E^y(\mathbb{V}_y^{\top} \yh))
\]
and the corresponding reconstruction errors (fourth row). In particular, despite a dimensionality reduction of $757$ times with respect to $N_h^y$, the low-dimensional features are able to reconstruct accurately the full-order snapshots. Moreover, Figure~\ref{fig:vacuum_state} displays the forward model predictions (third row) 
\[
\tildeyh(t_j) = \mathbb{V}_y \varphi_D^y(\varphi_N(\yn(t_{j-1}), \un(t_{j-1}), \mus)) \quad \forall j=1,...,N_t
\]
and the corresponding prediction errors (last row) related to the same test scenario.

\begin{figure}
\centering
\subfloat{\includegraphics[scale = 0.2]{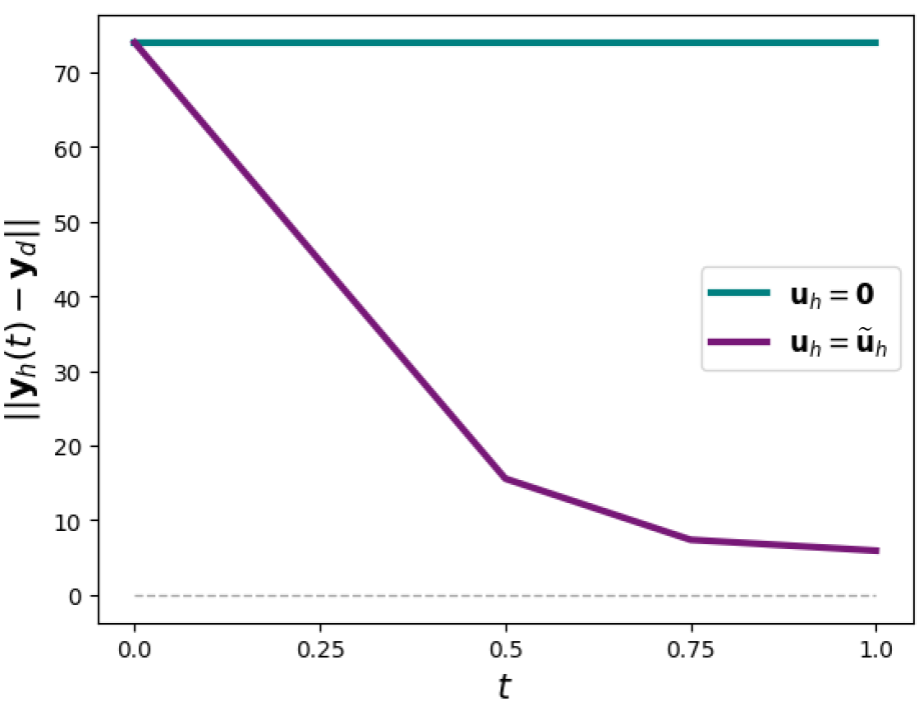}} \qquad \subfloat{\includegraphics[scale = 0.2]{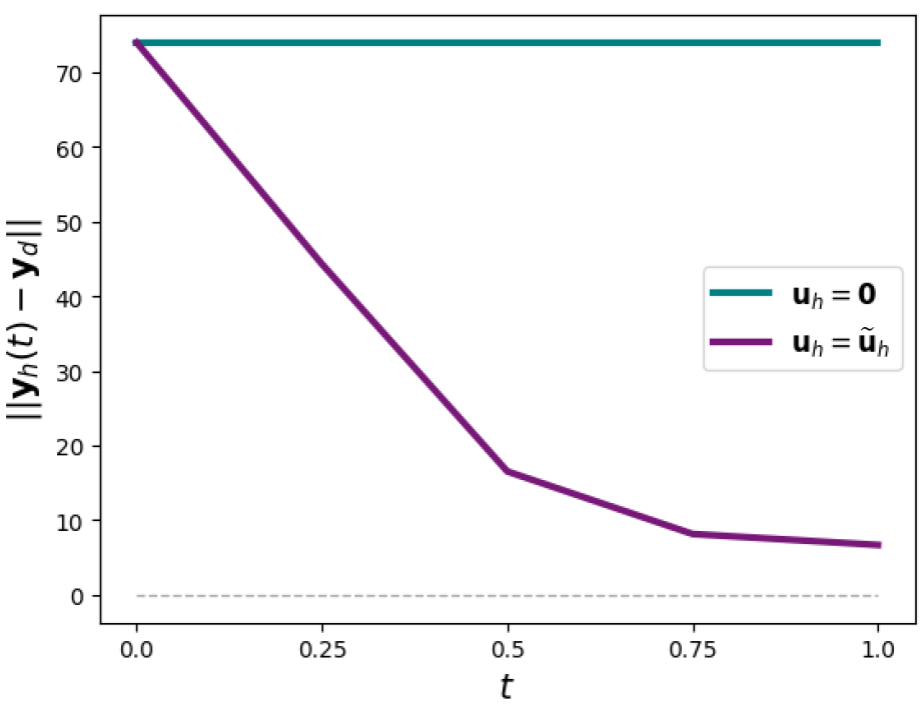}}\vspace{-1em}

\subfloat{\includegraphics[height = 0.195\linewidth]{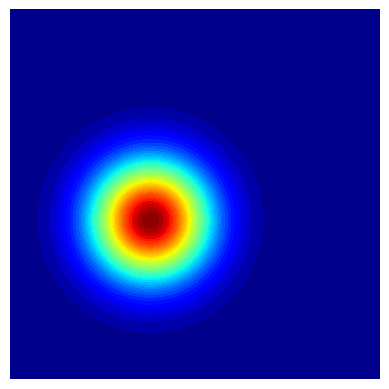}} 
\subfloat{\includegraphics[height = 0.195\linewidth]{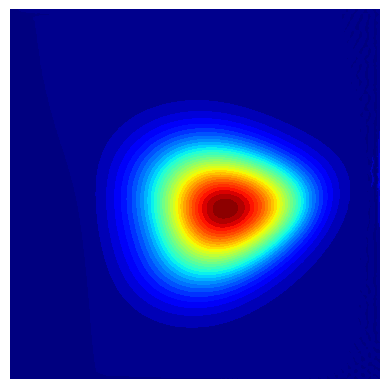}}
\subfloat{\includegraphics[height = 0.195\linewidth]{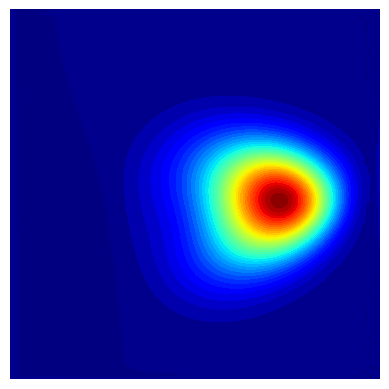}} 
\subfloat{\includegraphics[height = 0.195\linewidth]{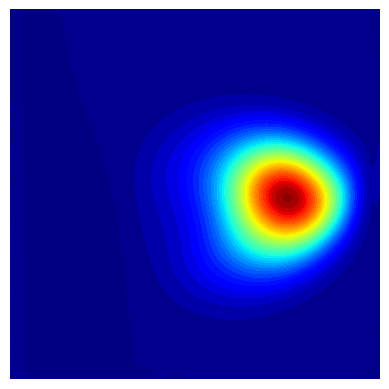}}
\subfloat{\includegraphics[height = 0.195\linewidth]{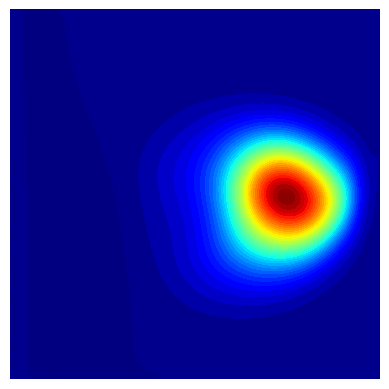}}\vspace{-1em}

\subfloat{\includegraphics[height = 0.195\linewidth]{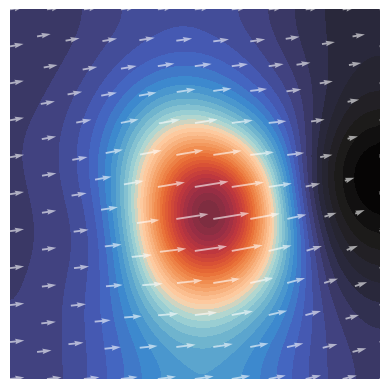}}
\subfloat{\includegraphics[height = 0.195\linewidth]{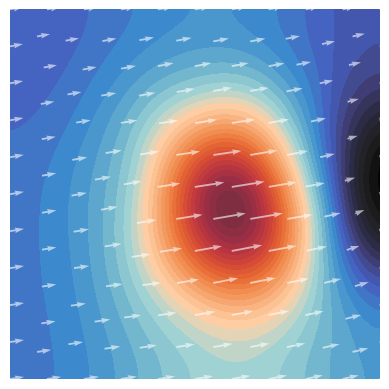}}
\subfloat{\includegraphics[height = 0.195\linewidth]{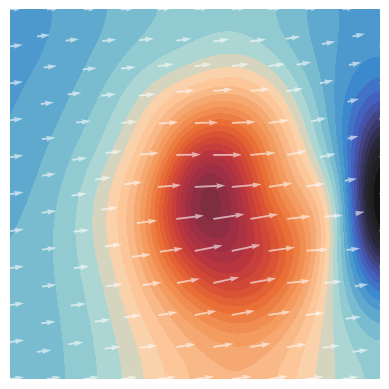}} 
\subfloat{\includegraphics[height = 0.195\linewidth]{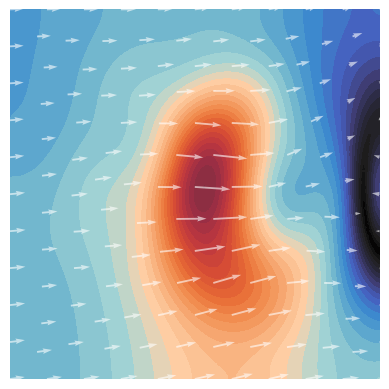}} 
\hphantom{\subfloat{\includegraphics[height = 0.195\linewidth]{Images/Vacuum/Control_test1_t=0.75.png}}}

\subfloat{\includegraphics[height = 0.195\linewidth]{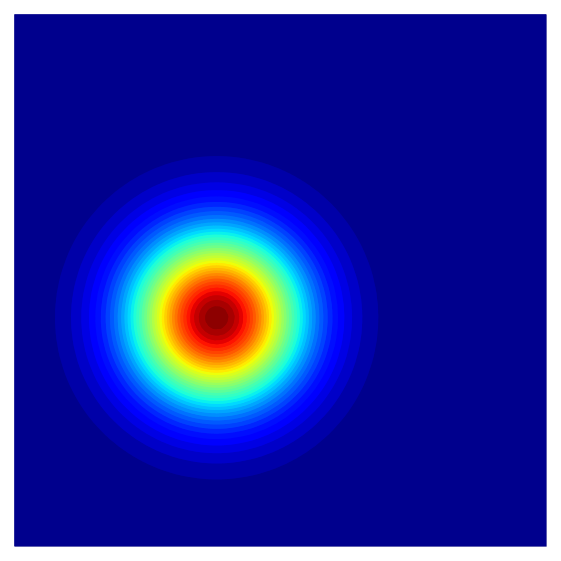}}
\subfloat{\includegraphics[height = 0.195\linewidth]{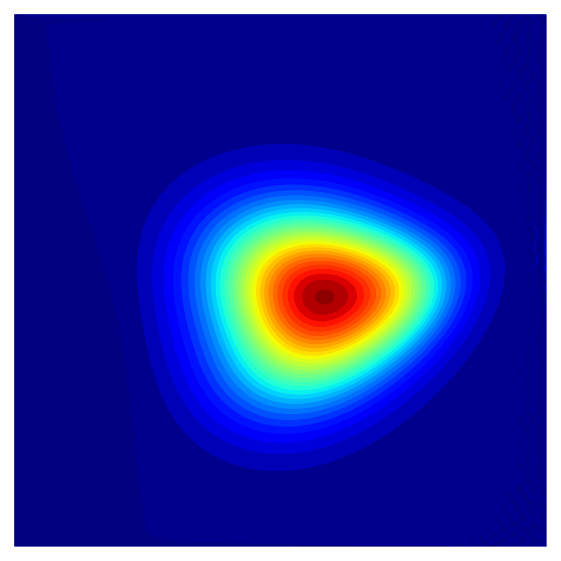}}
\subfloat{\includegraphics[height = 0.195\linewidth]{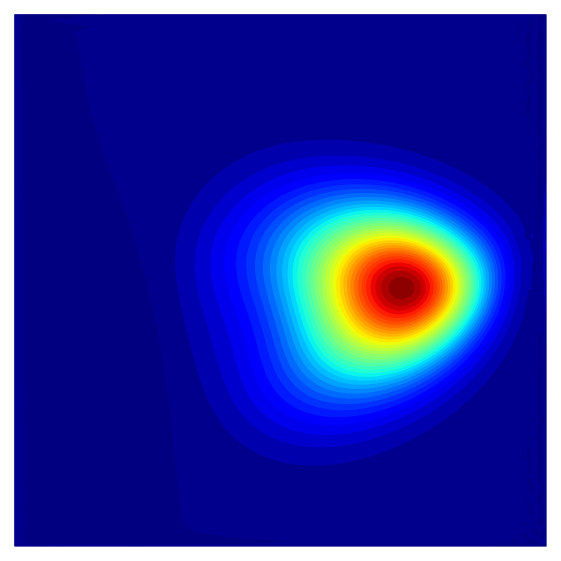}}
\subfloat{\includegraphics[height = 0.195\linewidth]{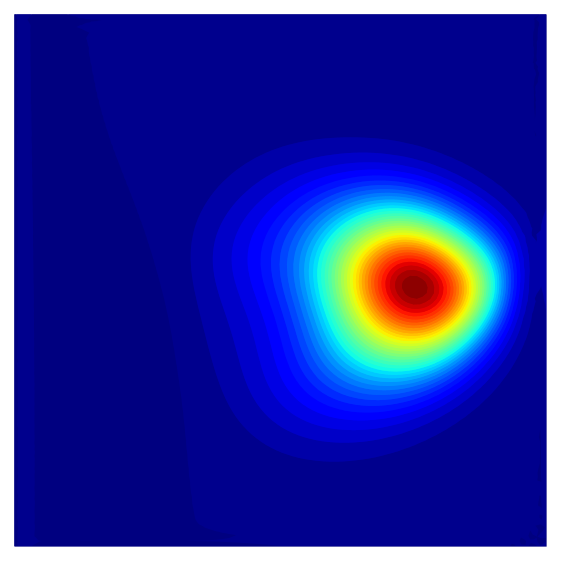}}
\subfloat{\includegraphics[height = 0.195\linewidth]{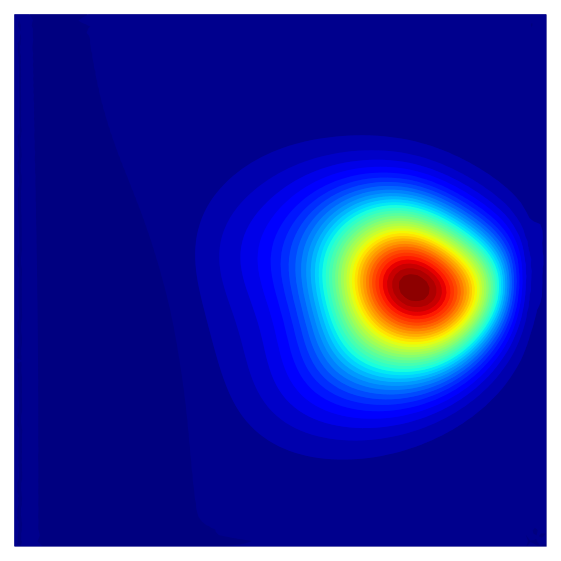}}\vspace{-1em}

\subfloat{\includegraphics[height = 0.195\linewidth]{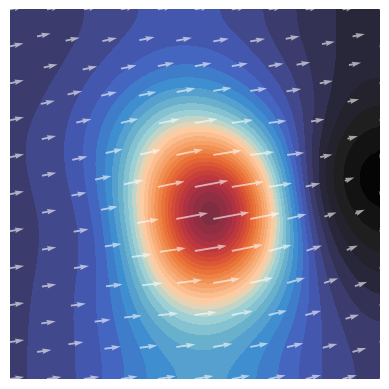}}
\subfloat{\includegraphics[height = 0.195\linewidth]{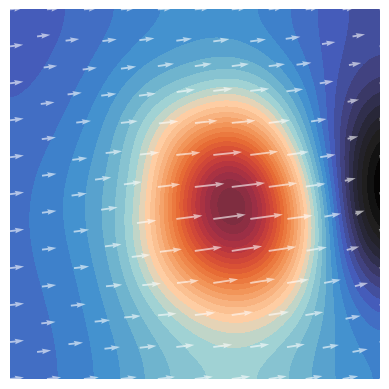}}
\subfloat{\includegraphics[height = 0.195\linewidth]{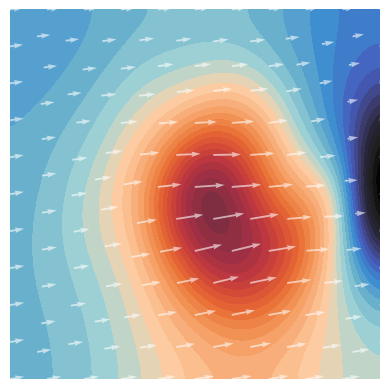}} 
\subfloat{\includegraphics[height = 0.195\linewidth]{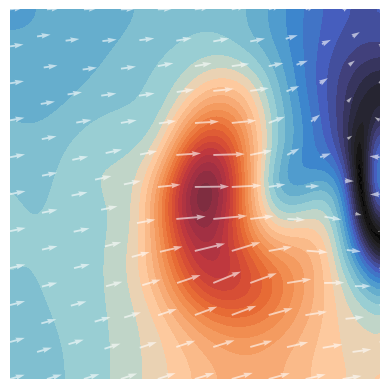}} 
\hphantom{\subfloat{\includegraphics[height = 0.195\linewidth]{Images/Vacuum/Control_test2_t=0.75.png}}}

\caption{\textit{Test 1.1}. Optimal transport in a vacuum. First row: discrepancy between the target configuration $\mathbf{y}_d$ centered at $(\mu_1^d, \mu_2^d) = (0.48, -0.03)$ and the state $\yh(t)$ considering the deep learning-based reduced order feedback controller (left) and the latent feedback loop (right) at different time instants in the uncontrolled setting ($\uh = 0$) and when exploiting the optimal control predicted by the policy ($\uh = \tildeuh$). Other panels: system evolution driven by policy-based controls related to an initial state centered at $(\mu_1^0, \mu_2^0) = (-0.24, -0.14)$ and scenario parameters $\mus = (0.48, -0.03)$, while exploiting the loop closure at the full-order (second and third rows) and latent (fourth and fifth rows) level. The control velocity fields on $\Omega$ are depicted through vector fields, with the underlying colours corresponding to their magnitude.}
\label{fig:vacuum_test}
\end{figure}

After data generation, dimensionality reduction and neural networks training in the offline phase, we can now control the dynamical system taking into account new initial configurations and scenarios, which have not been seen during training. Figure~\ref{fig:vacuum_test} displays the state evolution and the corresponding control actions related to a test case with initial state position $(\mu_1^0, \mu_2^0) = (-0.24, -0.14)$ and scenario parameters $\mus = (0.48, -0.03)$. Specifically, both the feedback loops at the full-order and latent levels effectively steer the state density towards the desired target, with a decreasing discrepancy $||\yh(t)-\mathbf{y}_d||$ over time. As far as computational times are concerned, the proposed controller with model closure at full-order level requires $0.81$ seconds to provide control actions and simulate states for all the $N_t$ time steps in the test setting considered, while it reduces to $0.028$ seconds in the case of latent feedback loop, with a remarkably high speed-up with respect to full-order methods ($1000\times$ for the loop closure at the full-order level, $32000\times$ for the latent feedback loop). Note that both proposed control strategies remain effective even when considering different (sufficiently small) time steps $\Delta t$ and different final times $T$, thus allowing control actions to be applied to the system at different instants (according to the delay in the state computation) and an arbitrary number of times.

%%%%%%%%%%%%%%%%%%%%%%%%%%%%%%%%%%%%%%%%%%%%%%%%%%%%%%%%%%%%%%%%%%%%%%%%%%%%%%%%%%%%%%%%%%%%%%%%%%%%%%%%%%%%%%%%%%%%%%%%%%%%%%%%%%%%%%%%%%%%%%%%%%%%%%%%%%

\subsection{Optimal transport in a fluid}
\label{subsec:fluid}

In this application, we focus on a more challenging optimal transport problem where {\em (i)} a rounded obstacle is added in the middle of the square $(-1,1)^2$ -- that is, the domain now becomes $\Omega = (-1,1)^2 \setminus \mathcal{B}_{0.15}(0,0)$ where $\mathcal{B}_{0.15}(0,0)$ is the circle centered at $(0,0)$ with radius $0.15$ -- and {\em (ii)} an underlying fluid flow brings the state upwards out of the domain. In the context of robotic swarms, this external fluid may represent a wind field (for aerial vehicles) or a water current (for underwater applications), affecting particle movement. The state system is now described in terms of the following Fokker-Plank equation
\begin{equation}
\label{eq:FP_fluid}
\begin{cases}
      \dfrac{\partial y}{\partial t} + \nabla \cdot (- \nu \nabla y + \mathbf{u} y + \mathbf{v} y) = 0 \qquad
      & \text{in} \ \Omega \times (0,T]
      \\
      (-\nu \nabla y + \mathbf{u} y + \mathbf{v} y) \cdot \mathbf{n} = 0
      \qquad
      & \text{on} \ \partial \Omega  \times (0,T]
      \\
       y(0) = y_0(\mu_1^0, \mu_2^0)
       \qquad
       & \text{in} \ \Omega \times \{t = 0\}
\end{cases}
\end{equation}
where, differently from Equation~\eqref{eq:FP}, an additional transport term with velocity $\mathbf{v}: \Omega \to \R^2$ is considered in order to describe the fluid flow surrounding the density. Note that, while $\mathbf{v}$ steers the density upwards out of the domain complicating the control task, the state is passive with respect to this transport effect. Specifically, the velocity field $\mathbf{v}$ in $\Omega$ is modeled through (steady, for the sake of simplicity) Navier-Stokes equations, that are
\begin{equation}
    \begin{cases}
       -\mu \Delta \mathbf{v} + (\mathbf{v} \cdot \nabla) \mathbf{v} + \nabla p = 0  \qquad &\mathrm{in} \ \Omega \\
       \mathrm{div }\  \mathbf{v} = 0  \qquad &\mathrm{in} \ \Omega \\
       \mathbf{v} = \mathbf{0}  \qquad &\mathrm{on} \ \Gamma_{\mathrm{obs}}\\
       \mathbf{v} = \mathbf{v}_{\text{in}}(\gamma_{\text{in}}, \alpha_{\text{in}}) \qquad &\mathrm{on} \ \Gamma_{\mathrm{in}} \\
       \mathbf{v} \cdot \mathbf{n} = 0  \qquad &\mathrm{on} \ \Gamma_{\mathrm{walls}} \\
       (\mu \nabla \mathbf{v} - p)\mathbf{n} \cdot \mathbf{t} = 0  \qquad &\mathrm{on} \ \Gamma_{\mathrm{walls}} \\
       (\mu \nabla \mathbf{v} - p)\mathbf{n} = 0  \qquad &\mathrm{on} \ \Gamma_{\mathrm{out}}
\end{cases}
\label{eq:NS}
\end{equation}
where $\mu = 0.01$ is the kinematic viscosity, $p: \Omega \to \R$ is the pressure field and $\mathbf{t}$ is the tangential versor to the boundary $\partial \Omega$. The underlying fluid enters the domain $\Omega = (-1,1)^2 \setminus \mathcal{B}_{0.15}(0,0)$ from the bottom side (inflow) $\Gamma_{\text{in}} = \partial \Omega \cap \{ x_2 = -1 \}$ with a Dirichlet boundary datum
\[
\mathbf{v}_{\text{in}}(\gamma_{\text{in}}, \alpha_{\text{in}}) = \begin{bmatrix}(x_1 + 1)(1 - x_1) \gamma_{\text{in}}\sin(\alpha_{\text{in}}) \\ \gamma_{\text{in}}\cos(\alpha_{\text{in}})
\end{bmatrix}
\]
where $\gamma_{\text{in}}$ and $\alpha_{\text{in}}$ are, respectively, the inflow velocity intensity and the angle of attack, while the parabolic factor $(x_1 + 1)(1 - x_1)$ is helpful to avoid singularities at the domain corners. On the outflow, that is the top side $\Gamma_{\text{out}} = \partial \Omega \cap \{ x_2 = 1\}$ where the flow exits the domain, homogeneous Neumann boundary conditions are considered. Moreover, we employ free-slip condition on the sidewalls $\Gamma_{\text{wall}} = \partial \Omega \cap \{x_1 = \pm 1\}$ and no-slip conditions on the obstacle boundary $\Gamma_{\text{obs}}$.

\begin{figure}
\centering
\subfloat{\includegraphics[scale = 0.23]{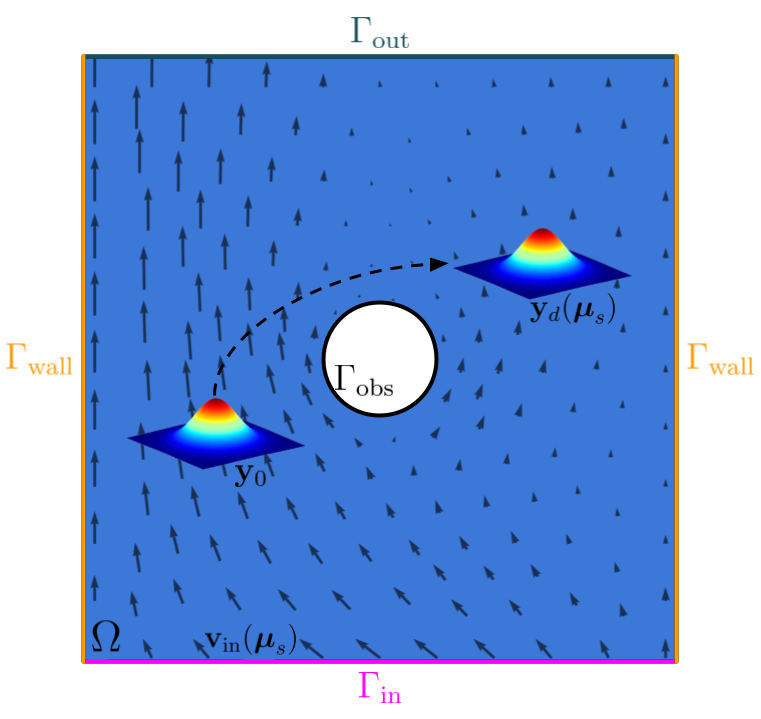}}

\subfloat{\includegraphics[height = 0.195\linewidth]{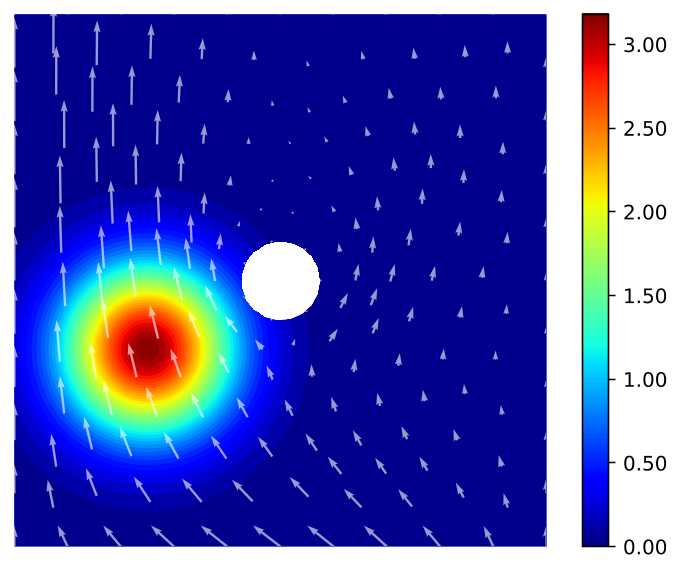}} \quad    \subfloat{\includegraphics[height = 0.195\linewidth]{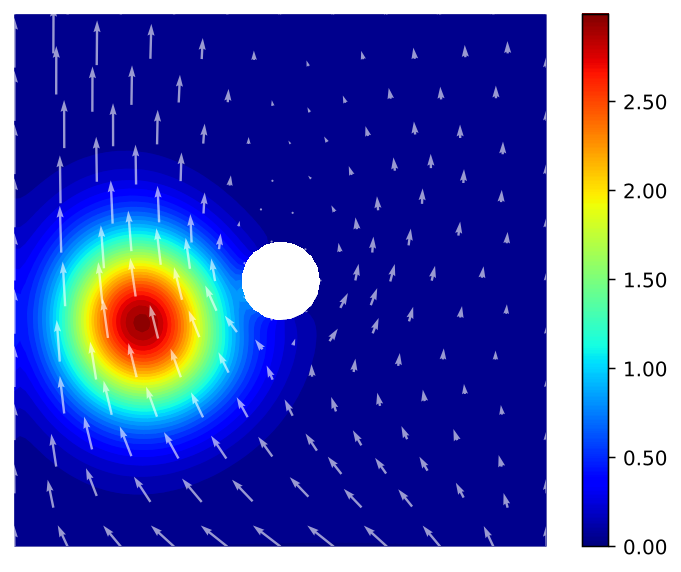}} \quad
\subfloat{\includegraphics[height = 0.195\linewidth]{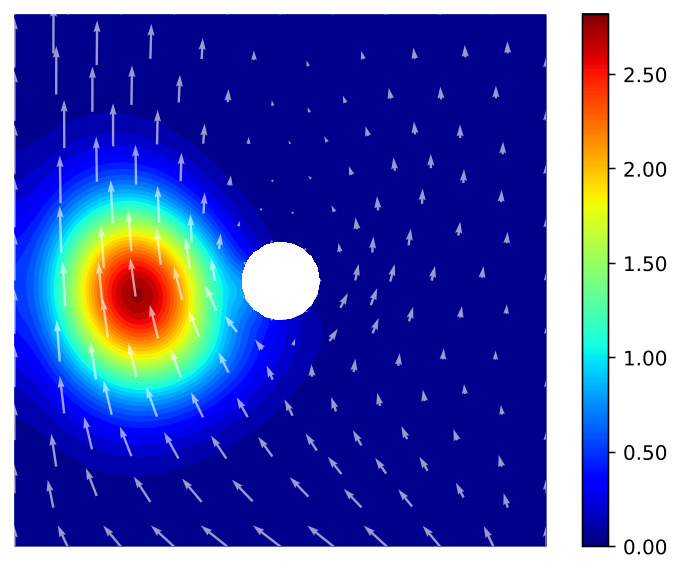}} \quad
\subfloat{\includegraphics[height = 0.195\linewidth]{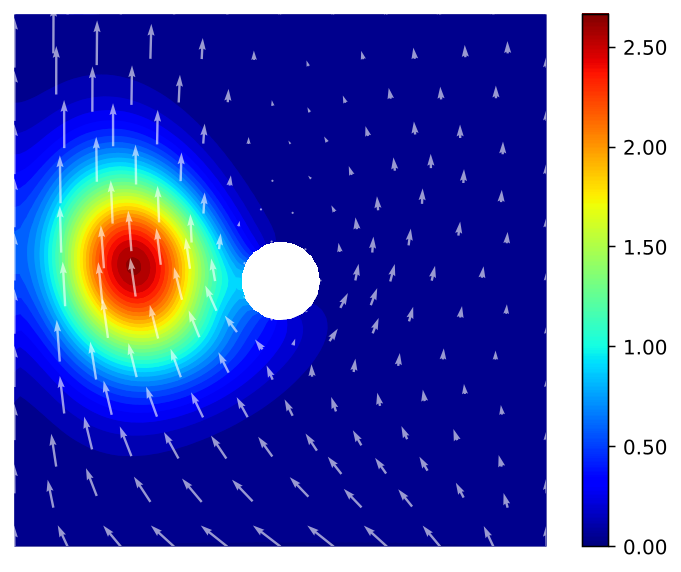}}
    
\subfloat{\includegraphics[height = 0.195\linewidth]{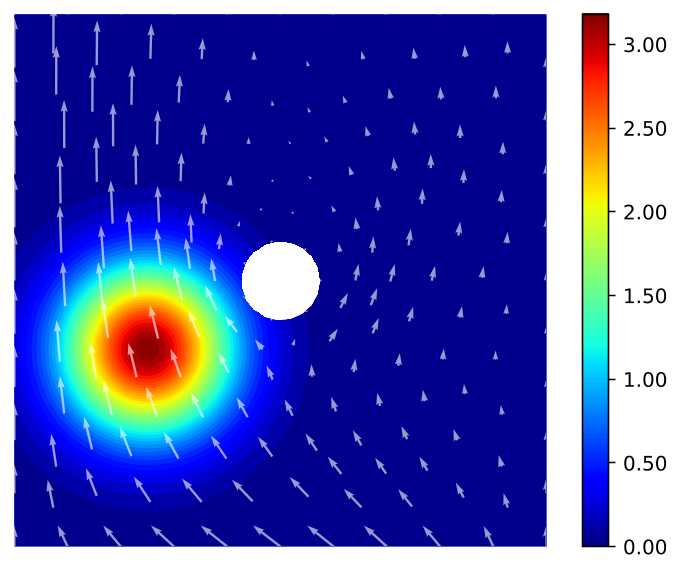}} \quad
\subfloat{\includegraphics[height = 0.195\linewidth]{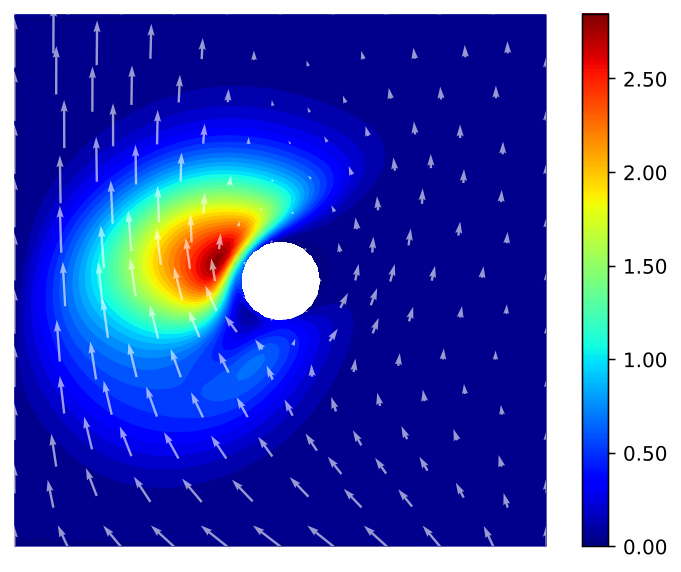}} \quad
\subfloat{\includegraphics[height = 0.195\linewidth]{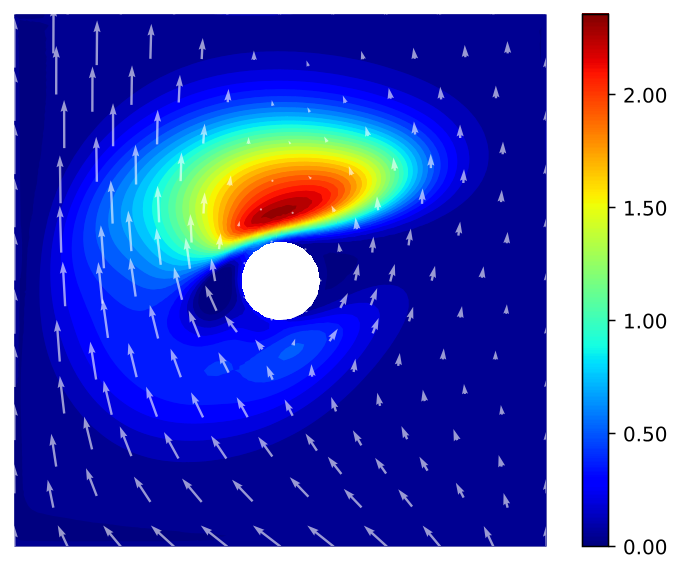}} \quad
\subfloat{\includegraphics[height = 0.195\linewidth]{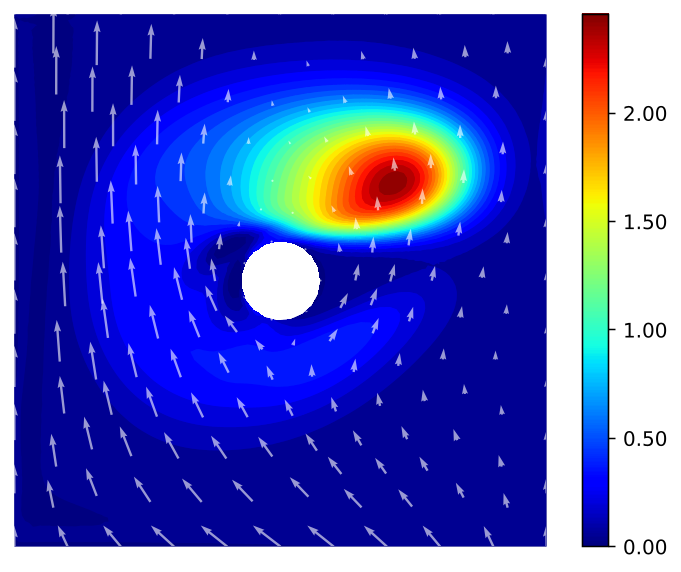}}

\subfloat{\includegraphics[height = 0.195\linewidth]{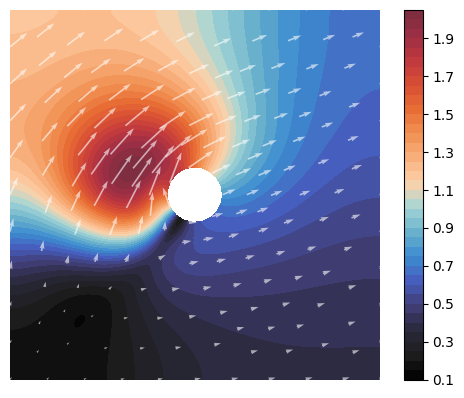}} \quad
\subfloat{\includegraphics[height = 0.195\linewidth]{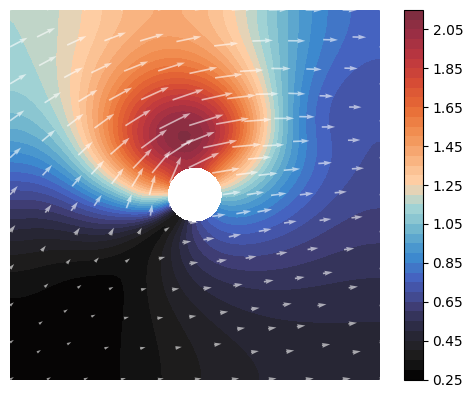}} \quad
\subfloat{\includegraphics[height = 0.195\linewidth]{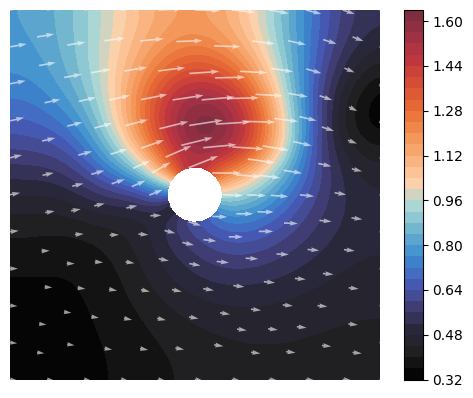}} \quad
\subfloat{\includegraphics[height = 0.195\linewidth]{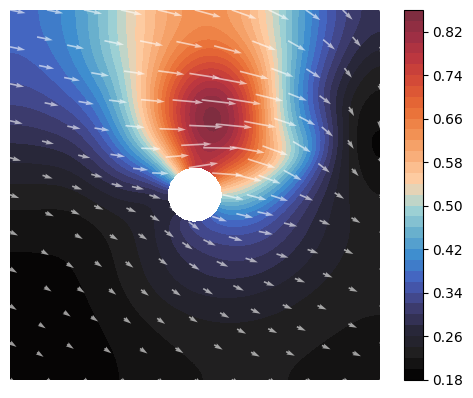}}

\caption{\textit{Test 1.2}. Optimal transport in a fluid. Top: representation of an optimal state trajectory in a fluid within the domain $\Omega$, where $\yinit$ stands for the initial density centered at $(\mu_1^0, \mu_2^0) = (-0.5, -0.25)$, $\mathbf{y}_d(\mus)$ represents the target configuration centered at $(\mu_1^d, \mu_2^d) =(0.5, 0.39)$ and $\mathbf{v}_{\text{in}}(\mus)$ is the inlet velocity with inflow intensity $\gamma_{\text{in}} = 0.40$ and angle of attack $\alpha_{\text{in}} = -0.92$, that is $\mus=(0.39, 0.40, -0.92)$. Other panels: space-varying uncontrolled state, optimal state and control at $t=0,0.25,0.5,0.75$ related to the scenario parameters $\mus = (0.39, 0.40, -0.92)$. The underlying fluid velocity vector field $\mathbf{v}_h$ on $\Omega$ is depicted together with the state. The control velocity fields on $\Omega$ are depicted through vector fields, with the underlying colours corresponding to their magnitude.}
\label{fig:fluid_setting}
\end{figure}

Similarly to the optimal transport problem in a vacuum detailed in Section~\ref{subsec:vacuum}, by minimizing the loss function in Equation~\eqref{eq:J} with $\beta = \beta_g = 0.2$, we aim to find the optimal control action that brings the state density from an initial configuration $y_0$ to the target destination $y_d$, that are
\begin{align*} 
y_0(-0.5, \mu_2^0) &= \dfrac{10}{\pi} \exp{(- 10(x_1 + 0.5)^2 - 10(x_2 - \mu_2^0)^2)}
\\
y_d(0.5,\mu_2^d) &= \dfrac{10}{\pi} \exp{(- 10(x_1 - 0.5)^2 - 10(x_2 - \mu_2^d)^2)}
\end{align*}
Note that, due to the presence of the obstacle, we set the mean $x_1$-coordinates of $y_0$ and $y_d$ equal to $-0.5$ and $0.5$, respectively. Note also that the integral on $\partial \Omega$ in the loss functional $J$ helps to avoid collisions both with domain boundaries and the circular obstacle, as typically crucial in safety critical applications related to robotic swarms. In this setting, we consider the vector of scenario parameters $\mus = (\mu_2^d,\gamma_{\text{in}}, \alpha_{\text{in}})$ in order to deal with both different final coordinates and different underlying flows. Figure~\ref{fig:fluid_setting} shows the test case setting, along with a comparison between the uncontrolled and the optimal state trajectories.

As far as the finite element discretization is concerned, we consider a mesh discretizing the domain $\Omega$ with $7681$ nodes and step size $h = 0.05$. While $\mathbb{P}_1$ finite elements are taken into account when discretizing state and control in Equation~\eqref{eq:FP_fluid}, Taylor-Hood $\mathbb{P}_2$-$\mathbb{P}_1$ elements are considered in Equation~\eqref{eq:NS} to guarantee the inf-sup stability and the well-posedness of the problem \cite{Quarteroni2017}. The degrees of freedom of the discrete state $\yh$, control $\uh$, velocity $\mathbf{v}_h$ and pressure $\mathbf{p}_h$ variables end up being, respectively, $N_h^y = 7681$, $N_h^u = 15362$, $N_h^v = 60732$ and $N_h^p = 7681$. Note that, to keep the Péclet number always smaller than $1$, thus avoiding instabilities while solving Equation~\eqref{eq:FP_fluid}, an artificial diffusion equal to $0.5 h \gamma_{\text{in}}$ is added to the diffusion coefficient $\nu$ \cite{Quarteroni2017}. Instead, regarding the time discretization, we consider an evenly spaced time grid spanning $[0,T]$ with time step $\Delta t = 0.25$, where the final time $T$ is set equal to $1.5$ and $N_t = 6$.

In the data generation step, we simulate $N_s = 150$ trajectories through the adjoint method implemented in \texttt{dolfin-adjoint}, exploiting L-BFGS-B as optimization algorithm, a tolerance equal to $10^{-6}$ and $500$ as maximum number of iterations. Every full-order optimal trajectory -- which requires, on average, $22$ minutes to be computed -- is related to scenario parameters and initial $x_2$-coordinate $\mu_2^0$ randomly sampled, respectively, in the parameter space $\mathcal{P} = (-0.5, 0.5) \times (0.1, 1.0) \times (-1.0, 1.0)$ and in the interval $(-0.5, 0.5)$. Note that, since $\gamma_{\text{in}} \in (0.1, 1.0)$, the Reynolds numbers taken into account range from $20$ to $200$. The $N_s N_t = 900$ snapshots are then split into training and test set with a $90:10$ ratio. 

\begin{figure}[!ht]
\centering
\begin{minipage}{0.44\linewidth}
\centering
\includegraphics[scale = 0.5]{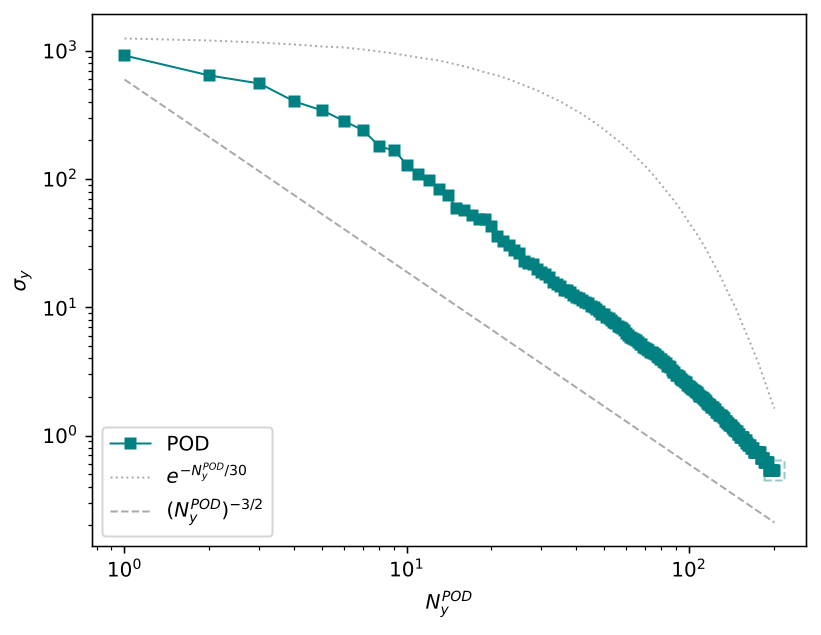}
\end{minipage}
\begin{minipage}{0.55\linewidth}
\centering
\hfill
\vspace{1.05cm}
    \subfloat{\includegraphics[height = 0.38\linewidth]{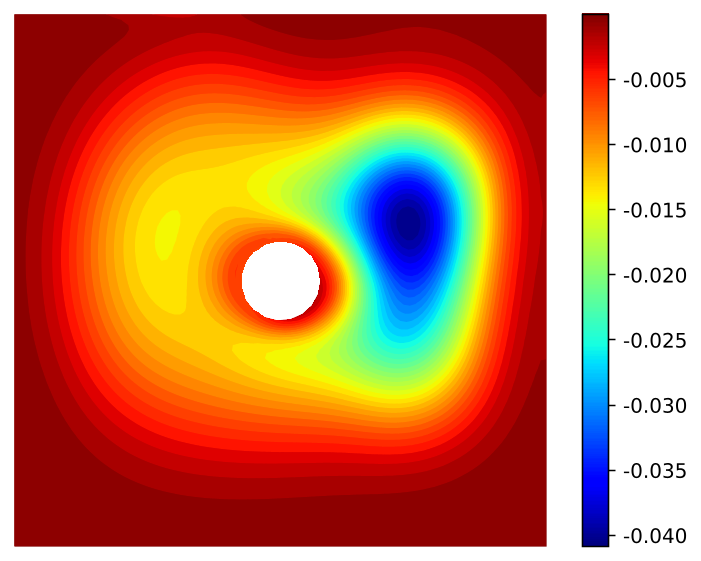}} \,
    \subfloat{\includegraphics[height = 0.38\linewidth]{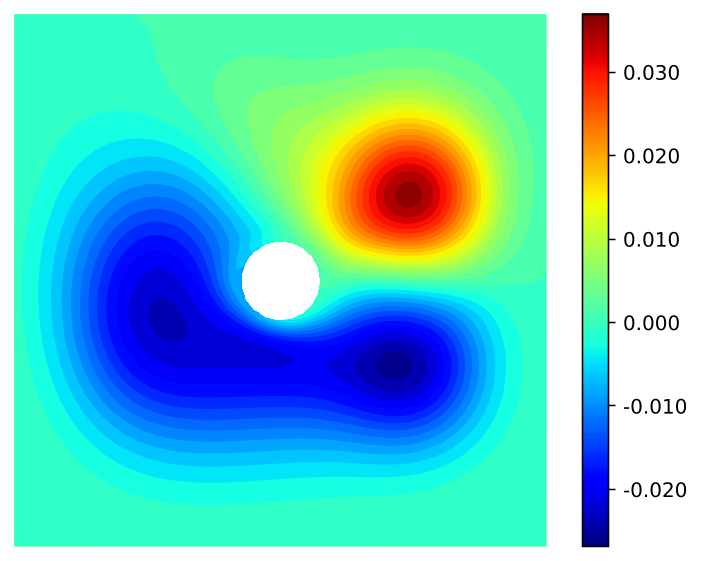}}
\end{minipage}

\begin{minipage}[c]{0.44\linewidth}
\centering
    \includegraphics[scale = 0.5]{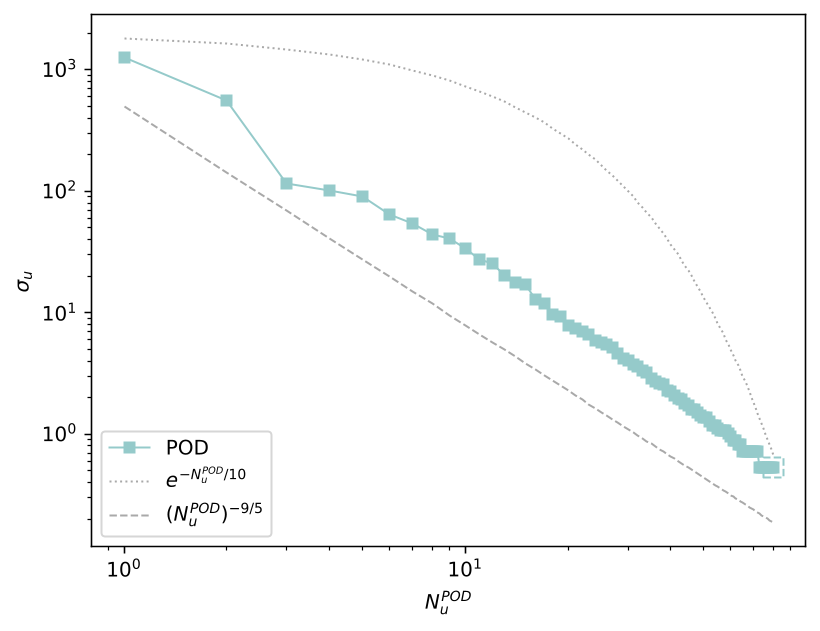}
\end{minipage}
\begin{minipage}[c]{0.55\linewidth}
\centering
\hfill
\vspace{1.05cm}
    \subfloat{\includegraphics[height = 0.38\linewidth]{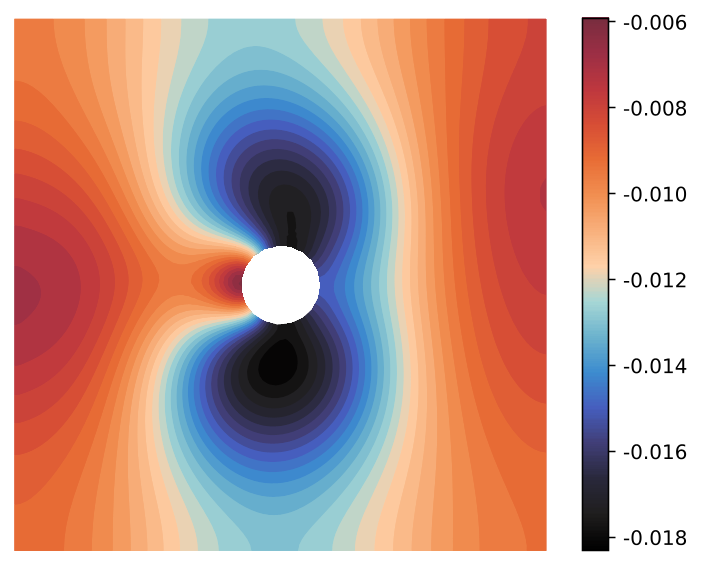}} \,
    \subfloat{\includegraphics[height = 0.38\linewidth]{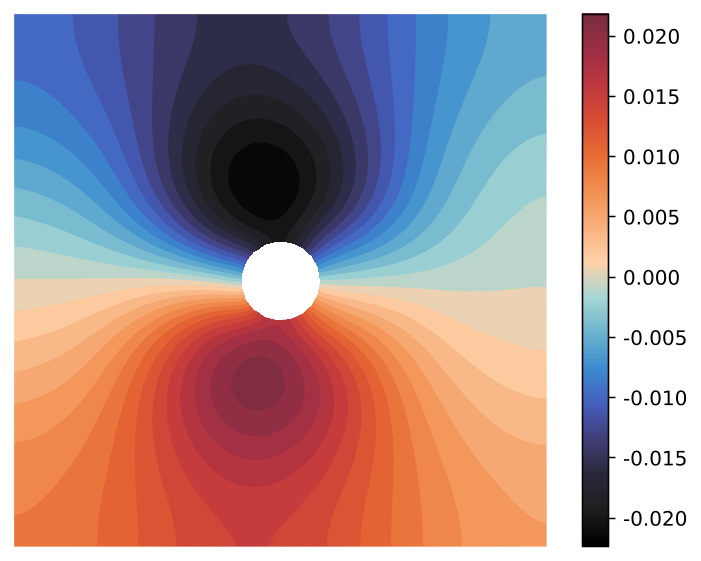}}
\end{minipage}

\begin{minipage}[c]{0.44\linewidth}
\centering
    \includegraphics[scale = 0.5]{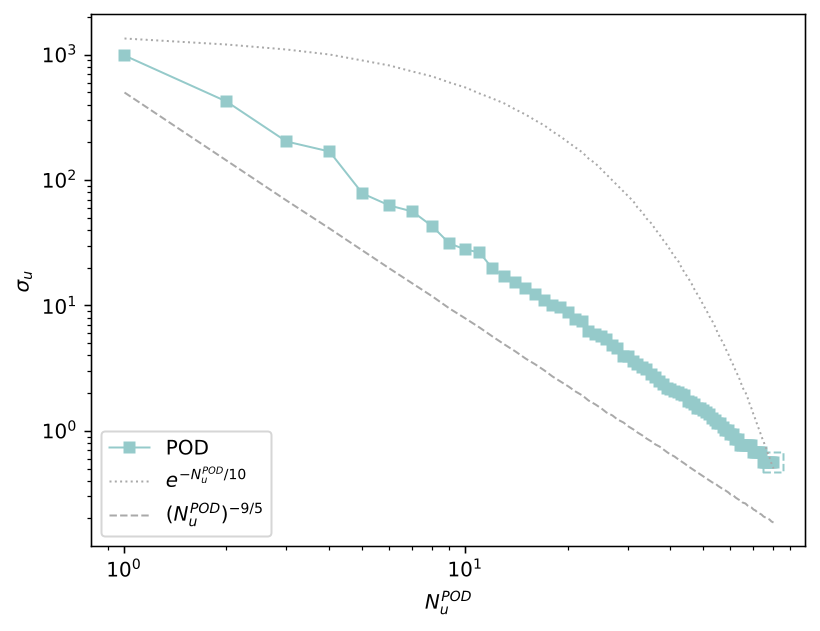}
\end{minipage}
\begin{minipage}[c]{0.55\linewidth}
\centering
\hfill
\vspace{1.05cm}
    \subfloat{\includegraphics[height = 0.38\linewidth]{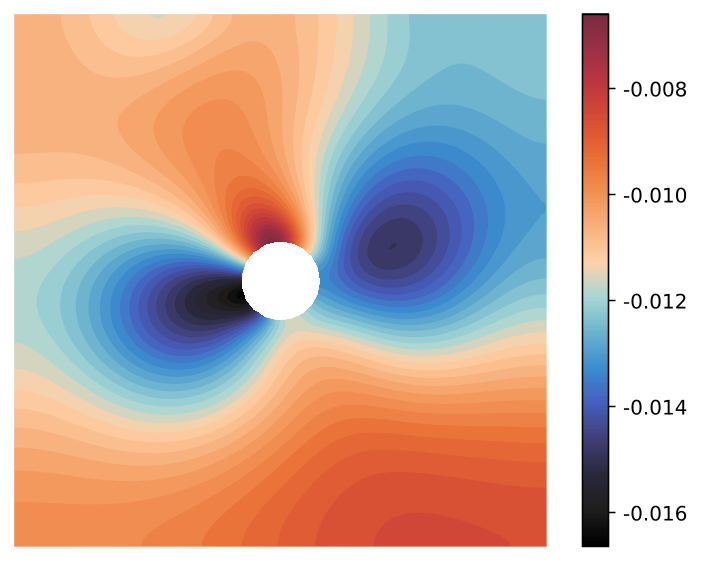}} \,
    \subfloat{\includegraphics[height = 0.38\linewidth]{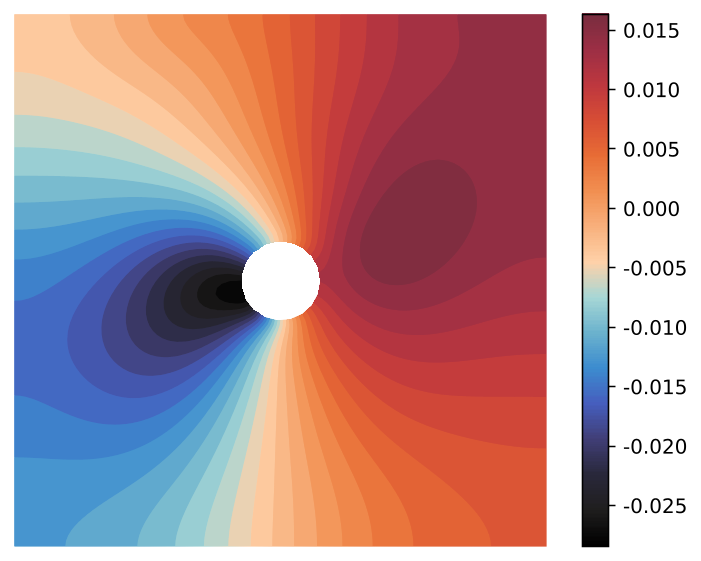}}
\end{minipage}

\caption{\textit{Test 1.2}. Optimal transport in a fluid. Singular values decay in log-log scale along with the two most energetic POD modes related to the state (top), $x_1$ component (center) and $x_2$ component (bottom) of the control.}
\label{fig:PODfluid}
\end{figure}

The state and control snapshots are then reduced through POD, looking at the singular values decays and the reconstruction errors in order to select the number of
modes to retain. Specifically, $200$ and $160$ modes are considered for state and control, respectively, ending up with reduction errors equal to $\varepsilon^y_{\mathrm{rel}} = 0.29\%$ and $\varepsilon^u_{\mathrm{rel}} = 0.29\%$. As already highlighted in the previous section, a linear projection is not enough in this context to achieve low-dimensional latent spaces for both state and control, as confirmed by the polynomial singular values decays in Figure~\ref{fig:PODfluid}.

\begin{figure}
\centering
\subfloat{\includegraphics[height = 0.195\linewidth]{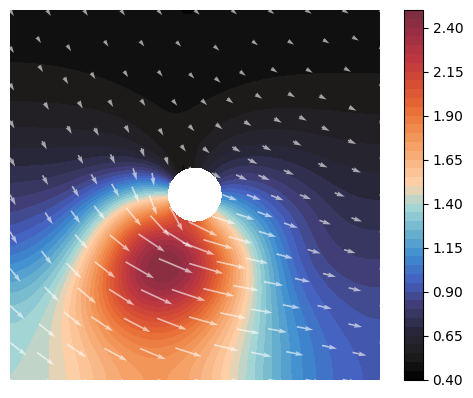}} \quad
\subfloat{\includegraphics[height = 0.195\linewidth]{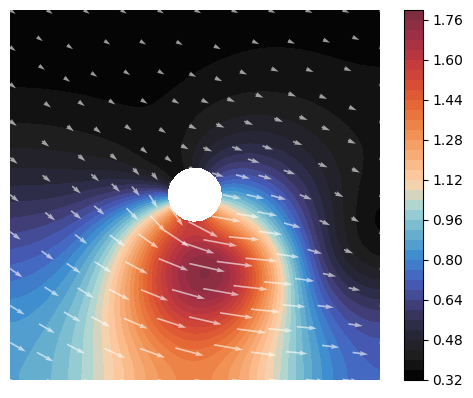}} \quad
\subfloat{\includegraphics[height = 0.195\linewidth]{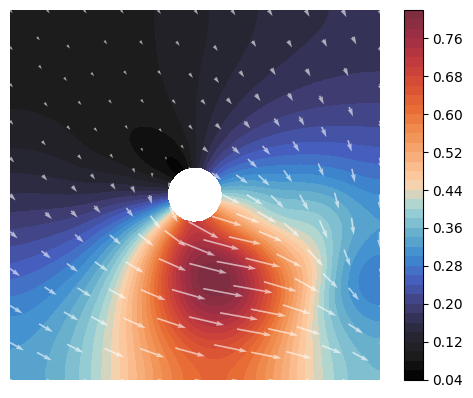}} \quad
\subfloat{\includegraphics[height = 0.195\linewidth]{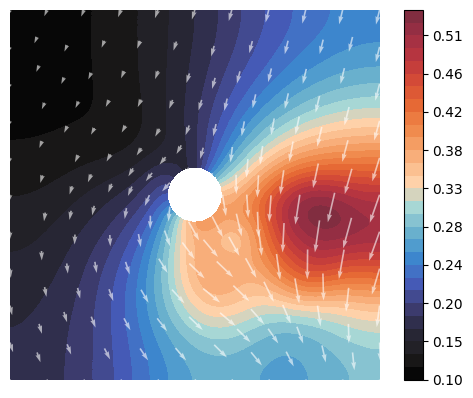}}

\subfloat{\includegraphics[height = 0.195\linewidth]{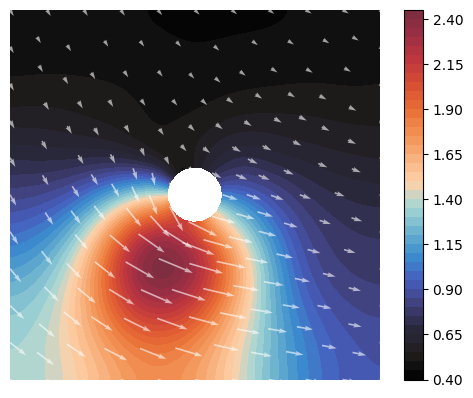}} \quad
\subfloat{\includegraphics[height = 0.195\linewidth]{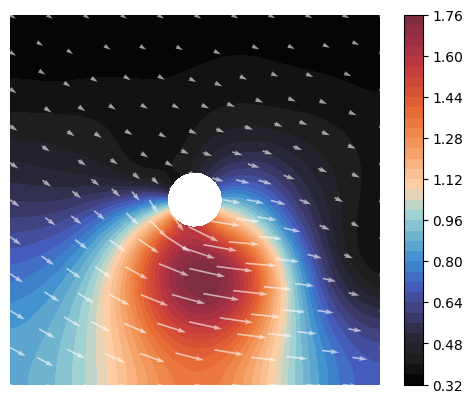}} \quad
\subfloat{\includegraphics[height = 0.195\linewidth]{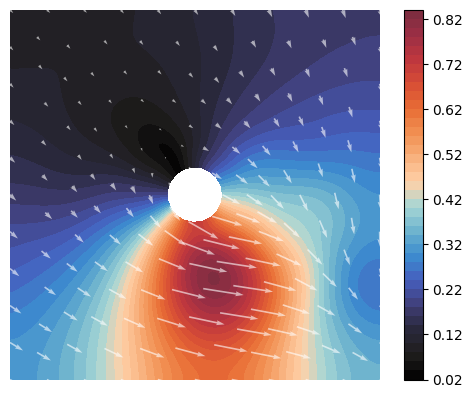}} \quad
\subfloat{\includegraphics[height = 0.195\linewidth]{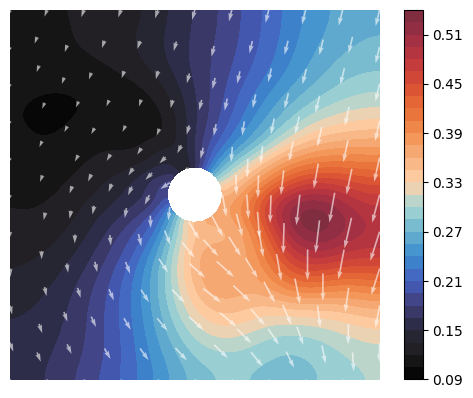}}

\subfloat{\includegraphics[height = 0.195\linewidth]{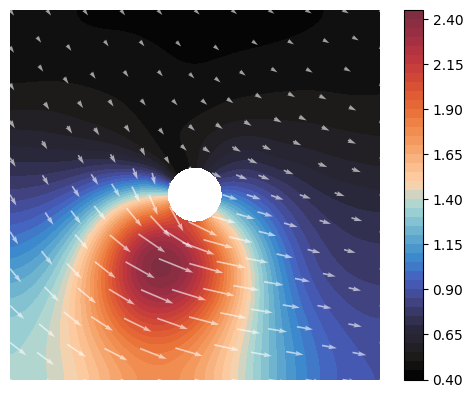}} \quad
\subfloat{\includegraphics[height = 0.195\linewidth]{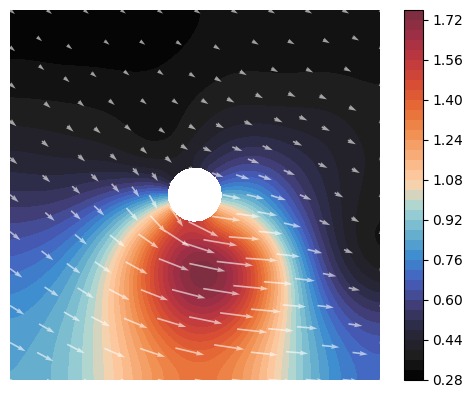}} \quad
\subfloat{\includegraphics[height = 0.195\linewidth]{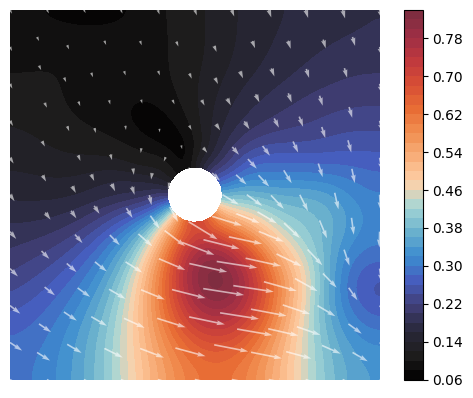}} \quad
\subfloat{\includegraphics[height = 0.195\linewidth]{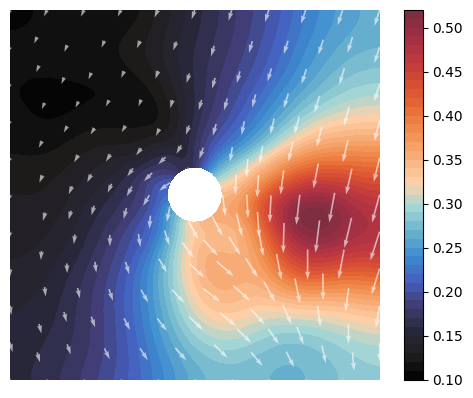}}

\subfloat{\includegraphics[height = 0.195\linewidth]{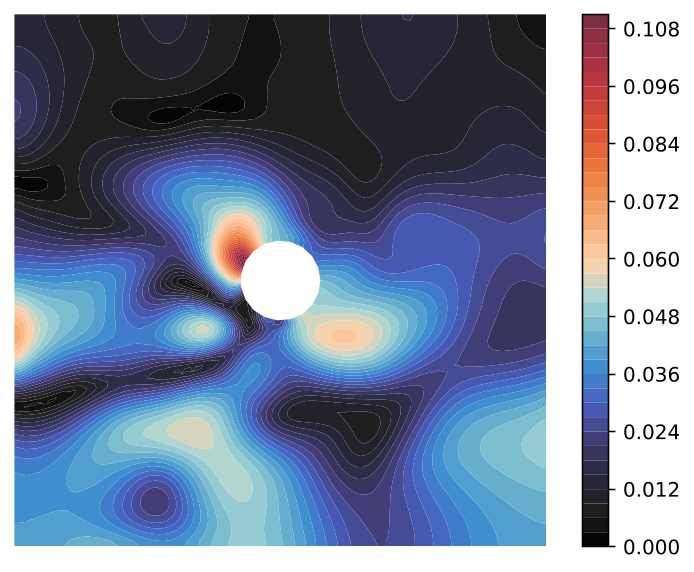}} \hspace{0.55em}
\subfloat{\includegraphics[height = 0.195\linewidth]{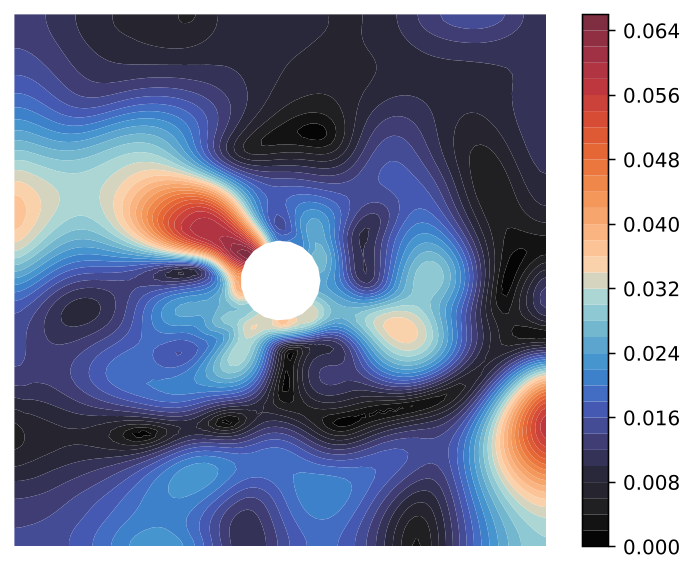}} \hspace{0.55em}
 \subfloat{\includegraphics[height = 0.195\linewidth]{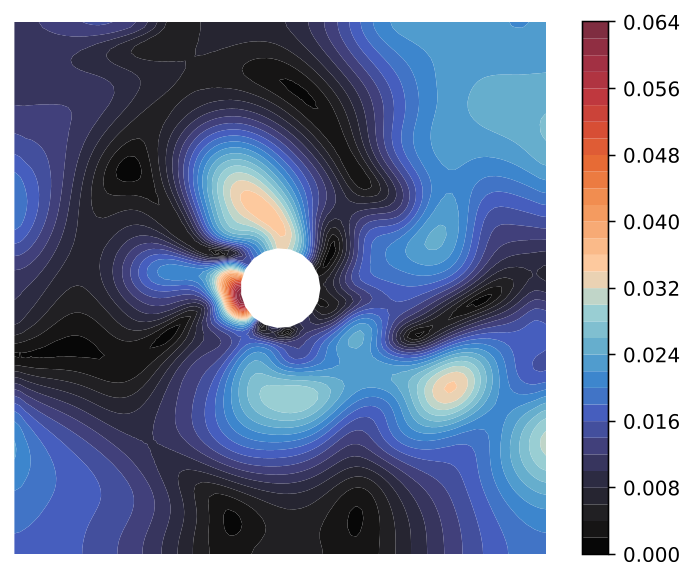}} \hspace{0.55em}
\subfloat{\includegraphics[height = 0.195\linewidth]{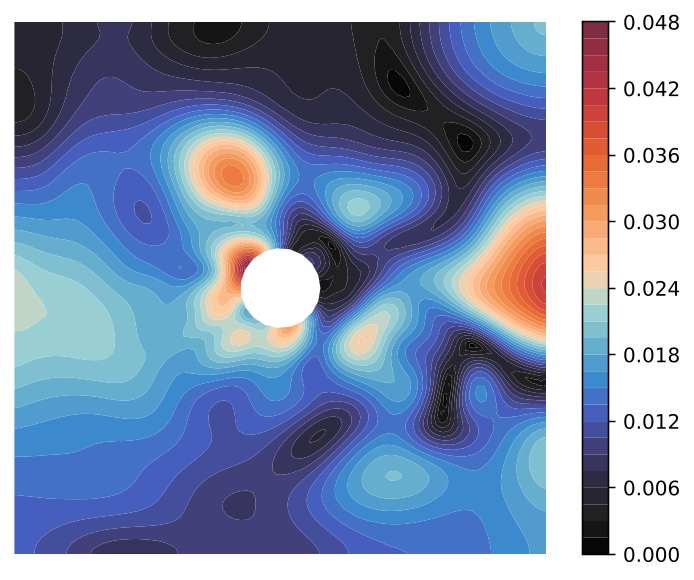}}

\subfloat{\includegraphics[height = 0.195\linewidth]{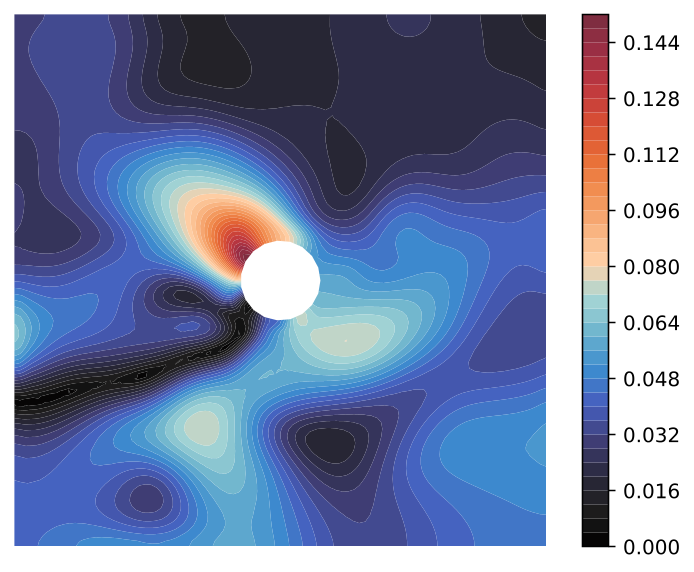}} \hspace{0.55em}
\subfloat{\includegraphics[height = 0.195\linewidth]{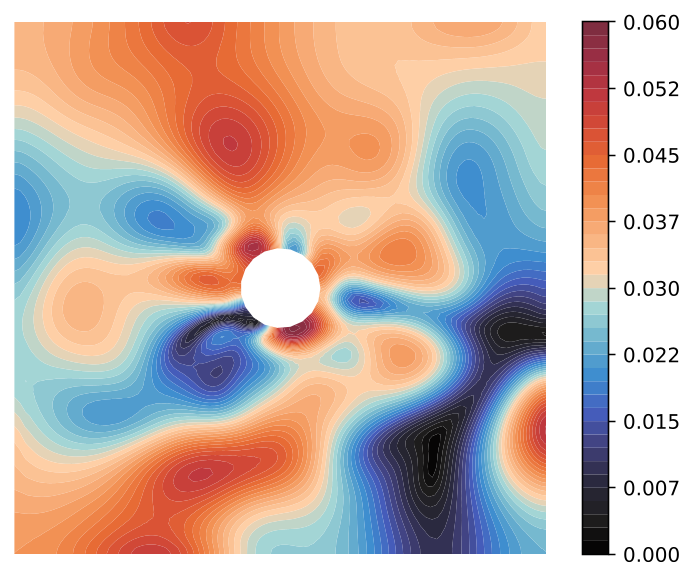}} \hspace{0.55em}
\subfloat{\includegraphics[height = 0.195\linewidth]{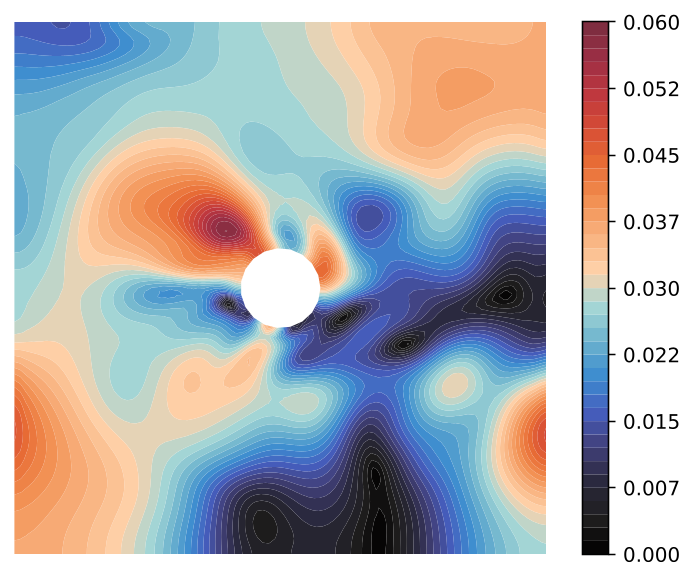}} \hspace{0.55em}
\subfloat{\includegraphics[height = 0.195\linewidth]{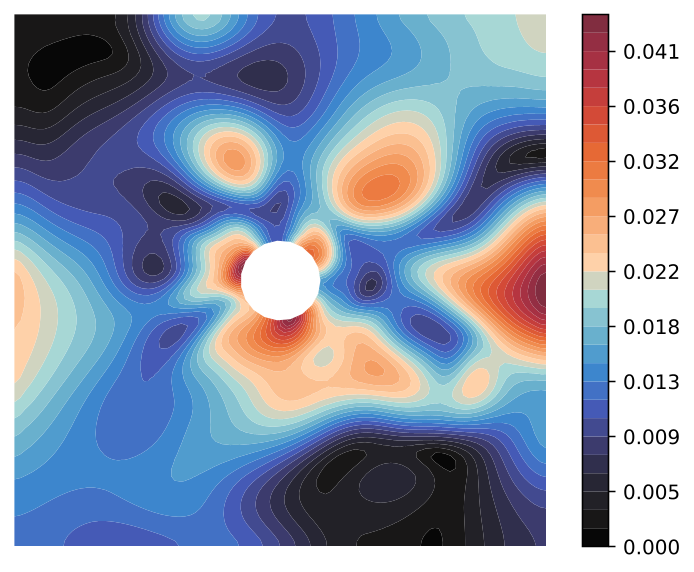}}

\caption{\textit{Test 1.2}. Optimal transport in a fluid. High-fidelity optimal control trajectory (first row), POD+AE reconstructions (second row), policy predictions (third row), POD+AE reconstruction errors (fourth row) and policy prediction errors (fifth row) at $t = 0, 0.25, 0.5, 0.75$ related to the test scenario parameters $\mus = (-0.30, 0.49, 0.59)$. The control velocity fields on $\Omega$ are depicted through vector fields, with the underlying colours corresponding to their magnitude.}
\label{fig:fluid_control}
\end{figure}

To achieve lower-dimensional subspaces, allowing for a lighter policy that is faster to train and evaluate online, we employ the POD+AE reduction strategy proposed by \cite{Fresca2022}. In particular, thanks to the nonlinear projection provided by the autoencoders with leaky Relu as activation function, we take into account latent spaces of dimension $N_y = N_u = 14$. Both the state and control encoders are made of $1$ hidden layer having $100$ neurons, while the corresponding decoders consist of $2$ hidden layers with $100$ neurons each. The low-dimensional policy $\pi_N$ is instead modeled through a deep feedforward neural network with $3$ hidden layers of $50$ neurons each and leaky Relu as activation function. Note that, together with the latent state coordinates and the scenario parameters, the input of $\pi_N$ is enriched by considering meaningful problem-driven quantities, such as $\tan(\mu_2^d)$, $\tan(\alpha_{\text{in}})$, $\gamma_{\text{in}} \sin(\alpha_{\text{in}})$ and $\gamma_{\text{in}} \cos(\alpha_{\text{in}})$. Starting from the initialization of the weights proposed by \cite{He2015}, we train the networks in $1$ hour and $25$ minutes minimizing the cumulative loss function $J_{\mathrm{NN}}$ with $\lambda_1 = \lambda_2 = \lambda_3 = 0.001$ through the L-BFGS optimization algorithm. After training, it is possible to accurately compress and reconstruct the state and control trajectories, with reconstruction errors equal to $\varepsilon^y_{\mathrm{rel}} = 4.39\%$ and $\varepsilon^u_{\mathrm{rel}} = 4.43\%$. Similarly, the low-dimensional policy $\pi_N$ is now able to correctly predict the optimal control actions starting from state data in the test set, with prediction error equal to $3.62\%$ at the latent level and $7.08\%$ after POD+AE decoding. To visually assess the accuracy and the generalization capabilities of the trained networks on test scenarios, Figure~\ref{fig:fluid_control} displays a control test trajectory along with the corresponding POD+AE reconstruction and policy prediction.

\begin{figure}
\centering
\subfloat{\includegraphics[height = 0.195\linewidth]{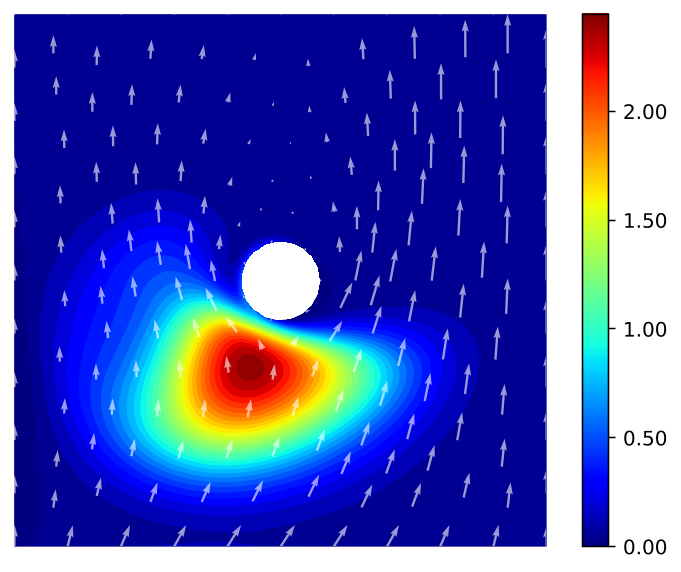}} \quad
\subfloat{\includegraphics[height = 0.195\linewidth]{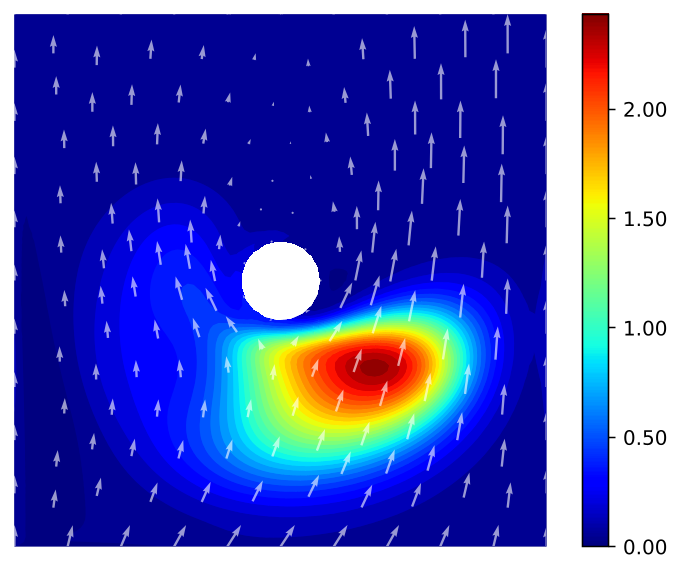}} \quad
\subfloat{\includegraphics[height = 0.195\linewidth]{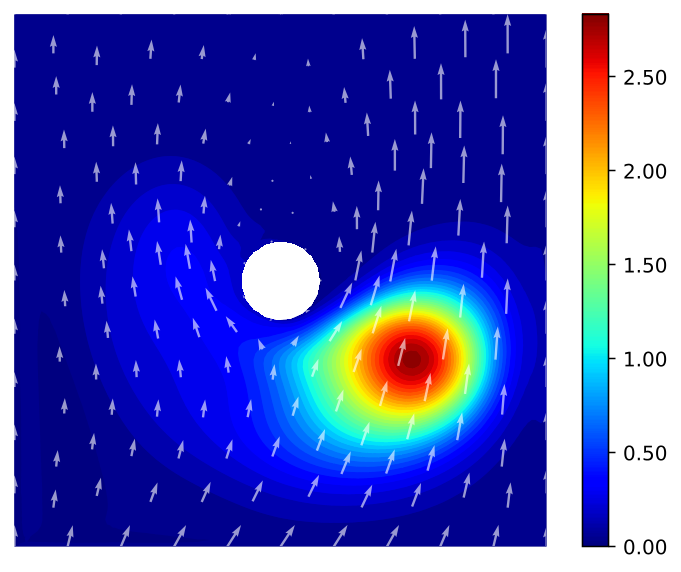}} \quad
\subfloat{\includegraphics[height = 0.195\linewidth]{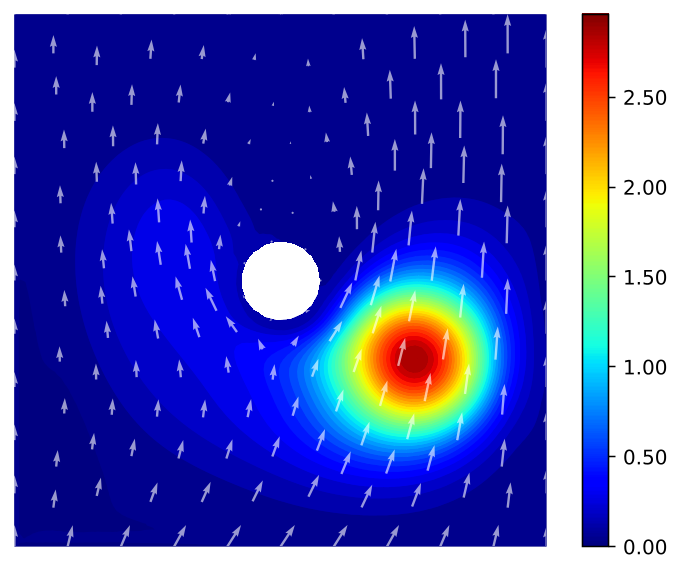}}

\subfloat{\includegraphics[height = 0.195\linewidth]{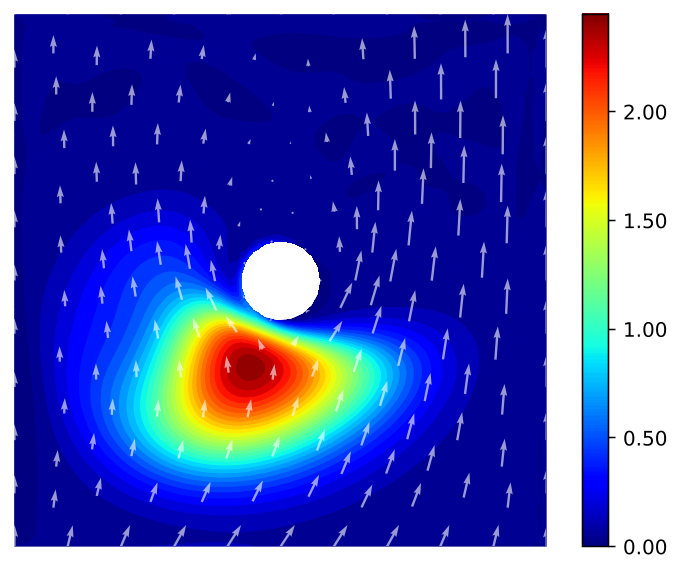}} \quad
\subfloat{\includegraphics[height = 0.195\linewidth]{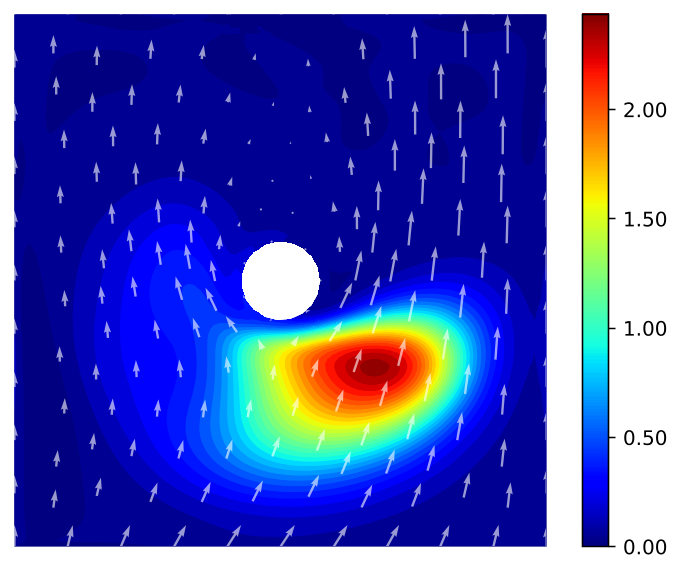}} \quad
\subfloat{\includegraphics[height = 0.195\linewidth]{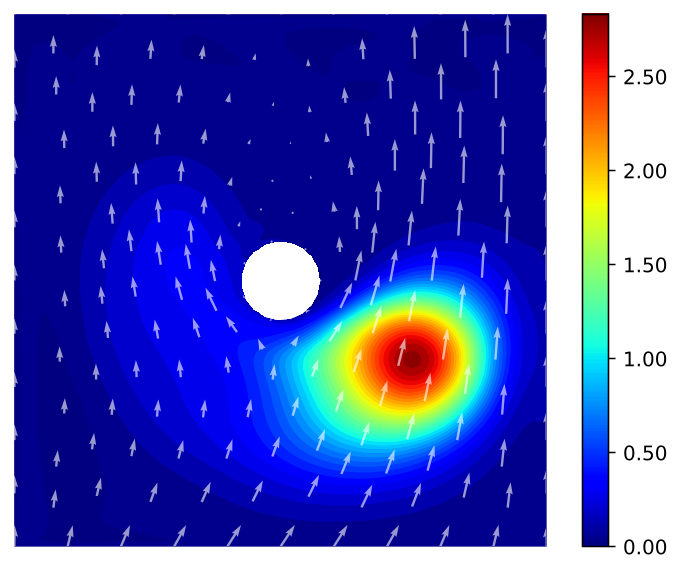}} \quad
\subfloat{\includegraphics[height = 0.195\linewidth]{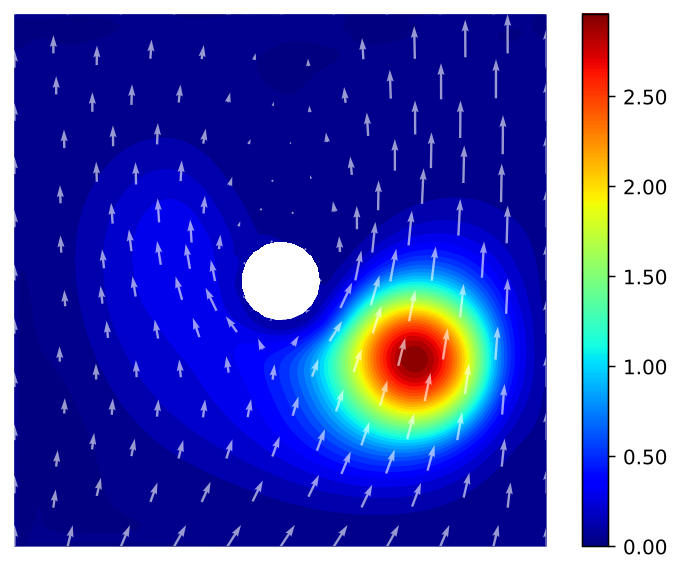}}

\subfloat{\includegraphics[height = 0.195\linewidth]{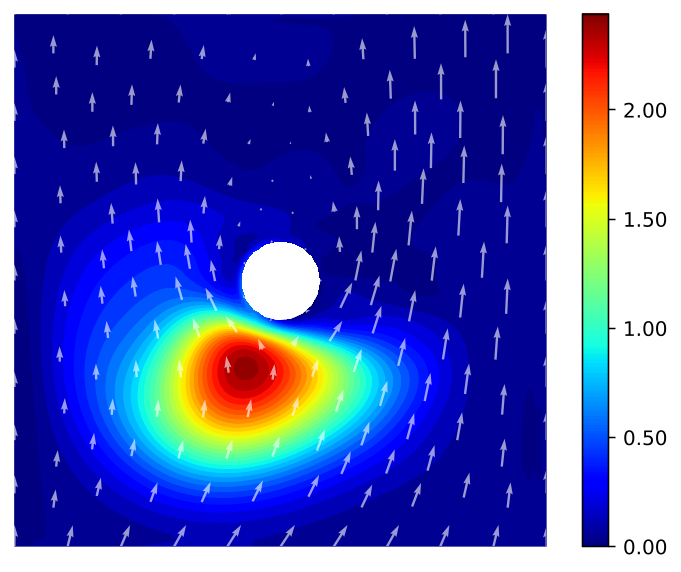}} \quad
\subfloat{\includegraphics[height = 0.195\linewidth]{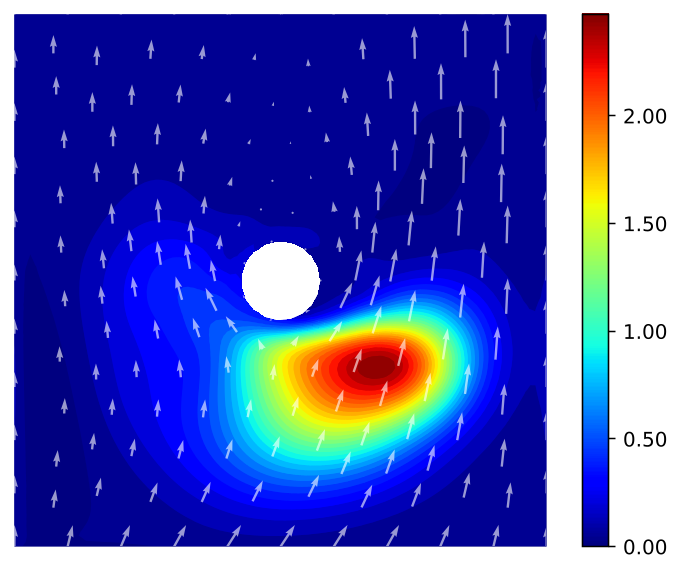}} \quad
\subfloat{\includegraphics[height = 0.195\linewidth]{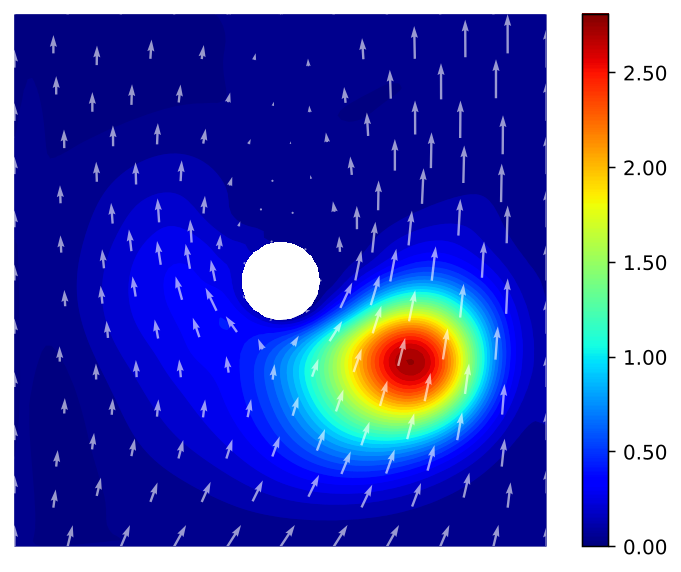}} \quad
\subfloat{\includegraphics[height = 0.195\linewidth]{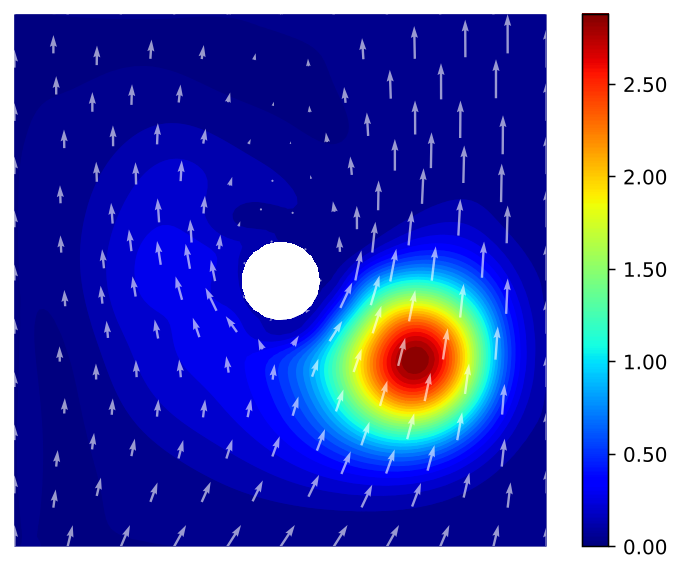}}

\subfloat{\includegraphics[height = 0.195\linewidth]{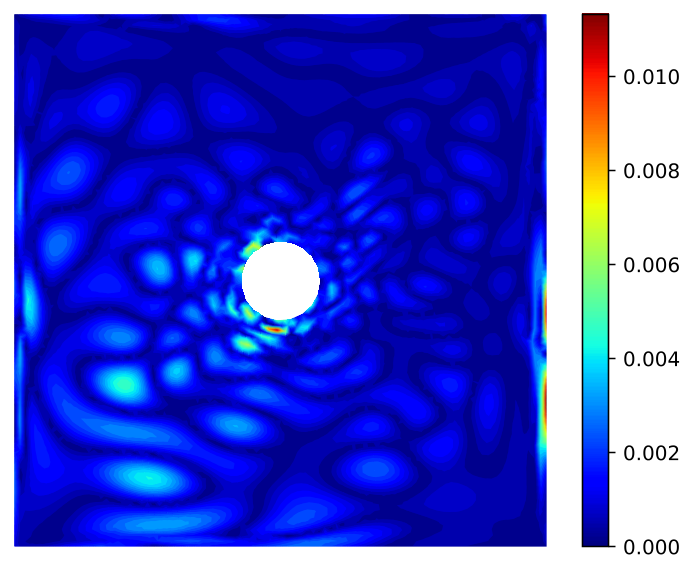}} \hspace{0.47em}
\subfloat{\includegraphics[height = 0.195\linewidth]{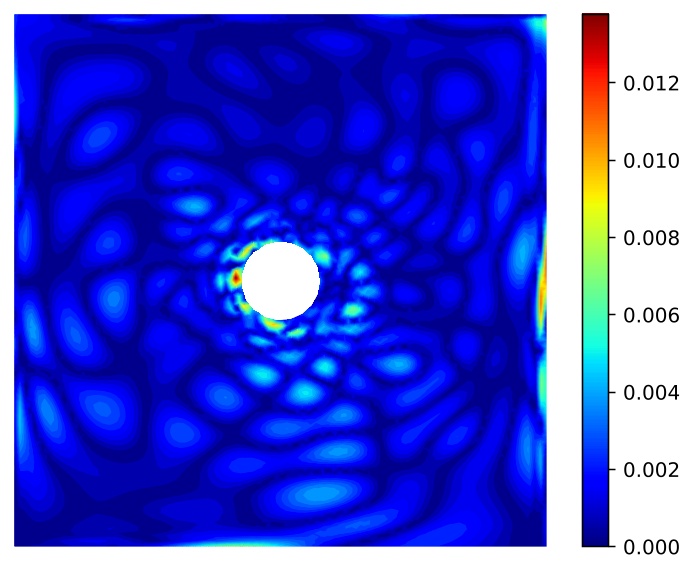}} \hspace{0.47em}
\subfloat{\includegraphics[height = 0.195\linewidth]{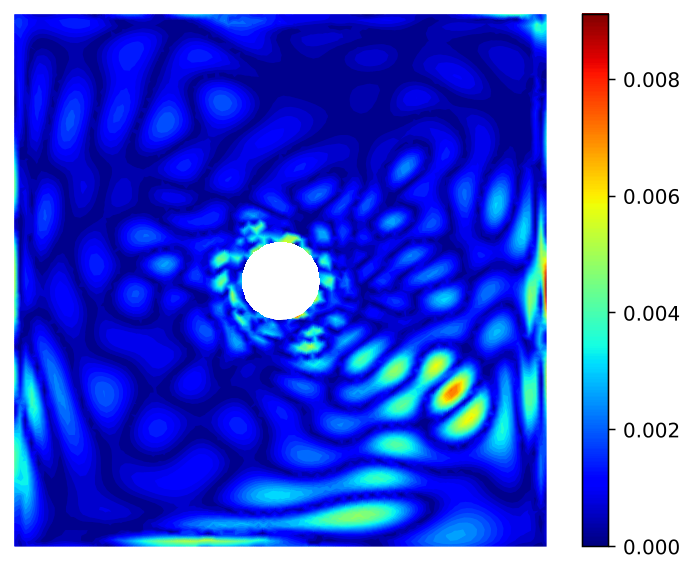}} \hspace{0.47em}
\subfloat{\includegraphics[height = 0.195\linewidth]{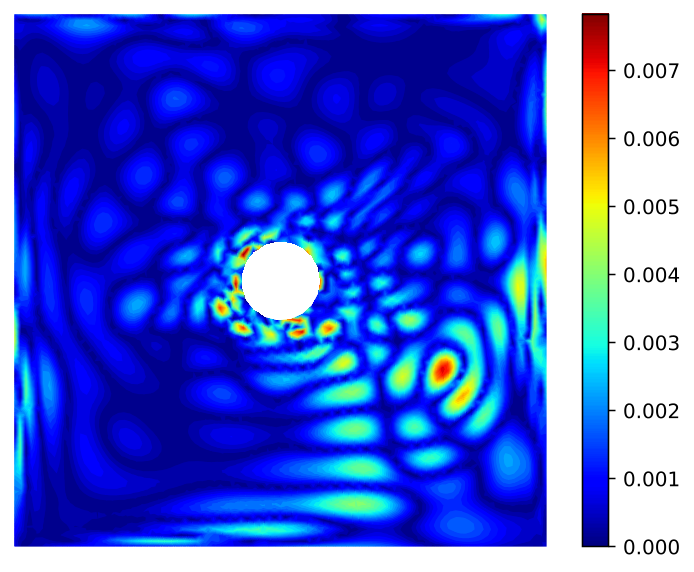}}

\subfloat{\includegraphics[height = 0.195\linewidth]{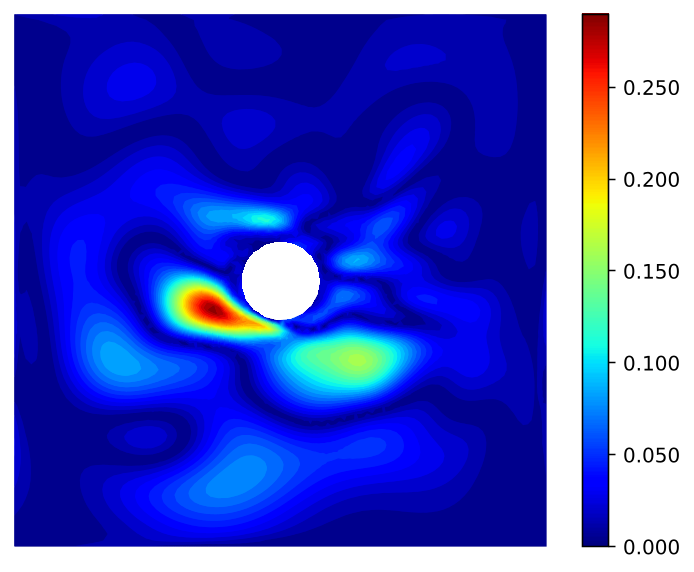}} \hspace{0.47em}
\subfloat{\includegraphics[height = 0.195\linewidth]{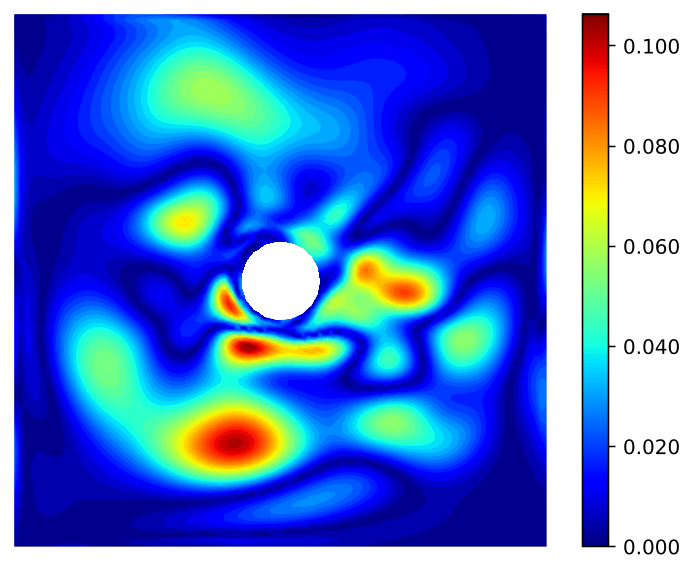}} \hspace{0.47em}
\subfloat{\includegraphics[height = 0.195\linewidth]{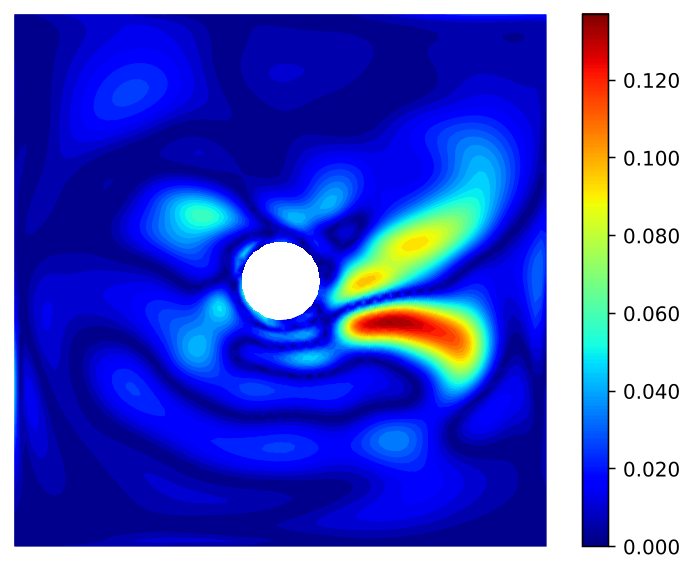}} \hspace{0.47em}
\subfloat{\includegraphics[height = 0.195\linewidth]{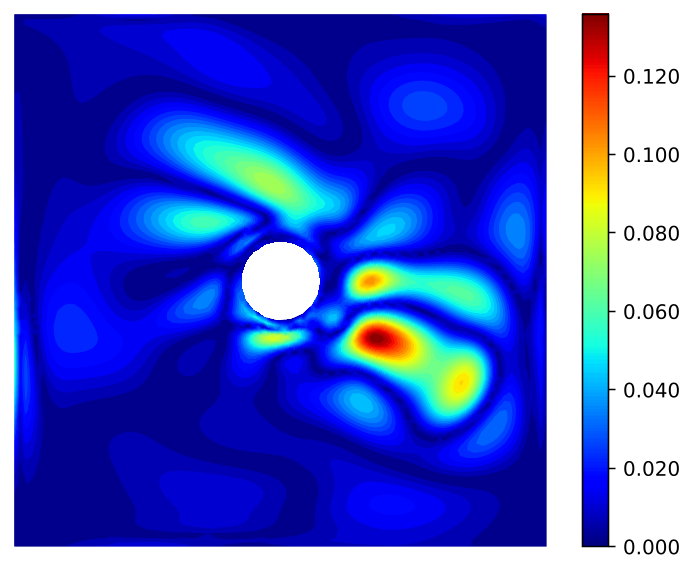}}

\caption{\textit{Test 1.2}. Optimal transport in a fluid. High-fidelity optimal state trajectory (first row), POD+AE reconstructions (second row), forward model predictions (third row), POD+AE reconstruction errors (fourth row) and forward model prediction errors (fifth row) at $t = 0.25, 0.5, 0.75, 1.0$ related to the test scenario parameters $\mus = (-0.30, 0.49, 0.59)$. The underlying fluid velocity vector field $\mathbf{v}_h$ on $\Omega$ is depicted together with the state.}
\label{fig:fluid_state}
\end{figure}

Whenever only a few lagging state data are available online -- in the worst-case scenario, only the initial state is known -- the latent feedback loop can be trained and exploited to continuously control the system anyway. This is possible if we take into account a low-dimensional surrogate model for the time-advancing scheme of the system at the latent level, that is we model $\varphi_N$ with a deep feedforward neural network having $3$ hidden layers of $50$ neurons each and leaky Relu as activation function. While $\varphi_N$ is initialized following the procedure in \cite{He2015}, the autoencoders and the policy inherit the initial weights from the previous training. As detailed in Section~\ref{subsec:latentfeedbackloop}, the cumulative loss function $J_{\mathrm{NN}}$ with $\lambda_1 = \lambda_2 = \lambda_3 = \lambda_6 = 0.0001$ and $\lambda_4 = \lambda_5 = 0.01$ is minimized through the L-BFGS optimization algorithm in $1$ hour and $31$ minutes. The POD+AE reconstruction errors on test data are equal to $\varepsilon^y_{\mathrm{rel}} = 3.71\%$ and $\varepsilon^u_{\mathrm{rel}} = 4.40\%$, while the policy prediction error after-decoding is $7.11\%$. Instead, the prediction-from-data and the prediction-from-policy of the forward model $\varphi_N$ are accurate up to an error equal to, respectively, $2.32 \%$ and $2.48 \%$ at the latent level, while they increase to $7.29\%$ and $7.63 \%$ after POD+AE decoding. The POD+AE state reconstructions and the forward model predictions related to a trajectory in the test set are visualized in Figure~\ref{fig:fluid_state}.

\begin{figure}[!ht]
\centering
\subfloat{\includegraphics[height = 0.195\linewidth]{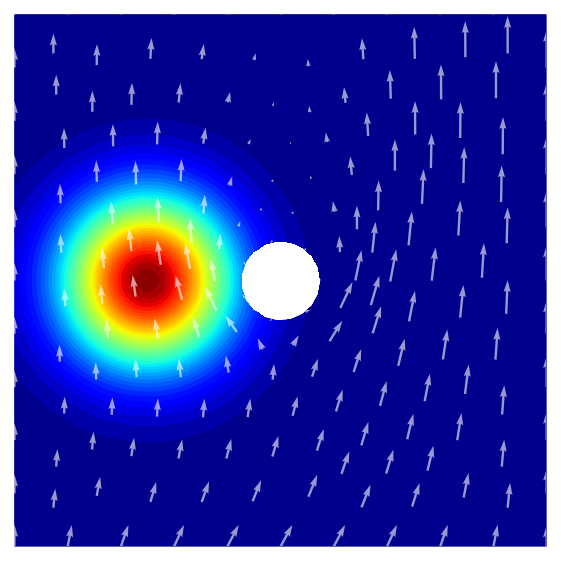}} 
\subfloat{\includegraphics[height = 0.195\linewidth]{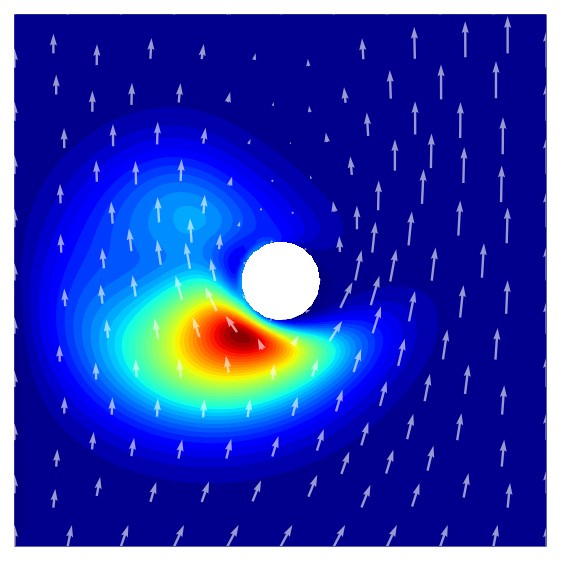}}
\subfloat{\includegraphics[height = 0.195\linewidth]{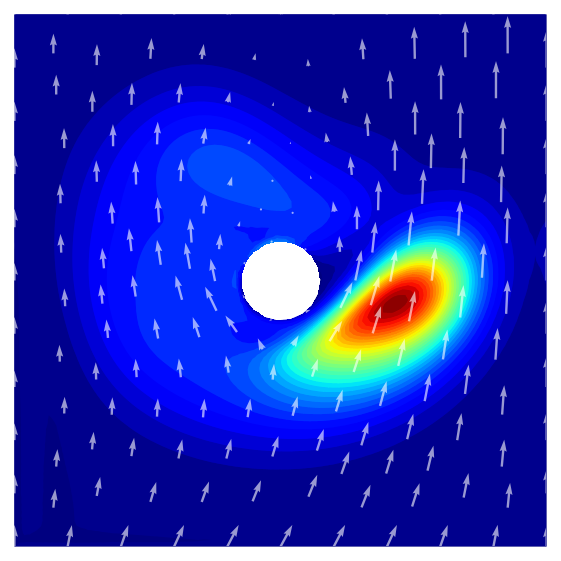}}
\subfloat{\includegraphics[height = 0.195\linewidth]{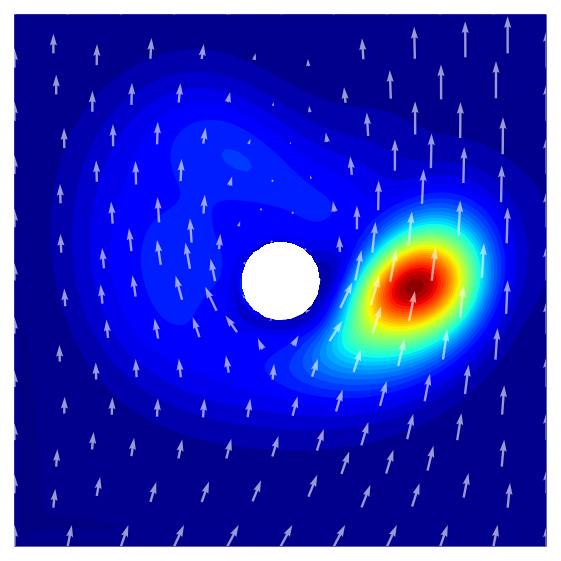}}
\subfloat{\includegraphics[height = 0.195\linewidth]{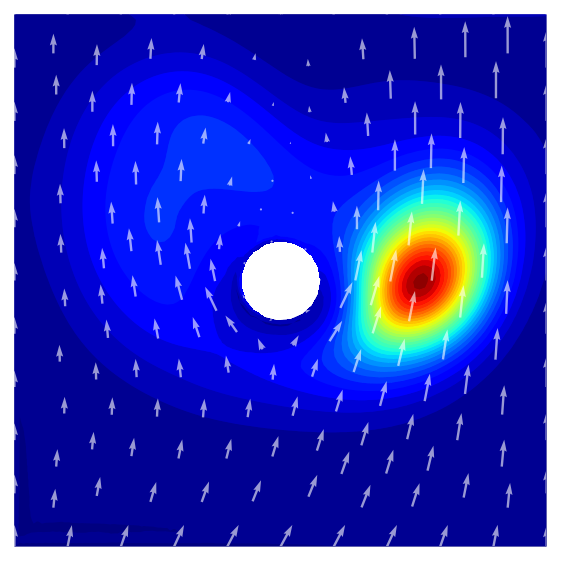}}\vspace{-1em}

\subfloat{\includegraphics[height = 0.195\linewidth]{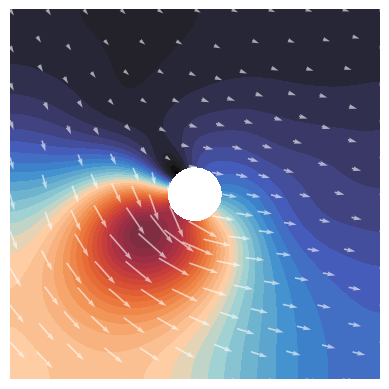}}
\subfloat{\includegraphics[height = 0.195\linewidth]{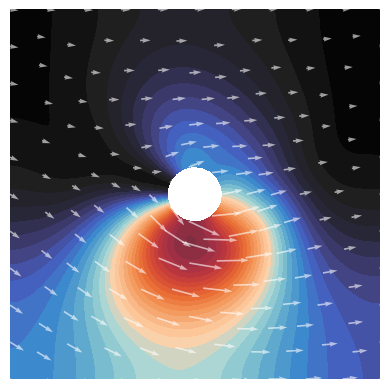}}
\subfloat{\includegraphics[height = 0.195\linewidth]{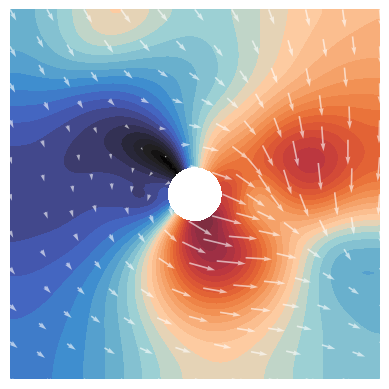}} 
\subfloat{\includegraphics[height = 0.195\linewidth]{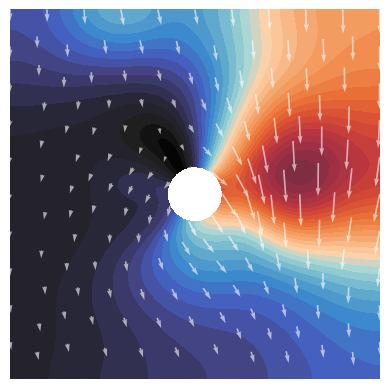}} 
\hphantom{\subfloat{\includegraphics[height = 0.195\linewidth]{Images/Fluid/Control_test1_t=1.png}}}

\subfloat{\includegraphics[height = 0.195\linewidth]{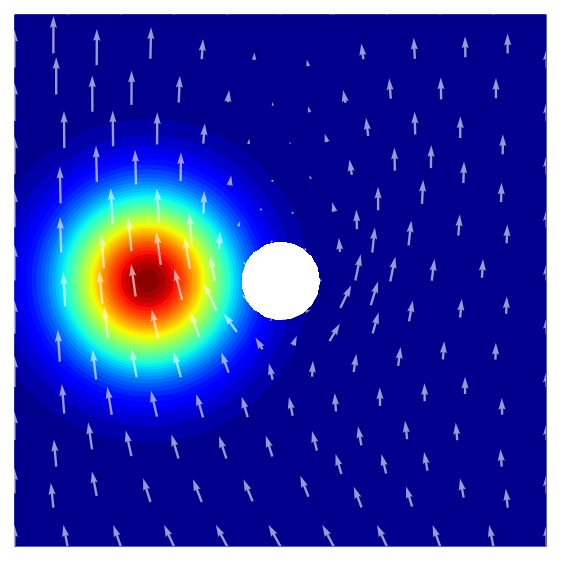}}
\subfloat{\includegraphics[height = 0.195\linewidth]{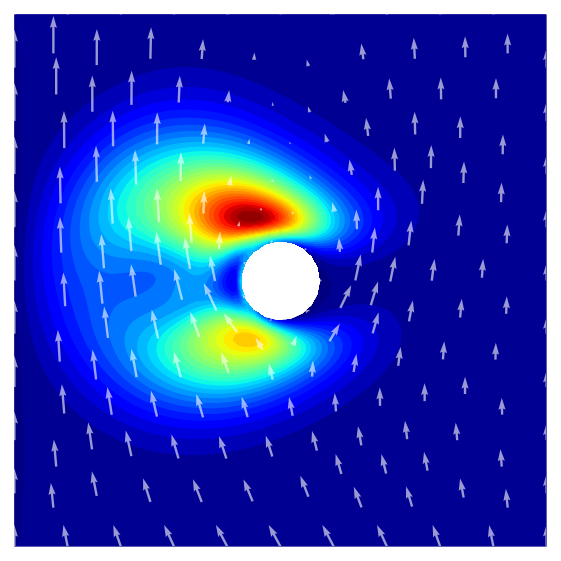}}
\subfloat{\includegraphics[height = 0.195\linewidth]{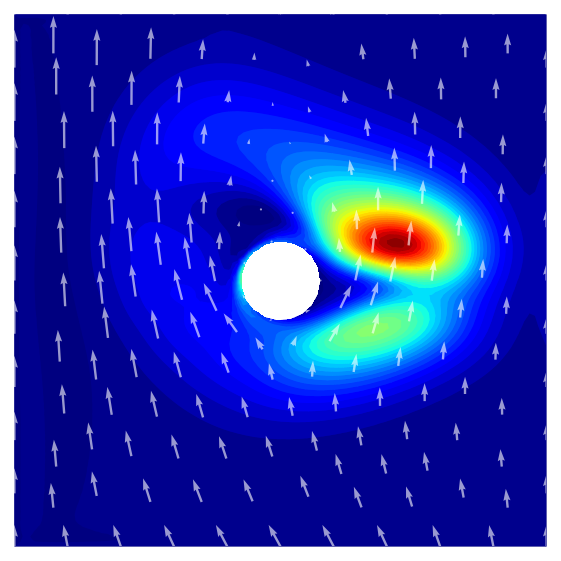}}
\subfloat{\includegraphics[height = 0.195\linewidth]{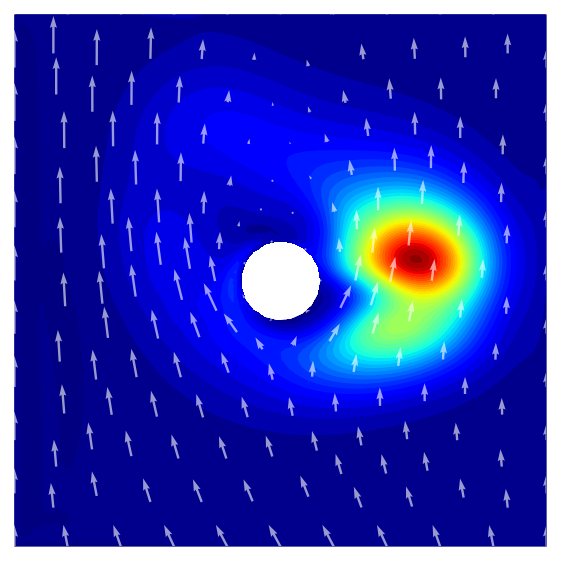}}
\subfloat{\includegraphics[height = 0.195\linewidth]{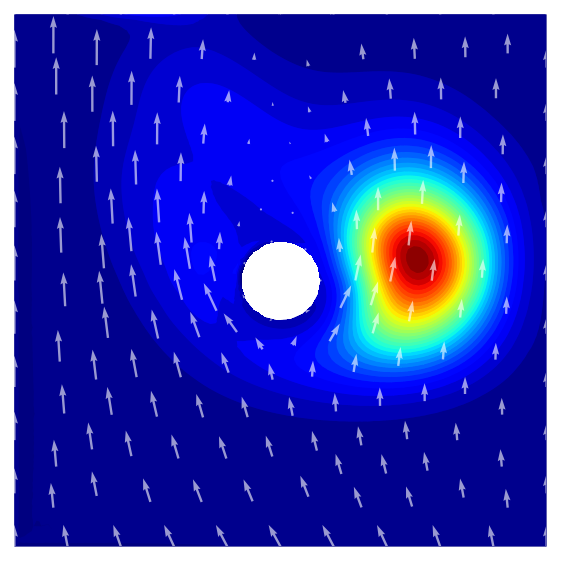}}\vspace{-1em}

\subfloat{\includegraphics[height = 0.195\linewidth]{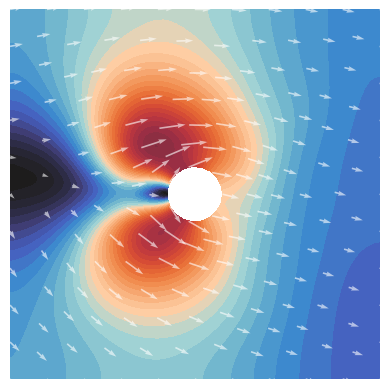}}
\subfloat{\includegraphics[height = 0.195\linewidth]{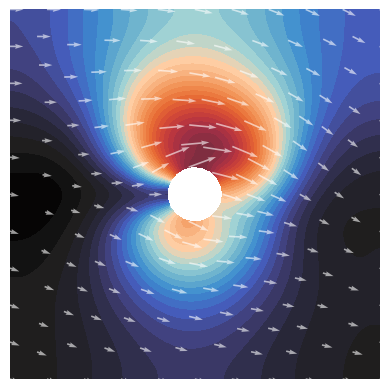}}
\subfloat{\includegraphics[height = 0.195\linewidth]{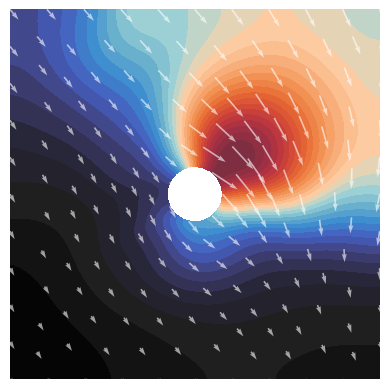}} 
\subfloat{\includegraphics[height = 0.195\linewidth]{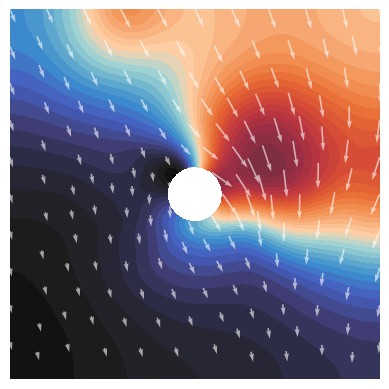}} 
\hphantom{\subfloat{\includegraphics[height = 0.195\linewidth]{Images/Fluid/Control_test2_t=1.png}}}

\caption{\textit{Test 1.2}. Optimal transport in a fluid. First and second rows: system evolution  driven by policy-based controls at $t = 0, 0.25, 0.75, 1.0, 1.5$ related to an initial state centered at $(\mu_1^0, \mu_2^0) = (-0.5, 0.0)$ and a vector of scenario parameters $\mus = (0.0, 0.5, 0.5)$, while exploiting the loop closure at the full-order level. Third and fourth rows: system evolution  driven by policy-based controls at $t = 0, 0.25, 0.75, 1.0, 1.5$ related to an initial state centered at $(\mu_1^0, \mu_2^0) = (-0.5, 0.0)$ and a vector of scenario parameters $\mus = (0.0, 0.5, -0.5)$, while exploiting the latent feedback loop. The underlying fluid velocity vector fields $\mathbf{v}_h$ on $\Omega$ is depicted together with the state. The control velocity fields on $\Omega$ are depicted through vector fields, with the underlying colours corresponding to their magnitude.}
\label{fig:fluid_test}
\end{figure}

Thanks to the trained deep-learning based reduced order feedback controller, we can now move the state for new initial configurations and scenarios, which have not been seen during training. Figure~\ref{fig:fluid_test} presents the state evolutions related to two different test cases with same initial state position $(\mu_1^0, \mu_2^0) = (-0.5, 0.0)$, target destination $(\mu_1^d, \mu_2^d) = (0.5, 0.0)$ and inflow velocity intensity $\gamma_{\text{in}} =0.5$, but different angles of attack ($\alpha_{\text{in}} =0.5$ for the first test case, while $\alpha_{\text{in}} =-0.5$ in the second setting). In the first scenario, we assume that we can continuously monitor the system, that is we have access to high-dimensional state data online. Instead, in the second test case, the latent feedback loop is necessary since only the initial configuration is known. By a visual inspection, it is possible to assess that in both cases the state density is moved towards the target destination. In particular, the reduced order feedback controller is capable of exploiting the underlying fluid flow in order to steer the system with the less expensive and energetic control action, as required by the regularization terms in the loss function $J$. Indeed, when the angle of attack is oriented to the right-hand side of the domain, the state transition occurs below the obstacle, taking advantage of the underlying current to reach the target location. In contrast, when the angle of attack is left-oriented, the controller exploits very different control velocity fields, splitting the state density into two clusters, where the bigger one moves above the obstacle, thus mainly avoiding countercurrent routes that would require more expensive control actions. Moreover, it is possible to note that, once the state reaches the target position, the optimal control is concentrated in the top-right region of $\Omega$ in order to balance the upward thrust provided by $\mathbf{v}_h$. As far as computational times are concerned, the deep learning-based reduced order feedback controller with model closure at the full-order level requires $1.21$ seconds to provide the control actions and simulate the states for all the $N_t$ time steps in the test case considered, while it boils down to $0.05$ seconds in the case of latent feedback loop, with a remarkably high speed-up with respect to full-order methods ($1100 \times$ for the loop closure at the full-order level, $26500 \times$ for the latent feedback loop). As already highlighted in the previous test case, the proposed controllers remain effective even when taking into account different (possibly inhomogeneous) time discretizations than the one exploited offline.

\begin{figure}
\centering
\subfloat{\includegraphics[height = 0.35\linewidth]{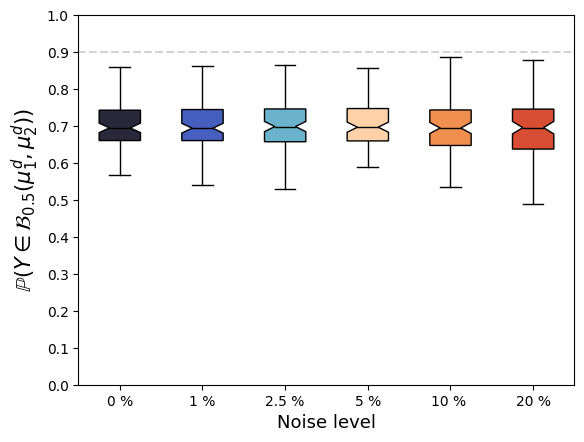}} \quad 
\subfloat{\includegraphics[height = 0.35\linewidth]{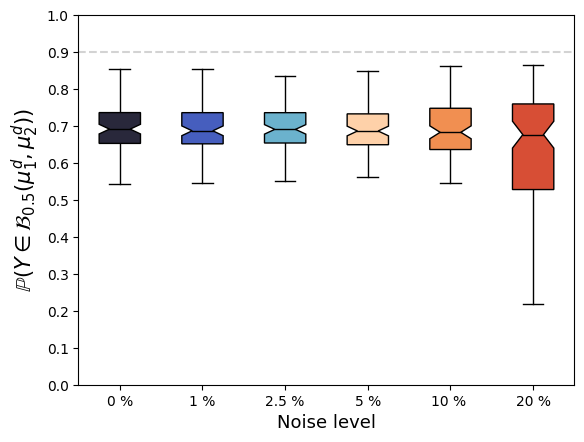}}

\subfloat{\includegraphics[height = 0.195\linewidth]{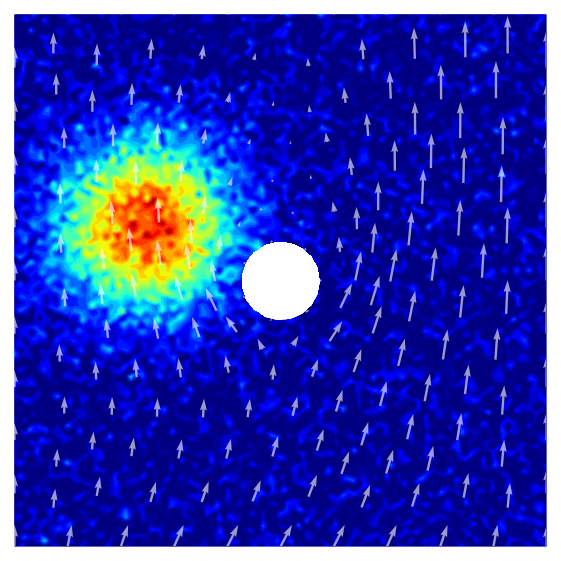}}
\subfloat{\includegraphics[height = 0.195\linewidth]{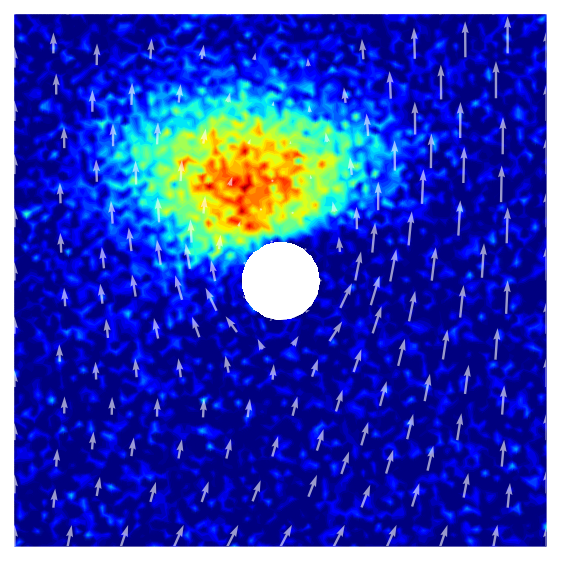}}
\subfloat{\includegraphics[height = 0.195\linewidth]{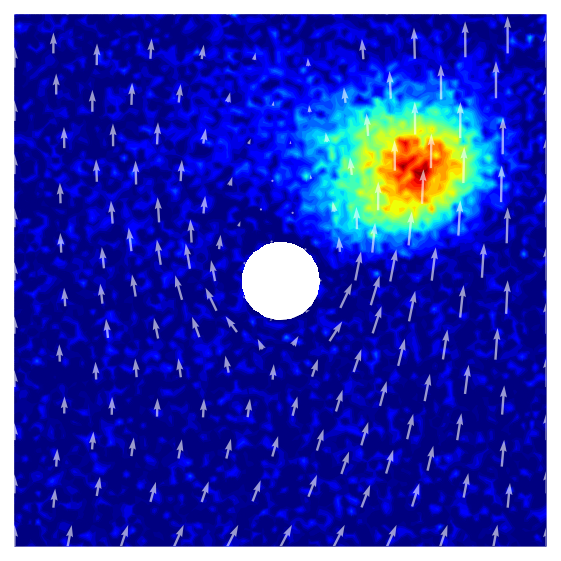}} 
\subfloat{\includegraphics[height = 0.195\linewidth]{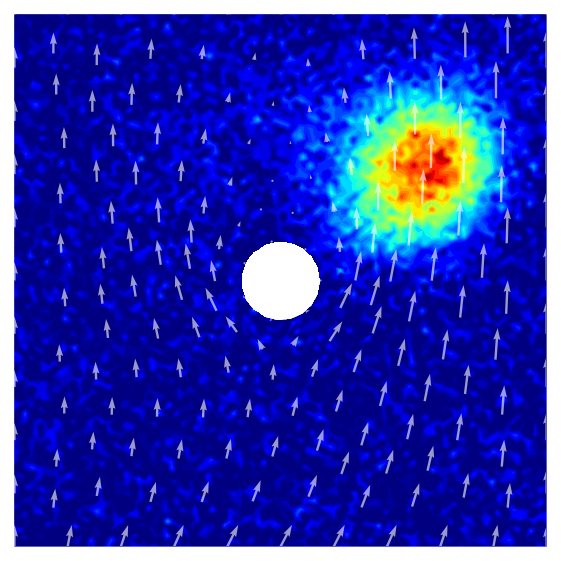}} 
\hphantom{\subfloat{\includegraphics[height = 0.195\linewidth]{Images/Fluid/State_noise_t=1.png}}}
\vspace{-1em}

\subfloat{\includegraphics[height = 0.195\linewidth]{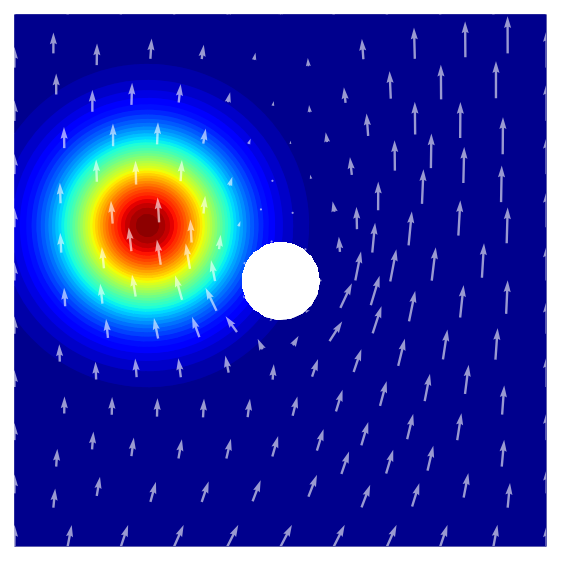}}
\subfloat{\includegraphics[height = 0.195\linewidth]{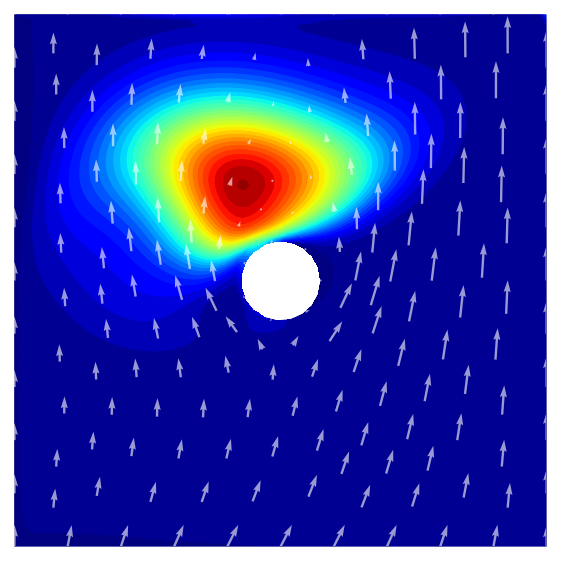}}
\subfloat{\includegraphics[height = 0.195\linewidth]{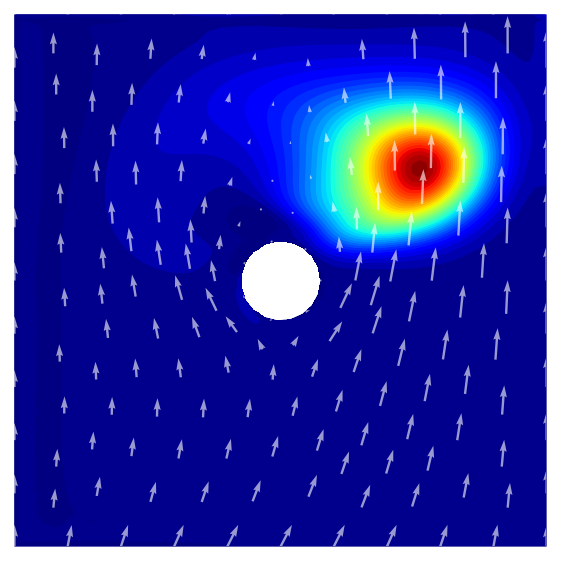}}
\subfloat{\includegraphics[height = 0.195\linewidth]{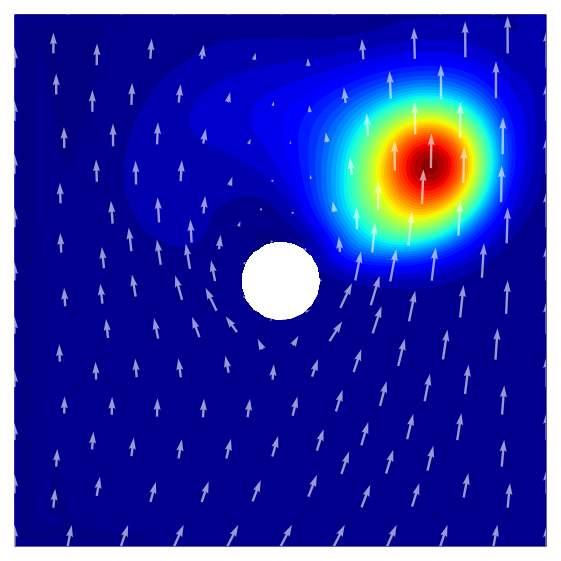}}
\subfloat{\includegraphics[height = 0.195\linewidth]{Images/Fluid/State_noise1_t=1.png}}\vspace{-1em}

\subfloat{\includegraphics[height = 0.195\linewidth]{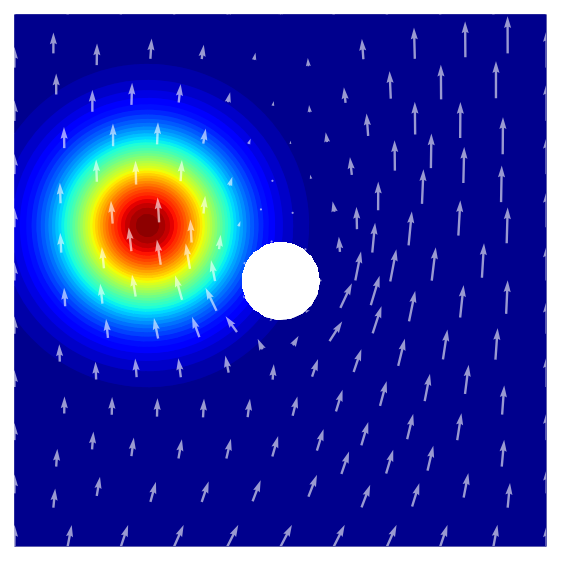}}
\subfloat{\includegraphics[height = 0.195\linewidth]{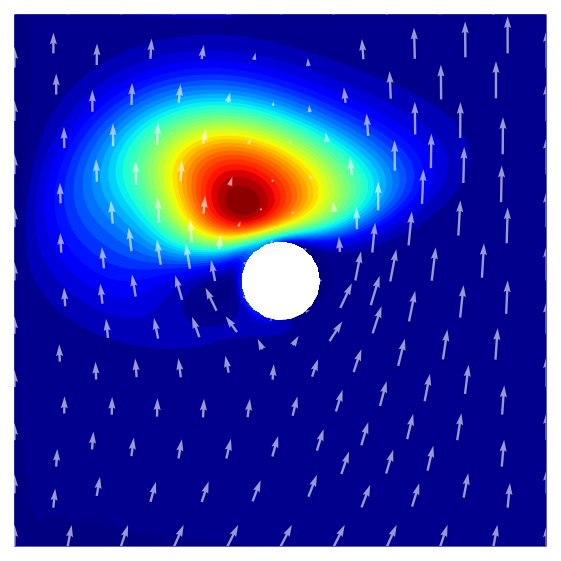}}
\subfloat{\includegraphics[height = 0.195\linewidth]{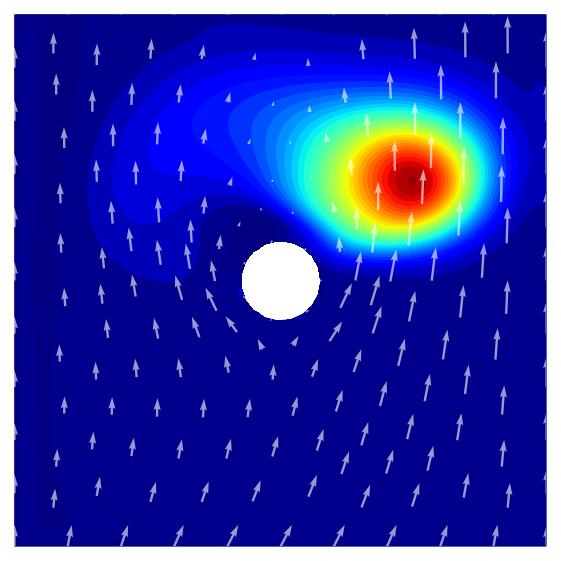}}
\subfloat{\includegraphics[height = 0.195\linewidth]{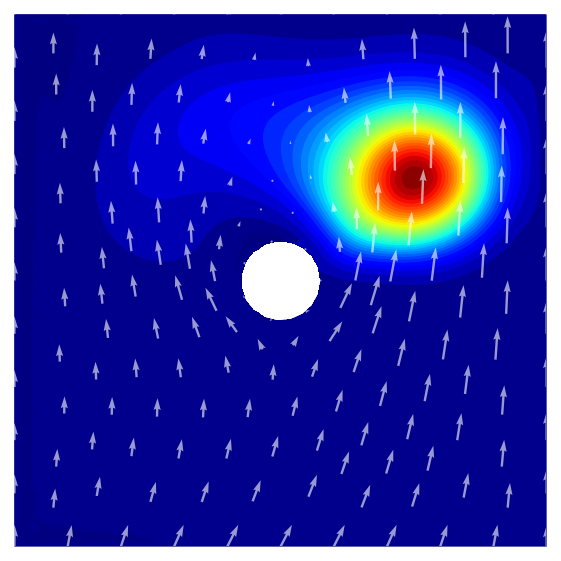}}
\subfloat{\includegraphics[height = 0.195\linewidth]{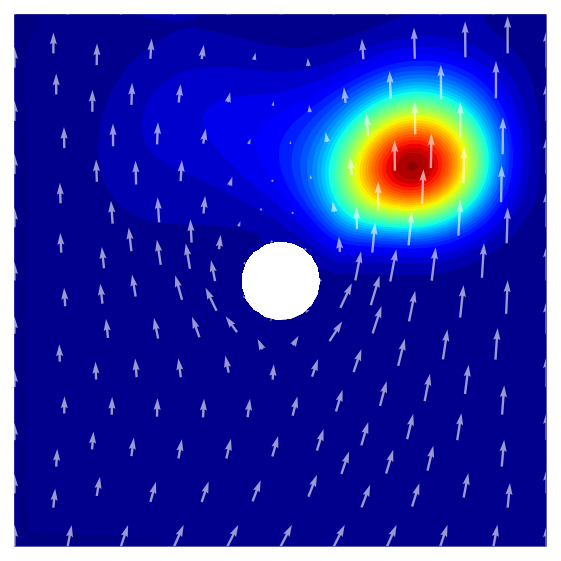}}

\caption{\textit{Test 1.2}. Optimal transport in a fluid. First row: boxplots of the probabilities of arrival in $100$ random scenarios for different noise levels when considering the deep learning-based reduced order feedback controller at full-order level (left) and the latent feedback loop (right). Second row: online state data at $t = 0, 0.25, 0.75, 1.0$ corrupted by Gaussian random noise with standard deviation $\sigma = 0.3$ related to an initial state centered at $(\mu_1^0, \mu_2^0) = (-0.5, 0.2)$ and a vector of scenario parameters $\mus = (0.42, 0.39, 0.48)$. Other panels: system evolution at $t = 0, 0.25, 0.75, 1.0, 1.5$ driven by policy-based controls exploiting state data corrupted by Gaussian random noise with standard deviation $\sigma = 0.3$ related to an initial state centered at $(\mu_1^0, \mu_2^0) = (-0.5, 0.2)$ and a vector of scenario parameters $\mus = (0.42, 0.39, 0.48)$, while exploiting the loop closure at the full-order (third row) and latent (fourth row) level. The underlying fluid velocity vector field $\mathbf{v}_h$ on $\Omega$ is depicted together with the state.}
\label{fig:fluid_noise}
\end{figure}

In this work, state information in the online phase is generated through full-order, high-fidelity simulations of the underlying dynamics, which is assumed to be known. However, in several practical settings, the state snapshots are captured through sensors, which may be affected by noise. To assess the generalization capabilities of our controller against disturbances, we try to optimally transport the state density in a parametric fluid while considering online state information corrupted by noise. We recall that noisy data are exploited here only in the online evaluation phase, while the training is entirely performed with the aforementioned noise-free snapshots. In particular, the noise is modeled through a Gaussian distribution, that is $\boldsymbol{\varepsilon} \sim \mathcal{N}(\boldsymbol{0}, \sigma^2 I)$, being $\boldsymbol{0}$ the $N_h^y$-dimensional zero vector and $I$ the $N_h^y \times N_h^y$ identity matrix, where $5$ different noise levels are taken into account, namely the standard deviations $\sigma = 0.03, 0.075, 0.15, 0.3, 0.6$ (which correspond approximately to, respectively, $1 \%, 2.5 \%, 5 \%, 10 \%, 20 \%$ of the range of state values). Figure~\ref{fig:fluid_noise} compares the performances of the deep learning-based reduced order feedback controllers in $100$ random scenarios across different noise levels, along with the noise-free case. In particular, for every test case, we compute the probability of arrival 
\[
\mathbb{P}(Y \in \mathcal{B}_{0.5}(\mu_1^d, \mu_2^d)) \approx \int_{\mathcal{B}_{0.5}(\mu_1^d, \mu_2^d)} y(T) d\Omega
\]
where $Y$ is a random variable with probability density function equal to $y(T)$, that is the final state obtained by exploiting our controllers, extended on the whole $\R^2$. Instead, the circle centered at $(\mu_1^d, \mu_2^d)$ with radius $0.5$, namely $\mathcal{B}_{0.5}(\mu_1^d, \mu_2^d)$, corresponds approximately to the $90 \%$-confidence region of the target density, that is $\mathbb{P}(Y_d \in \mathcal{B}_{0.5}(\mu_1^d, \mu_2^d)) \approx 0.9$, where $Y_d$ is a random variable having $y_d$ as probability density function. Despite the noise added to the state data, it is possible to assess that both the feedback loop at the full-order and latent level are capable of steering the state density towards the final target, with results comparable to the setting without noise. This is mainly due to the POD+AE reduction, which is able to extract the relevant features for control design, while discarding erratic and high-frequency disturbances. When the noise level is remarkably high, the accuracy of the feedback latent loop is lower since its optimal transport predictions are entirely based on a single noisy state snapshot at $t = 0$. Instead, the full-order counterpart is still able to steer the state towards the final destination taking advantage of the multiple corrupted state data received at every time step. Note that our reduced order feedback controllers perform all the $600$ simulations considered in this analysis in $720$ and $30$ seconds, respectively, while a full-order solver based on, e.g., FEM would require approximately $8$ days of computations. An example of a state trajectory driven by our controllers when taking into account a standard deviation $\sigma$ equal to $0.3$, a random initial setting and a random scenario is available in Figure~\ref{fig:fluid_noise}. 

\section{Conclusions}\label{sec:conclusions}

In this work, we present a deep learning-based reduced order feedback controller capable of steering dynamical systems very rapidly. Unlike several control strategies available in the literature, the proposed framework is capable of dealing with both {\em (i)} complex and parametrized dynamics modeled via (possibly nonlinear) time-dependent PDEs, {\em (ii)} high-dimensional state observations and {\em (iii)} distributed control actions. This allows us to rapidly control complicated systems in multiple scenarios unseen during training, as often required in applications, amortizing the offline cost due to data generation and networks training. Although only full-order synthetic data have been exploited in this work, the proposed framework can be easily extended to (even low- and high-dimensional) sensor data or videos capturing the dynamics, paving the way for cheap and portable control devices. Indeed, as shown in Section~\ref{subsec:fluid} when dealing with the optimal transport problem in a fluid, the proposed architectures are capable of dealing with noisy data, thanks to the reduction carried out and the feedback signal considered. To handle the high-dimensionality of the data, we extract low-dimensional features relevant for control design through very accurate and efficient non-intrusive reduced order models, such as POD, AE and POD+AE. Note that we consider a unified framework with different reduction strategies in order to exploit the most effective one for each problem at hand, extending the current state-of-the-art on non-intrusive ROMs to closed-loop control problems. Thanks to data compression, it is now possible to learn a low-dimensional surrogate model for the policy bridging state and control latent spaces, which is faster to train and evaluate online.

Inspired by the MPC approach, we consider a model closure at the latent level, here referred to as latent feedback loop, embedding a surrogate model of the reduced order dynamics in the controller. Our controller is therefore able to continuously control the dynamics of interest even in the absence of online state data, avoiding losses of optimality, performance or stability, while overcoming the necessity of continuous monitoring of the dynamical systems. 

Throughout the optimal transport test cases presented, we demonstrate the accuracy of our approach in very challenging settings, characterized by high-dimensional variables, transport-dominated trajectories and complex parameters-to-solution dependencies. The speed-up of our controllers with respect to full-order high-fidelity models based on, e.g., FEM is remarkably high. Indeed, the proposed control strategies consist of efficient forward passes through the considered networks, which usually exploit light architectures due to the dimensionality reduction performed. Moreover, after training, the same architecture may be recycled to obtain rapidly different control actions related to different scenarios of interest.

The proposed controller may be extended in future works in multiple directions. For instance, as already mentioned, different data sources may be considered, such as sensors or cameras recording the system evolution, possibly corrupted by noise. In addition, multi-agent reinforcement learning strategies may be exploited, overcoming the necessity of a high-fidelity full-order OCP solver to generate training data, while still considering distributed parametrized systems. Another possible improvement may be dedicated to speed up the offline phase. For instance, both data-driven or physics-informed surrogate models for the state and adjoint equations may be introduced in order to rapidly augment the dataset of (possibly few) high-fidelity optimal snapshots with many low-fidelity data.

%\section*{Supporting information}

\section*{Acknowledgments}
AM acknowledges the Project “Reduced Order Modeling and Deep Learning for the real- time approximation of PDEs (DREAM)” (Starting Grant No. FIS00003154), funded by the Italian Science Fund (FIS) - Ministero dell'Università e della Ricerca and the project FAIR (Future Artificial Intelligence Research), funded by the NextGenerationEU program within the PNRR-PE-AI scheme (M4C2, Investment 1.3, Line on Artificial Intelligence). MT and AM are members of the Gruppo Nazionale Calcolo Scientifico-Istituto Nazionale di Alta Matematica (GNCS-INdAM).

\bibliographystyle{abbrv}
\bibliography{references}  %%% Uncomment this line and comment out the ``thebibliography'' section below to use the external .bib file (using bibtex) .

\end{document}